\documentclass[10pt, a4paper]{article}
\usepackage{titling}
\title{Finite Element Approximation of \protect\\the Fractional Porous Medium Equation}
\newcommand{\titlePDF}{Finite Element Scheme for the Fractional Porous Medium Equation
with Fractional Pressure}
\newcommand{\authorPDF}{Carrillo, Fronzoni, Suli.}
\newcommand{\subjectPDF}{}
\newcommand{\keywordsPDF}{}

 \usepackage{authblk}
\usepackage{verbatim}

\author[1]{Jos\'{e} A. Carrillo}
\author[1]{\authorcr Stefano Fronzoni}
\author[1]{Endre S\"{u}li}

\affil[1]{
	Mathematical Institute, University of Oxford
}
\affil[ ]{
	OX2 6GG Oxford, United Kingdom
}
\affil[ ]{\textit{carrillo@maths.ox.ac.uk, fronzoni@maths.ox.ac.uk, suli@maths.ox.ac.uk}}

\affil[ ]{}
 
\usepackage[left=0.7in, right=0.7in, top=1.1in, bottom=1.1in]{geometry}

\makeatletter
\let\newtitle\@title
\let\newauthor\@author
\let\newdate\@date
\makeatother
\usepackage{fancyhdr}
\usepackage[
	hypertexnames=false,
	pdftex,
	pdftitle={\titlePDF. \authorPDF},
	pdfauthor={\authorPDF},
	pdfsubject={\subjectPDF},
	pdfkeywords={\keywordsPDF},
]{hyperref}
\hypersetup{
	colorlinks=true,
	linkcolor=blue,
	citecolor=blue,
	urlcolor=blue,
	linktocpage
}

\usepackage[utf8]{inputenc}
\usepackage[T1]{fontenc}
\usepackage{amsmath} 
\usepackage{amsthm}
\usepackage{amsfonts}
\usepackage{amssymb}
\usepackage{graphicx}
\usepackage{mathtools}
\usepackage{enumitem}
\usepackage{ifthen}
\usepackage{dsfont}
\usepackage{float}
\usepackage{lettrine}
\usepackage{lipsum}
\usepackage{subcaption}

\usepackage[UKenglish]{babel}

\usepackage[dvipsnames]{xcolor}

\definecolor{ppGreen}{HTML}{008000}
\definecolor{ppBlue}{HTML}{0000FF}
\definecolor{ppRed}{HTML}{FF0000}
\definecolor{ppPurple}{HTML}{800080}
\definecolor{lightblue}{rgb}{0.145,0.6666,1}
\definecolor{grey52}{RGB}{52,52,52}
\definecolor{color1}{RGB}{0,62,116}
\definecolor{color2}{RGB}{152,152,152}
\definecolor{color3}{RGB}{52,52,52}
\definecolor{color4}{RGB}{100,100,100}

\definecolor{imperialnavy}{RGB}{0,33,71}
\definecolor{imperialblue}{RGB}{0,62,116}
\definecolor{imperialgrey}{RGB}{235,238,238}
\definecolor{imperialcoolgrey}{RGB}{157,157,157}

\newcounter{review}

\newcommand{\ntcreview}[3]{\refstepcounter{review}{\color{#2}{\textbf{[#1]}: #3}}}

\newcommand{\creview}[3]{\ntcreview{#1}{#2}{#3}\addcontentsline{tor}{subsection}{\thereview~\textbf{[#1]}:~#3
	}}

\newcommand{\review}[2]{\creview{#1}{blue}{#2}}

\makeatletter

\newcommand\listreviewname{List of Reviews}
\newcommand\listofreviews{\section*{\listreviewname}\@starttoc{tor}}
\makeatother

\usepackage[all]{nowidow}

\newcommand{\subjectclassification}[1]{

	{\small\textbf{\textit{AMS Subject Classification --- }} #1}

}
\newcommand{\keywords}[1]{

	{\small\textbf{\textit{Keywords --- }} #1}

}

\usepackage{tikz}

\usepackage{pgfplots}
\pgfplotsset{compat=1.16}

\usepackage{float}
\usepackage{subcaption}
\usepackage[section]{placeins}
\usepackage{rotating}

\usepackage{array}
\newcolumntype{L}[1]{>{\raggedright\let\newline\\\arraybackslash\hspace{0pt}}m{#1}}
\newcolumntype{C}[1]{>{\centering\let\newline\\\arraybackslash\hspace{0pt}}m{#1}}
\newcolumntype{R}[1]{>{\raggedleft\let\newline\\\arraybackslash\hspace{0pt}}m{#1}}
\usepackage{booktabs}
\usepackage{tabularx}
\usepackage{multirow}

\usepackage{titling}

\usepackage[algoruled]{algorithm2e}

\usepackage[capitalise]{cleveref}

 \usepackage{accents}

\newcommand\term\emph

\numberwithin{equation}{section}

\usepackage{autobreak}

\usepackage{enumerate}

\makeatletter
\def\@maketitle{\newpage
	\begin{center}\let \footnote \thanks
		{\LARGE\bfseries \@title \par}\vskip 2.5em{\large
				\lineskip .5em\begin{tabular}[t]{c}\@author
				\end{tabular}\par}\vskip 1em{\large \@date}\end{center}\par
	\vskip 1.5em}
\makeatother



\usepackage{amsthm}
\theoremstyle{plain}
\newtheorem{theo}{Theorem}[section]
\newtheorem{lemma}[theo]{Lemma}
\newtheorem{prop}[theo]{Proposition}
\newtheorem{cor}[theo]{Corollary}

\theoremstyle{remark}
\newtheorem{remark}[theo]{\bf Remark}

\usepackage{esint}

\def\XXint#1#2#3{{\setbox0=\hbox{$#1{#2#3}{\int}$ }
			\vcenter{\hbox{$#2#3$ }}\kern-.6\wd0}}

\DeclarePairedDelimiter{\norm}{\|}{\|}

\DeclarePairedDelimiter{\inn}{\langle}{\rangle}

\usepackage{suffix}

\newcommand{\inner}[2]{\inn{#1,#2}}
\WithSuffix\newcommand\inner*[2]{\inn*{#1,#2}}

\DeclarePairedDelimiter{\positive}{(}{)^{+}}
\DeclarePairedDelimiter{\negative}{(}{)^{-}}

\newcommand\pos\positive
\renewcommand\neg\negative

\WithSuffix\newcommand\pos*{\positive*}
\WithSuffix\newcommand\neg*{\negative*}

\newcommand{\R}{{\mathbb{R}}}

\newcommand{\pnorm}[2]{\norm{#2}_{\L{#1}}}
\WithSuffix\newcommand\pnorm*[2]{\norm*{#2}_{\L{#1}}}

\newcommand{\psnorm}[3]{\norm{#3}_{\L{#1}(#2)}}
\WithSuffix\newcommand\psnorm*[3]{\norm*{#3}_{\L{#1}(#2)}}

\newcommand{\pnormp}[2]{\pnorm{#1}{#2}^{#1}}
\WithSuffix\newcommand\pnormp*[2]{\pnorm*{#1}{#2}^{#1}}

\newcommand{\psnormp}[3]{\psnorm{#1}{#2}{#3}^{#1}}
\WithSuffix\newcommand\psnormp*[3]{\psnorm*{#1}{#2}{#3}^{#1}}

\newcommand\svec\vec

\renewcommand{\vec}{\mathbf}
\renewcommand{\svec}{\boldsymbol}

\newcommand{\conv}{\ast}
\renewcommand{\star}{\conv}

\renewcommand{\d}{\mathrm{d}}
\newcommand{\dd}{\mathop{}\!\d}

\newcommand{\dt}{\dd t}

\newcommand{\dx}{\dd x}
\newcommand{\dy}{\dd y}

\newcommand{\grad}{\nabla}

\newcommand{\ppr}{(r)}

\newcommand{\Wr}{^{W,\,\ppr}}

\newlength{\dhatheight}

\ifthenelse{\isundefined{\Wr}}{
	\newcommand{\Wr}{^{W,\,\ppr}}
}{
	\renewcommand{\Wr}{^{W,\,\ppr}}
}

\newcommand{\s}{s}

\newcommand{\Hs}{\mathbb{H}^{\s}}
\newcommand{\Hsast}{\mathbb{H}^{\s}_{\ast}}
\newcommand{\rhodL}{\rho_{\delta, L}}
\newcommand{\cdL}{c_{\delta,L}}
\newcommand{\rhodLh}{\rho_{h, \delta, L}}
\newcommand{\cdLh}{c_{h, \delta, L}}
\newcommand{\rhoL}{\rho_{L}}
\newcommand{\rhoLDt}{\rho^{\Delta t}_{L}}
\newcommand{\rhoLDtp}{\rho^{\Delta t, +}_{L}}
\newcommand{\rhoLDtm}{\rho^{\Delta t, -}_{L}}
\newcommand{\rhoLDtpm}{\rho^{\Delta t(, \pm)}_{L}}
\newcommand{\rhoDt}{\rho^{\Delta t}}
\newcommand{\rhoDtp}{\rho^{\Delta t, +}}
\newcommand{\rhoDtm}{\rho^{\Delta t, -}}
\newcommand{\rhoDtpm}{\rho^{\Delta t(, \pm)}}
\newcommand{\cL}{c_{L}}

\newcommand{\dtau}{\mathrm{d}\tau}

\newcommand*\circled[1]{\tikz[baseline=(char.base)]{
            \node[shape=circle,draw,inner sep=0.5pt] (char) {#1};}}

 \usepackage{todonotes}

\newif\ifskiptable

\pgfplotsset{colormap={hsv}{
			hsb(0.00cm)=(0.00,0,0.95);
			hsb(0.05cm)=(0.05,1,1);
			hsb(0.10cm)=(0.10,1,1);
			hsb(0.15cm)=(0.15,1,1);
			hsb(0.20cm)=(0.20,1,1);
			hsb(0.25cm)=(0.25,1,1);
			hsb(0.30cm)=(0.30,1,1);
			hsb(0.35cm)=(0.35,1,1);
			hsb(0.40cm)=(0.40,1,1);
			hsb(0.45cm)=(0.45,1,1);
			hsb(0.50cm)=(0.50,1,1);
			hsb(0.55cm)=(0.55,1,1);
			hsb(0.60cm)=(0.60,1,1);
			hsb(0.65cm)=(0.65,1,1);
			hsb(0.70cm)=(0.70,1,1);
			hsb(0.75cm)=(0.75,1,1);
			hsb(0.80cm)=(0.80,1,1);
			hsb(0.85cm)=(0.85,1,1);
			hsb(0.90cm)=(0.90,1,1);
			hsb(0.95cm)=(0.95,1,1);
			hsb(1.00cm)=(1.00,1,1);
		}
}

\pgfplotsset{colormap={hsvSoft}{
			hsb(0.00cm)=(0.00,0,0.95);
			hsb(0.05cm)=(0.05,1,1);
			hsb(0.10cm)=(0.10,1,1);
			hsb(0.15cm)=(0.15,1,1);
			hsb(0.20cm)=(0.20,1,1);
			hsb(0.25cm)=(0.25,1,1);
			hsb(0.30cm)=(0.30,1,1);
			hsb(0.35cm)=(0.35,1,1);
			hsb(0.40cm)=(0.40,1,1);
			hsb(0.45cm)=(0.45,1,1);
			hsb(0.50cm)=(0.50,1,1);
			hsb(0.55cm)=(0.55,1,1);
			hsb(0.60cm)=(0.60,1,1);
			hsb(0.65cm)=(0.65,1,1);
			hsb(0.70cm)=(0.70,1,1);
			hsb(0.75cm)=(0.75,1,1);
			hsb(0.80cm)=(0.80,1,1);
			hsb(0.85cm)=(0.85,1,1);
			hsb(0.90cm)=(0.90,1,1);
			hsb(0.95cm)=(0.95,1,1);
			hsb(1.00cm)=(0.00,0,0.95);
		}
}

\pgfplotsset{colormap={viridisFull}{
			rgb=(0.26700401, 0.00487433, 0.32941519)
			rgb=(0.26851048, 0.00960483, 0.33542652)
			rgb=(0.26994384, 0.01462494, 0.34137895)
			rgb=(0.27130489, 0.01994186, 0.34726862)
			rgb=(0.27259384, 0.02556309, 0.35309303)
			rgb=(0.27380934, 0.03149748, 0.35885256)
			rgb=(0.27495242, 0.03775181, 0.36454323)
			rgb=(0.27602238, 0.04416723, 0.37016418)
			rgb=(0.2770184 , 0.05034437, 0.37571452)
			rgb=(0.27794143, 0.05632444, 0.38119074)
			rgb=(0.27879067, 0.06214536, 0.38659204)
			rgb=(0.2795655 , 0.06783587, 0.39191723)
			rgb=(0.28026658, 0.07341724, 0.39716349)
			rgb=(0.28089358, 0.07890703, 0.40232944)
			rgb=(0.28144581, 0.0843197 , 0.40741404)
			rgb=(0.28192358, 0.08966622, 0.41241521)
			rgb=(0.28232739, 0.09495545, 0.41733086)
			rgb=(0.28265633, 0.10019576, 0.42216032)
			rgb=(0.28291049, 0.10539345, 0.42690202)
			rgb=(0.28309095, 0.11055307, 0.43155375)
			rgb=(0.28319704, 0.11567966, 0.43611482)
			rgb=(0.28322882, 0.12077701, 0.44058404)
			rgb=(0.28318684, 0.12584799, 0.44496 )
			rgb=(0.283072 , 0.13089477, 0.44924127)
			rgb=(0.28288389, 0.13592005, 0.45342734)
			rgb=(0.28262297, 0.14092556, 0.45751726)
			rgb=(0.28229037, 0.14591233, 0.46150995)
			rgb=(0.28188676, 0.15088147, 0.46540474)
			rgb=(0.28141228, 0.15583425, 0.46920128)
			rgb=(0.28086773, 0.16077132, 0.47289909)
			rgb=(0.28025468, 0.16569272, 0.47649762)
			rgb=(0.27957399, 0.17059884, 0.47999675)
			rgb=(0.27882618, 0.1754902 , 0.48339654)
			rgb=(0.27801236, 0.18036684, 0.48669702)
			rgb=(0.27713437, 0.18522836, 0.48989831)
			rgb=(0.27619376, 0.19007447, 0.49300074)
			rgb=(0.27519116, 0.1949054 , 0.49600488)
			rgb=(0.27412802, 0.19972086, 0.49891131)
			rgb=(0.27300596, 0.20452049, 0.50172076)
			rgb=(0.27182812, 0.20930306, 0.50443413)
			rgb=(0.27059473, 0.21406899, 0.50705243)
			rgb=(0.26930756, 0.21881782, 0.50957678)
			rgb=(0.26796846, 0.22354911, 0.5120084 )
			rgb=(0.26657984, 0.2282621 , 0.5143487 )
			rgb=(0.2651445 , 0.23295593, 0.5165993 )
			rgb=(0.2636632 , 0.23763078, 0.51876163)
			rgb=(0.26213801, 0.24228619, 0.52083736)
			rgb=(0.26057103, 0.2469217 , 0.52282822)
			rgb=(0.25896451, 0.25153685, 0.52473609)
			rgb=(0.25732244, 0.2561304 , 0.52656332)
			rgb=(0.25564519, 0.26070284, 0.52831152)
			rgb=(0.25393498, 0.26525384, 0.52998273)
			rgb=(0.25219404, 0.26978306, 0.53157905)
			rgb=(0.25042462, 0.27429024, 0.53310261)
			rgb=(0.24862899, 0.27877509, 0.53455561)
			rgb=(0.2468114 , 0.28323662, 0.53594093)
			rgb=(0.24497208, 0.28767547, 0.53726018)
			rgb=(0.24311324, 0.29209154, 0.53851561)
			rgb=(0.24123708, 0.29648471, 0.53970946)
			rgb=(0.23934575, 0.30085494, 0.54084398)
			rgb=(0.23744138, 0.30520222, 0.5419214 )
			rgb=(0.23552606, 0.30952657, 0.54294396)
			rgb=(0.23360277, 0.31382773, 0.54391424)
			rgb=(0.2316735 , 0.3181058 , 0.54483444)
			rgb=(0.22973926, 0.32236127, 0.54570633)
			rgb=(0.22780192, 0.32659432, 0.546532 )
			rgb=(0.2258633 , 0.33080515, 0.54731353)
			rgb=(0.22392515, 0.334994 , 0.54805291)
			rgb=(0.22198915, 0.33916114, 0.54875211)
			rgb=(0.22005691, 0.34330688, 0.54941304)
			rgb=(0.21812995, 0.34743154, 0.55003755)
			rgb=(0.21620971, 0.35153548, 0.55062743)
			rgb=(0.21429757, 0.35561907, 0.5511844 )
			rgb=(0.21239477, 0.35968273, 0.55171011)
			rgb=(0.2105031 , 0.36372671, 0.55220646)
			rgb=(0.20862342, 0.36775151, 0.55267486)
			rgb=(0.20675628, 0.37175775, 0.55311653)
			rgb=(0.20490257, 0.37574589, 0.55353282)
			rgb=(0.20306309, 0.37971644, 0.55392505)
			rgb=(0.20123854, 0.38366989, 0.55429441)
			rgb=(0.1994295 , 0.38760678, 0.55464205)
			rgb=(0.1976365 , 0.39152762, 0.55496905)
			rgb=(0.19585993, 0.39543297, 0.55527637)
			rgb=(0.19410009, 0.39932336, 0.55556494)
			rgb=(0.19235719, 0.40319934, 0.55583559)
			rgb=(0.19063135, 0.40706148, 0.55608907)
			rgb=(0.18892259, 0.41091033, 0.55632606)
			rgb=(0.18723083, 0.41474645, 0.55654717)
			rgb=(0.18555593, 0.4185704 , 0.55675292)
			rgb=(0.18389763, 0.42238275, 0.55694377)
			rgb=(0.18225561, 0.42618405, 0.5571201 )
			rgb=(0.18062949, 0.42997486, 0.55728221)
			rgb=(0.17901879, 0.43375572, 0.55743035)
			rgb=(0.17742298, 0.4375272 , 0.55756466)
			rgb=(0.17584148, 0.44128981, 0.55768526)
			rgb=(0.17427363, 0.4450441 , 0.55779216)
			rgb=(0.17271876, 0.4487906 , 0.55788532)
			rgb=(0.17117615, 0.4525298 , 0.55796464)
			rgb=(0.16964573, 0.45626209, 0.55803034)
			rgb=(0.16812641, 0.45998802, 0.55808199)
			rgb=(0.1666171 , 0.46370813, 0.55811913)
			rgb=(0.16511703, 0.4674229 , 0.55814141)
			rgb=(0.16362543, 0.47113278, 0.55814842)
			rgb=(0.16214155, 0.47483821, 0.55813967)
			rgb=(0.16066467, 0.47853961, 0.55811466)
			rgb=(0.15919413, 0.4822374 , 0.5580728 )
			rgb=(0.15772933, 0.48593197, 0.55801347)
			rgb=(0.15626973, 0.4896237 , 0.557936 )
			rgb=(0.15481488, 0.49331293, 0.55783967)
			rgb=(0.15336445, 0.49700003, 0.55772371)
			rgb=(0.1519182 , 0.50068529, 0.55758733)
			rgb=(0.15047605, 0.50436904, 0.55742968)
			rgb=(0.14903918, 0.50805136, 0.5572505 )
			rgb=(0.14760731, 0.51173263, 0.55704861)
			rgb=(0.14618026, 0.51541316, 0.55682271)
			rgb=(0.14475863, 0.51909319, 0.55657181)
			rgb=(0.14334327, 0.52277292, 0.55629491)
			rgb=(0.14193527, 0.52645254, 0.55599097)
			rgb=(0.14053599, 0.53013219, 0.55565893)
			rgb=(0.13914708, 0.53381201, 0.55529773)
			rgb=(0.13777048, 0.53749213, 0.55490625)
			rgb=(0.1364085 , 0.54117264, 0.55448339)
			rgb=(0.13506561, 0.54485335, 0.55402906)
			rgb=(0.13374299, 0.54853458, 0.55354108)
			rgb=(0.13244401, 0.55221637, 0.55301828)
			rgb=(0.13117249, 0.55589872, 0.55245948)
			rgb=(0.1299327 , 0.55958162, 0.55186354)
			rgb=(0.12872938, 0.56326503, 0.55122927)
			rgb=(0.12756771, 0.56694891, 0.55055551)
			rgb=(0.12645338, 0.57063316, 0.5498411 )
			rgb=(0.12539383, 0.57431754, 0.54908564)
			rgb=(0.12439474, 0.57800205, 0.5482874 )
			rgb=(0.12346281, 0.58168661, 0.54744498)
			rgb=(0.12260562, 0.58537105, 0.54655722)
			rgb=(0.12183122, 0.58905521, 0.54562298)
			rgb=(0.12114807, 0.59273889, 0.54464114)
			rgb=(0.12056501, 0.59642187, 0.54361058)
			rgb=(0.12009154, 0.60010387, 0.54253043)
			rgb=(0.11973756, 0.60378459, 0.54139999)
			rgb=(0.11951163, 0.60746388, 0.54021751)
			rgb=(0.11942341, 0.61114146, 0.53898192)
			rgb=(0.11948255, 0.61481702, 0.53769219)
			rgb=(0.11969858, 0.61849025, 0.53634733)
			rgb=(0.12008079, 0.62216081, 0.53494633)
			rgb=(0.12063824, 0.62582833, 0.53348834)
			rgb=(0.12137972, 0.62949242, 0.53197275)
			rgb=(0.12231244, 0.63315277, 0.53039808)
			rgb=(0.12344358, 0.63680899, 0.52876343)
			rgb=(0.12477953, 0.64046069, 0.52706792)
			rgb=(0.12632581, 0.64410744, 0.52531069)
			rgb=(0.12808703, 0.64774881, 0.52349092)
			rgb=(0.13006688, 0.65138436, 0.52160791)
			rgb=(0.13226797, 0.65501363, 0.51966086)
			rgb=(0.13469183, 0.65863619, 0.5176488 )
			rgb=(0.13733921, 0.66225157, 0.51557101)
			rgb=(0.14020991, 0.66585927, 0.5134268 )
			rgb=(0.14330291, 0.66945881, 0.51121549)
			rgb=(0.1466164 , 0.67304968, 0.50893644)
			rgb=(0.15014782, 0.67663139, 0.5065889 )
			rgb=(0.15389405, 0.68020343, 0.50417217)
			rgb=(0.15785146, 0.68376525, 0.50168574)
			rgb=(0.16201598, 0.68731632, 0.49912906)
			rgb=(0.1663832 , 0.69085611, 0.49650163)
			rgb=(0.1709484 , 0.69438405, 0.49380294)
			rgb=(0.17570671, 0.6978996 , 0.49103252)
			rgb=(0.18065314, 0.70140222, 0.48818938)
			rgb=(0.18578266, 0.70489133, 0.48527326)
			rgb=(0.19109018, 0.70836635, 0.48228395)
			rgb=(0.19657063, 0.71182668, 0.47922108)
			rgb=(0.20221902, 0.71527175, 0.47608431)
			rgb=(0.20803045, 0.71870095, 0.4728733 )
			rgb=(0.21400015, 0.72211371, 0.46958774)
			rgb=(0.22012381, 0.72550945, 0.46622638)
			rgb=(0.2263969 , 0.72888753, 0.46278934)
			rgb=(0.23281498, 0.73224735, 0.45927675)
			rgb=(0.2393739 , 0.73558828, 0.45568838)
			rgb=(0.24606968, 0.73890972, 0.45202405)
			rgb=(0.25289851, 0.74221104, 0.44828355)
			rgb=(0.25985676, 0.74549162, 0.44446673)
			rgb=(0.26694127, 0.74875084, 0.44057284)
			rgb=(0.27414922, 0.75198807, 0.4366009 )
			rgb=(0.28147681, 0.75520266, 0.43255207)
			rgb=(0.28892102, 0.75839399, 0.42842626)
			rgb=(0.29647899, 0.76156142, 0.42422341)
			rgb=(0.30414796, 0.76470433, 0.41994346)
			rgb=(0.31192534, 0.76782207, 0.41558638)
			rgb=(0.3198086 , 0.77091403, 0.41115215)
			rgb=(0.3277958 , 0.77397953, 0.40664011)
			rgb=(0.33588539, 0.7770179 , 0.40204917)
			rgb=(0.34407411, 0.78002855, 0.39738103)
			rgb=(0.35235985, 0.78301086, 0.39263579)
			rgb=(0.36074053, 0.78596419, 0.38781353)
			rgb=(0.3692142 , 0.78888793, 0.38291438)
			rgb=(0.37777892, 0.79178146, 0.3779385 )
			rgb=(0.38643282, 0.79464415, 0.37288606)
			rgb=(0.39517408, 0.79747541, 0.36775726)
			rgb=(0.40400101, 0.80027461, 0.36255223)
			rgb=(0.4129135 , 0.80304099, 0.35726893)
			rgb=(0.42190813, 0.80577412, 0.35191009)
			rgb=(0.43098317, 0.80847343, 0.34647607)
			rgb=(0.44013691, 0.81113836, 0.3409673 )
			rgb=(0.44936763, 0.81376835, 0.33538426)
			rgb=(0.45867362, 0.81636288, 0.32972749)
			rgb=(0.46805314, 0.81892143, 0.32399761)
			rgb=(0.47750446, 0.82144351, 0.31819529)
			rgb=(0.4870258 , 0.82392862, 0.31232133)
			rgb=(0.49661536, 0.82637633, 0.30637661)
			rgb=(0.5062713 , 0.82878621, 0.30036211)
			rgb=(0.51599182, 0.83115784, 0.29427888)
			rgb=(0.52577622, 0.83349064, 0.2881265 )
			rgb=(0.5356211 , 0.83578452, 0.28190832)
			rgb=(0.5455244 , 0.83803918, 0.27562602)
			rgb=(0.55548397, 0.84025437, 0.26928147)
			rgb=(0.5654976 , 0.8424299 , 0.26287683)
			rgb=(0.57556297, 0.84456561, 0.25641457)
			rgb=(0.58567772, 0.84666139, 0.24989748)
			rgb=(0.59583934, 0.84871722, 0.24332878)
			rgb=(0.60604528, 0.8507331 , 0.23671214)
			rgb=(0.61629283, 0.85270912, 0.23005179)
			rgb=(0.62657923, 0.85464543, 0.22335258)
			rgb=(0.63690157, 0.85654226, 0.21662012)
			rgb=(0.64725685, 0.85839991, 0.20986086)
			rgb=(0.65764197, 0.86021878, 0.20308229)
			rgb=(0.66805369, 0.86199932, 0.19629307)
			rgb=(0.67848868, 0.86374211, 0.18950326)
			rgb=(0.68894351, 0.86544779, 0.18272455)
			rgb=(0.69941463, 0.86711711, 0.17597055)
			rgb=(0.70989842, 0.86875092, 0.16925712)
			rgb=(0.72039115, 0.87035015, 0.16260273)
			rgb=(0.73088902, 0.87191584, 0.15602894)
			rgb=(0.74138803, 0.87344918, 0.14956101)
			rgb=(0.75188414, 0.87495143, 0.14322828)
			rgb=(0.76237342, 0.87642392, 0.13706449)
			rgb=(0.77285183, 0.87786808, 0.13110864)
			rgb=(0.78331535, 0.87928545, 0.12540538)
			rgb=(0.79375994, 0.88067763, 0.12000532)
			rgb=(0.80418159, 0.88204632, 0.11496505)
			rgb=(0.81457634, 0.88339329, 0.11034678)
			rgb=(0.82494028, 0.88472036, 0.10621724)
			rgb=(0.83526959, 0.88602943, 0.1026459 )
			rgb=(0.84556056, 0.88732243, 0.09970219)
			rgb=(0.8558096 , 0.88860134, 0.09745186)
			rgb=(0.86601325, 0.88986815, 0.09595277)
			rgb=(0.87616824, 0.89112487, 0.09525046)
			rgb=(0.88627146, 0.89237353, 0.09537439)
			rgb=(0.89632002, 0.89361614, 0.09633538)
			rgb=(0.90631121, 0.89485467, 0.09812496)
			rgb=(0.91624212, 0.89609127, 0.1007168 )
			rgb=(0.92610579, 0.89732977, 0.10407067)
			rgb=(0.93590444, 0.8985704 , 0.10813094)
			rgb=(0.94563626, 0.899815 , 0.11283773)
			rgb=(0.95529972, 0.90106534, 0.11812832)
			rgb=(0.96489353, 0.90232311, 0.12394051)
			rgb=(0.97441665, 0.90358991, 0.13021494)
			rgb=(0.98386829, 0.90486726, 0.13689671)
			rgb=(0.99324789, 0.90615657, 0.1439362 )
		}
}

\pgfplotsset{colormap={viridisSoft}{
			rgb255=(242, 242, 242);
rgb=(0.28026,0.1657,0.4765);
			rgb=(0.26366,0.23763,0.51877);
			rgb=(0.23744,0.3052,0.54192);
			rgb=(0.20862,0.36775,0.55267);
			rgb=(0.18225,0.42618,0.55711);
			rgb=(0.1592,0.48224,0.55807);
			rgb=(0.13777,0.53749,0.5549);
			rgb=(0.12115,0.59274,0.54465);
			rgb=(0.12808,0.64775,0.5235);
			rgb=(0.18065,0.7014,0.48819);
			rgb=(0.27415,0.75198,0.4366);
			rgb=(0.39517,0.79747,0.36775);
			rgb=(0.53561,0.83578,0.2819);
			rgb=(0.68895,0.86545,0.18272);
			rgb=(0.84557,0.88733,0.0997);
			rgb=(0.99324,0.90616,0.14394)
		}
}

\pgfplotsset{colormap={cellRed}{
			rgb255=(242.0,242.0,242.0);
			rgb255=(241.63157894736844,234.47368421052633,234.47368421052633);
			rgb255=(241.26315789473685,226.94736842105266,226.94736842105266);
			rgb255=(240.89473684210526,219.42105263157893,219.42105263157893);
			rgb255=(240.5263157894737,211.89473684210526,211.89473684210526);
			rgb255=(240.1578947368421,204.3684210526316,204.3684210526316);
			rgb255=(239.78947368421052,196.84210526315792,196.84210526315792);
			rgb255=(239.42105263157896,189.31578947368422,189.31578947368422);
			rgb255=(239.05263157894737,181.78947368421052,181.78947368421052);
			rgb255=(238.6842105263158,174.26315789473688,174.26315789473688);
			rgb255=(238.31578947368422,166.73684210526315,166.73684210526315);
			rgb255=(237.94736842105263,159.21052631578948,159.21052631578948);
			rgb255=(237.57894736842104,151.68421052631578,151.68421052631578);
			rgb255=(237.21052631578948,144.1578947368421,144.1578947368421);
			rgb255=(236.84210526315792,136.63157894736844,136.63157894736844);
			rgb255=(236.47368421052633,129.10526315789474,129.10526315789474);
			rgb255=(236.10526315789474,121.57894736842107,121.57894736842107);
			rgb255=(235.73684210526318,114.05263157894737,114.05263157894737);
			rgb255=(235.3684210526316,106.52631578947368,106.52631578947368);
			rgb255=(235.0,99.0,99.0);
		}
}

\pgfplotsset{colormap={cellGreen}{
			rgb255=(242.0,242.0,242.0);
			rgb255=(236.21052631578948,239.5263157894737,234.26315789473685);
			rgb255=(230.42105263157896,237.05263157894737,226.5263157894737);
			rgb255=(224.6315789473684,234.57894736842104,218.78947368421052);
			rgb255=(218.8421052631579,232.10526315789474,211.05263157894737);
			rgb255=(213.05263157894737,229.63157894736844,203.31578947368422);
			rgb255=(207.26315789473685,227.1578947368421,195.57894736842107);
			rgb255=(201.4736842105263,224.68421052631578,187.8421052631579);
			rgb255=(195.68421052631578,222.21052631578948,180.10526315789474);
			rgb255=(189.8947368421053,219.73684210526318,172.36842105263162);
			rgb255=(184.10526315789474,217.26315789473682,164.63157894736844);
			rgb255=(178.31578947368422,214.78947368421052,156.89473684210526);
			rgb255=(172.5263157894737,212.31578947368422,149.1578947368421);
			rgb255=(166.73684210526318,209.84210526315792,141.42105263157896);
			rgb255=(160.94736842105263,207.3684210526316,133.6842105263158);
			rgb255=(155.1578947368421,204.89473684210526,125.94736842105263);
			rgb255=(149.3684210526316,202.42105263157893,118.21052631578948);
			rgb255=(143.57894736842104,199.94736842105266,110.47368421052632);
			rgb255=(137.78947368421052,197.47368421052633,102.73684210526316);
			rgb255=(132.0,195.0,95.0);
		}
}

\pgfplotsset{colormap={cellRedSquared}{
			rgb255=(242.0,242.0,242.0);
			rgb255=(241.28254847645428,227.34349030470915,227.34349030470915);
			rgb255=(240.60387811634348,213.47922437673128,213.47922437673128);
			rgb255=(239.9639889196676,200.40720221606648,200.40720221606648);
			rgb255=(239.36288088642658,188.1274238227147,188.1274238227147);
			rgb255=(238.8005540166205,176.63988919667594,176.63988919667594);
			rgb255=(238.2770083102493,165.94459833795014,165.94459833795014);
			rgb255=(237.79224376731304,156.04155124653738,156.04155124653738);
			rgb255=(237.34626038781164,146.93074792243766,146.93074792243766);
			rgb255=(236.93905817174516,138.61218836565098,138.61218836565098);
			rgb255=(236.57063711911357,131.0858725761773,131.0858725761773);
			rgb255=(236.2409972299169,124.35180055401662,124.35180055401662);
			rgb255=(235.95013850415512,118.40997229916897,118.40997229916897);
			rgb255=(235.69806094182823,113.26038781163435,113.26038781163435);
			rgb255=(235.4847645429363,108.90304709141274,108.90304709141274);
			rgb255=(235.3102493074792,105.33795013850416,105.33795013850416);
			rgb255=(235.17451523545705,102.56509695290858,102.56509695290858);
			rgb255=(235.0775623268698,100.58448753462605,100.58448753462605);
			rgb255=(235.01939058171746,99.3961218836565,99.3961218836565);
			rgb255=(235.0,99.0,99.0);
		}
}
\pgfplotsset{colormap={cellGreenSquared}{
			rgb255=(242.0,242.0,242.0);
			rgb255=(230.7257617728532,237.18282548476455,226.93351800554018);
			rgb255=(220.06094182825484,232.62603878116343,212.6814404432133);
			rgb255=(210.00554016620498,228.32963988919667,199.2437673130194);
			rgb255=(200.5595567867036,224.29362880886427,186.62049861495845);
			rgb255=(191.7229916897507,220.5180055401662,174.8116343490305);
			rgb255=(183.49584487534625,217.0027700831025,163.81717451523545);
			rgb255=(175.87811634349032,213.74792243767314,153.63711911357342);
			rgb255=(168.86980609418282,210.75346260387812,144.27146814404432);
			rgb255=(162.47091412742384,208.01939058171746,135.72022160664818);
			rgb255=(156.68144044321332,205.54570637119116,127.98337950138506);
			rgb255=(151.50138504155126,203.33240997229916,121.06094182825484);
			rgb255=(146.9307479224377,201.37950138504155,114.95290858725764);
			rgb255=(142.96952908587255,199.68698060941827,109.65927977839334);
			rgb255=(139.61772853185596,198.25484764542935,105.18005540166205);
			rgb255=(136.8753462603878,197.0831024930748,101.51523545706371);
			rgb255=(134.74238227146813,196.17174515235456,98.66481994459834);
			rgb255=(133.21883656509695,195.5207756232687,96.62880886426592);
			rgb255=(132.30470914127426,195.13019390581718,95.40720221606648);
			rgb255=(132.0,195.0,95.0);
		}
} 
\usepgfplotslibrary{patchplots}

\pgfplotsset{every axis/.append style={
			grid=both,
			grid style={white, line width=.1pt},
			major grid style={white, line width=1.5pt},
			axis background/.style={fill=gray!10},
			axis line style={draw=none},
			tick style={draw=none},
			xlabel = $x$,
line width=1pt,
legend style={
					line width = 1pt,
					draw=none,
					/tikz/every even column/.append style={column sep=0.5cm}
				},
		}}

\usepgfplotslibrary{groupplots}

\usepackage{calc}

\definecolor{gg0}{HTML}{E24A33}
\definecolor{gg1}{HTML}{348ABD}
\definecolor{gg2}{HTML}{988ED5}
\definecolor{gg3}{HTML}{777777}
\definecolor{gg4}{HTML}{FBC15E}
\definecolor{gg5}{HTML}{8EBA42}
\definecolor{gg6}{HTML}{FFB5B8}

\pgfplotsset{
	/pgfplots/colormap={bright}{rgb255=(0,0,0) rgb255=(78,3,100) rgb255=(2,74,255)
			rgb255=(255,21,181) rgb255=(255,113,26) rgb255=(147,213,114) rgb255=(230,255,0)
			rgb255=(255,255,255)}
}

\newcommand{\addappendix}{\section*{\appendixname}\addcontentsline{toc}{section}{\appendixname}\counterwithin*{figure}{section}\stepcounter{section}\renewcommand{\thesection}{A}\renewcommand{\thefigure}{\thesection.\arabic{figure}}}

\def\bm{\boldsymbol}

\definecolor{brandeisblue}{rgb}{0.0, 0.44, 1.0}
\definecolor{lincolngreen}{rgb}{0.11, 0.35, 0.02}
\definecolor{indiagreen}{rgb}{0.07, 0.53, 0.03}
\definecolor{venetianred}{rgb}{0.78, 0.03, 0.08}
\definecolor{darkorange}{rgb}{1.0, 0.55, 0.0}
\definecolor{burntorange}{rgb}{0.8, 0.33, 0.0}
\definecolor{flame}{rgb}{0.89, 0.35, 0.13}
\definecolor{non-photoblue}{rgb}{0.64, 0.87, 0.93}

\usepackage{tikz}
\usepackage{tkz-base}
\usetikzlibrary{calc}
\usepackage{tkz-euclide}
\usepackage{xcolor}

\renewcommand{\review}[2]{}
\renewcommand{\creview}[3]{}
\renewcommand{\ntcreview}[3]{}

\renewcommand{\tableofcontents}{}
\renewcommand{\listofreviews}{}

\definecolor{revisionColourOne}{RGB}{180,0,0}
\definecolor{revisionColourTwo}{RGB}{0,0,180}

\newcommand{\compEmb}{\hookrightarrow \!\!\!\rightarrow}

\usepackage{fullpage}

\usepackage{setspace}
\begin{document}
\begin{singlespace}\maketitle\end{singlespace}
\begin{abstract}
We construct a finite element method for the numerical solution of a fractional porous medium equation on a bounded open Lipschitz polytopal domain $\Omega \subset \R^{d}$, where $d = 2$ or $3$. The pressure in the model is defined as the solution of a fractional Poisson equation, involving the fractional Neumann Laplacian in terms of its spectral definition. We perform a rigorous passage to the limit as the spatial and temporal discretization parameters tend to zero and show that a subsequence of the sequence of finite element approximations defined by the proposed numerical method converges to a bounded and nonnegative weak solution of the initial-boundary-value problem under consideration. This result can be therefore viewed as a constructive proof of the existence of a nonnegative, energy-dissipative, weak solution to the initial-boundary-value problem for the fractional porous medium equation under consideration, based on the Neumann Laplacian. The convergence proof relies on results concerning the finite element approximation of the spectral fractional Laplacian and compactness techniques for nonlinear partial differential equations, together with properties of the equation, which are shown to be inherited by the numerical method. We also prove that the total energy associated with the problem under consideration exhibits exponential decay in time. 
\end{abstract}
 \subjectclassification{\subjectPDF 35K55; 35R11; 65N30.}

\section{Introduction}\label{Sect:1}
The aim of this work is to construct and analyze a numerical scheme for a porous medium equation involving the fractional Laplacian, in order to model nonlocal features. Let us denote by $\rho$ a certain density distribution and by $\bm{v}$ a velocity field; both are functions of space and time, $\rho = \rho(x, t)$, $\bm{v}= \bm{v}(x,t)$, the velocity is a vector function and we assume the density to be nonnegative, $\rho(x, t) \geq 0$. The density and the velocity field are related by the continuity equation
\begin{subequations} \label{MainProb}
\begin{equation}
    \frac{\partial \rho}{\partial t} + \nabla \cdot (\rho \bm{v}) = 0.
\end{equation}
We postulate that the velocity $\bm{v}$ derives from a potential $c=c(x, t)$ through the equality 
\begin{equation}\label{eq:cv}
    \bm{v} =  -\nabla c.
\end{equation}

This happens for instance when $\rho$ is the density of a fluid and motion is modelled by Darcy's law, which relates velocity and pressure of the fluid (in that case $c$ denotes the pressure). The potential can then be related in a closed form to the density $c = c(\rho)$ through a suitable equation of state, leading to various different models. For example, the simplest case $c = 2k \rho$ in fluid flow modeling, with $k$ signifying a positive constant, pertains to filtration of a fluid through a porous stratum in the isothermal setting, that is the standard model for groundwater infiltration, Boussinesq's equation $\rho_{t} = k \; \Delta\rho^{2}$ (see \cite{bear2013dynamics}). More generally one has  the model $\rho_{t} = k \; \Delta \rho^{m}$ for $m>1$. 

The structure of the porous medium can be very complex, however. Different pores, tubes and water filaments can make the phenomenon difficult to study, and standard models fail in accurately  describing these complex features. A number of experimentalists have shown that in certain porous media, such as construction materials, classical diffusion cannot accurately predict the evolution of moisture. It was observed in particular that the wetting front may move in a significantly different way than  classical models predict. 
For instance, in \cite{SCHUMER200169} the authors discuss solute transport in aquifers and they explain that standard dispersion can only occur in homogeneous aquifers.  To address this issue, different models have been proposed in which nonlocal interactions and effects are present. For example, in \cite{mihir2012} the authors propose a porous medium modelled with a network of channels of all length-scales through which the water particles can
move over a very long distance. 

In recent years there has been considerable interest in the study of nonlocal operators with the aim to incorporate long-range interactions. The paradigm of such operators is the fractional Laplacian $(-\Delta)^{\s}$ with $s \in (0,1)$. The use of the fractional Laplacian in diffusive processes is discussed in a large body of literature and the reader may wish to refer to \cite{antil2017fractional, chen2022analysis, escudero2006fractional} for some examples of its application. 

We will consider as the object of our study a nonlocal porous medium equation, in which the long-range interactions are modelled by a fractional-order elliptic equation; more precisely, we will consider
\begin{equation} \label{FracPoiZeroZero}
    -(-\Delta)^{\s} c = \rho, \quad 0< s <1.
\end{equation}
\end{subequations}
If $s=0$ we recover the standard porous medium equation, while the case $s=1$ has been investigated in \cite{markowich1990,chapman1996mean, weinan1994dynamics}. 

The model (\ref{MainProb}) has been justified in the hydrological setting in \cite{MR3940331}, where a deterministic derivation of the equation was presented, subject to the assumption that the fluid in the porous medium can contribute to the flux at any point with long jumps depending only on the distance. The system (\ref{MainProb}) also arises in other contexts: a model similar to the case $\s = 1$ was considered in \cite{ambrosio2008gradient} to describe the evolution of vortex-density in superconductors, while the equation with $\s = \frac{1}{2}$ was suggested in \cite{head1972dislocation} to model the motion of dislocations in a continuum. These equations also find applications in modelling charged particle transport in materials, since they are fractional counterparts of classical drift-diffusion equations for semiconductors \cite{markowich1990}. For further discussion concerning models of porous medium diffusion with a nonlocal pressure of this type we refer to Section 5.6.1 of \cite{MR3588125}, where it is noted that the problem in a bounded domain with Dirichlet or Neumann boundary data has not been studied and computational techniques to deal with these kinds of nonlocal effects and the design of generic numerical schemes for this model are needed. 

Instead of the first-order hyperbolic equation \eqref{MainProb} we will consider a parabolic regularization of the problem, where standard Fickian diffusion is present with the positive diffusion coefficient scaled to $1$, because its specific value is of no relevance to the discussion herein. This leads to the following system:
\begin{equation} \label{FullProbWSpace}
\left \{
\begin{aligned}
&\frac{\partial \rho}{\partial t} = \Delta \rho - \nabla \cdot (\rho \nabla c), \\
& - (-\Delta)^{\s} c = \rho.
\end{aligned}
\right.
\end{equation}
The fractional porous medium equation on $\mathbb{R}^{d}$ is obtained by inverting the fractional operator $(-\Delta)^{\s}$ in terms of its fundamental solution $U_s$, i.e., $-c=U_s\ast \rho$, and plugging this into the first equation in \eqref{FullProbWSpace}, which then yields
\begin{equation} \label{FullProbWSpace2}
\frac{\partial \rho}{\partial t} = \Delta \rho + \nabla \cdot (\rho \nabla (U_s\ast\rho)). \\
\end{equation}
The previous equation has formally the structure of a Wassertein gradient flow of an energy functional $\mathcal{E}(\rho)$, as in \cite{carrillovillani2003,AGS,carrillovazquez2015} for instance, if mass is normalized, thanks to it being conserved. Consequently, it leads to an evolution of densities governed by the dissipation of the free energy functional
\begin{equation*} 
\mathcal{E}(\rho) := \int_{\mathbb{R}^{d}} \rho \log \rho \dx + \frac{1}{2} \int_{\mathbb{R}^{d}} (U_s\ast \rho) \rho \dx , 
\end{equation*}
where the first term is the Boltzmann entropy and the second term is the total interaction potential energy. 

It has been shown in \cite{caffarelli2010nonlinear,caffarelli2013regularity,chen2022analysis} that under suitable assumptions there exists a weak solution to \eqref{MainProb} and \eqref{FullProbWSpace} on the whole of $\mathbb{R}^{d}$, which has several properties, including conservation of mass, the dissipation of the free energy functional, and the decay of the spatial $L^{\infty}$ norm in time. We will show that our numerical scheme has some of these features as well. Moreover in the recent works \cite{caffarelli2010asymptotic, carrillo2015exponential} the asymptotic behaviour of solutions to \eqref{MainProb} has been investigated where the dissipation of the free energy functional plays a key role. The existence of a smooth solution on $\mathbb{R}^{d}$ to the equation  (\ref{FullProbWSpace}) or equivalently \eqref{FullProbWSpace2}, supplemented with an initial condition,  has been proved in \cite{choi2021classical}. Indeed, a large portion of the available literature for this type of equation considers the model on the whole of $\mathbb{R}^{d}$ written in the form \eqref{FullProbWSpace2} as it uses the integral definition of the fractional Laplacian as a convolution with the Riesz kernel; see \cite{bucur2016nonlocal} for further details. Many definitions of the fractional Laplacian are available, both on the whole of $\mathbb{R}^{d}$ and on a bounded domain $\Omega \subset \mathbb{R}^{d}$ (cf. \cite{lischke2020fractional}), and while the definitions are equivalent in the former case they are not interchangeable in the latter. 

The novelty of our work consists in the first place in considering the model (\ref{FullProbWSpace}) on a bounded domain $\Omega \subset \mathbb{R}^{d}$, and secondly in using the spectral fractional Laplacian for the potential instead of the integral definition. These two aspects, together with our use of a provably convergent numerical method for this fractional-order porous medium model, based on a spatial continuous piecewise affine finite element approximation and implicit Euler time-stepping, constitute the main contributions of our work.

The definition of the finite element approximation of the spectral fractional Laplacian (cf. Section \ref{Sect:3}) involves the complete set of eigenvalues and eigenfunctions of the discrete Laplacian, and its implementation is therefore impracticable, except in very special cases (e.g. when $\Omega$ is a square in $\mathbb{R}^2$ or a cube in $\mathbb{R}^3$, with a homogeneous Dirichlet, homogeneous Neumann, or a periodic boundary condition imposed on $\partial\Omega$), when the complete set of eigenvalues and eigenfunctions can be explicitly calculated. It turns out however that an accurate approximation of the discrete fractional spectral Laplacian can be efficiently computed by using a technique based on rational approximation, which only requires accurate estimates of the smallest and largest eigenvalue of the discrete Laplacian. The details of the implementation of this technique will be the focus of our forthcoming, computationally-oriented, publication \cite{carrillosuli2024}. For further details concerning the finite element approximation
   of the fractional Laplacian the reader may wish to consult \cite{bonito2018numerical, Bonito2019, NochettoBorthagaray2019, Bersetche2022}.

A key novelty of our approach is that we adapt to the case of a bounded domain the variational viewpoint that was successfully used when the spatial domain is the whole of $\mathbb{R}^d$, making an effective use of the \textit{a priori} estimates at all levels of discretization. The system \eqref{FullProbWSpace} has an analogous variational structure in bounded domains as when the spatial domain is the whole of $\mathbb{R}^d$. In fact, we will show that a related free energy functional of the form
\begin{equation*}
E(\rho) := \int_{\Omega} G(\rho)\dx - \frac{1}{2} \int_{\Omega} c \rho \dx, 
\end{equation*}
is dissipated by the system \eqref{FullProbWSpace} subject to Neumann boundary conditions imposed on $\rho$ and $c$, and with the spectral fractional Laplacian. Similar free energy functionals have been considered in the related Keller--Segel models for chemotaxis or gravitational models with Fickian diffusion in bounded domains and Neumann boundary conditions, to obtain regularity estimates and to analyse their asymptotic behavior; see \cite{bilernadzieja94,biler99-2} for instance. In our case, the repulsion effect due to the forcing term obtained from the spectral fractional Laplacian in \eqref{FullProbWSpace} allows us to obtain further crucial \textit{a priori} estimates. We also prove that, for any convex function $F \in C^{2}([0, \infty))$, 
\[ 
\frac{\mathrm{d}}{\mathrm{d}t} \int_{\Omega} F(\rho)\dx \leq 0, 
\]
where $\rho$ is a solution to the system \eqref{FullProbWSpace} on $\Omega$ subject to Neumann boundary conditions imposed on both $\rho$ and $c$.
This dissipation of the functional $\rho \mapsto \int_\Omega F(\rho)\dx$ along its evolution allows us to obtain pointwise control of the density uniformly in time from its initial datum at the semidiscrete level, which is of paramount importance in the proof of the existence of weak solutions, as well as in our study of the asymptotic behaviour of solutions to both the numerical scheme and the limiting PDE system \eqref{FullProbWSpace}.

The outline of the paper is as follows. In the next section we introduce the relevant function spaces and the statement of the problem under consideration. We propose a weak formulation of the initial-boundary-value problem for the fractional porous medium equation and its regularized version involving two positive cut-off parameters, $\delta$ and $L$, corresponding to cut-off from below and from above, respectively. We then discuss key properties of the model, including conservation of mass and the boundedness in time of the spatial $L^{\infty}$ norm. In Section \ref{Sect:3} we introduce the time-discretization and the spatial finite element approximation. We then show that a subsequence of this sequence of approximations converges to a bounded and nonnegative weak solution of the semidiscrete-in-time problem, as the spatial discretization parameter $h$ and the cut-off parameter $\delta$ tend to zero. 
In Section \ref{Sect:4}, after discussing further properties of the solution of the semidiscrete-in-time scheme, including in particular a bound on the numerical solution that is independent of the cut-off parameter $L$, we prove that a subsequence of this sequence of approximations converges to a bounded and nonnegative weak solution of our initial-boundary-value problem, as the time discretization parameter $\Delta t$ tends to zero. We also show that the total energy associated with the problem under consideration exhibits exponential decay in time from a specified nonnegative and bounded initial density. We close in Section \ref{Sect:5} with concluding remarks and a list of relevant open problems.

\section{The fractional porous medium model}\label{Sect:2}

We shall supplement the partial differential equations \eqref{FullProbWSpace} with an initial condition for $\rho$ and boundary conditions for $\rho$ and $c$; the most natural choice for the latter are no-flux boundary conditions for both $\rho$ and $c$.  Suppose that $\Omega \subset \mathbb{R}^d$ is a bounded open Lipschitz domain; in the next section, simply for ease of triangulating $\Omega$, we shall confine ourselves to the case when $\Omega$ is a polytopal Lipschitz domain in $\mathbb{R}^d$, $d=2,3$. The initial-boundary-value problem we shall be considering is the following:  

\begin{equation} \label{MainFullProb}
\left \{
\begin{aligned}
&\frac{\partial \rho}{\partial t} = \Delta \rho - \nabla \cdot (\rho \nabla c) & \textrm{in } & \Omega \times (0, \infty) , \\
& - (-\Delta)^{\s} c = \rho^{\ast} & \textrm{in } & \Omega \times (0, \infty), \\
&\partial_{n} \rho = 0, \quad \partial_{n} c = 0 & \textrm{on~\!} & \partial \Omega \times (0, \infty), \\
&\rho(x,0) = \rho_0(x)\quad \mbox{for all $x \in \Omega$}, &
\end{aligned}
\right.
\end{equation}
where $\rho_0 \in L^\infty(\Omega)$ is a given nonnegative initial datum, with normalized mass $\frac{1}{|\Omega|} \int_{\Omega} \rho_{0} \dx = 1$, the fractional Laplacian $(-\Delta)^{\s}$ is understood to be based on its spectral definition, which we will state later on, and $\rho^{\ast}$ stands for the projection of $\rho$ into the space of functions with zero integral average, namely 
\begin{equation} \label{ProjZeroAv} \rho^{\ast} := \rho - \overline{\rho}, \quad \text{with } \overline{\rho} := \frac{1}{|\Omega|} \int_{\Omega} \rho \dx . \end{equation}

We denote by $(\cdot, \cdot)$ the $L^{2}(\Omega)$ inner product and  by $\|\cdot\|_{L^{p}(\Omega)}$ and $\|\cdot\|_{H^{1}(\Omega)}$ the standard $L^{p}(\Omega)$ and $H^{1}(\Omega)$ norm, respectively. 
Moreover let $H^{-1}(\Omega)$ denote the dual space of $H^{1}(\Omega)$; we denote by $\langle \cdot , \cdot \rangle$ the associated duality pairing between $H^{-1}(\Omega)$ and $H^{1}(\Omega)$.

Let us denote by $L^{2}_{\ast}(\Omega)$ and  $H^{1}_{\ast}(\Omega)$ the closed linear subspace of $L^{2}(\Omega)$ and $H^{1}(\Omega)$ with zero integral average on $\Omega$, respectively; that is, 
$$L^{2}_{\ast}(\Omega) := \bigg\{ v \in L^{2}(\Omega) \text{ such that } \int_{\Omega}v \dx = 0 \bigg\}\quad
\mbox{and}  
\quad
H^{1}_{\ast}(\Omega) := \bigg\{ v \in H^{1}(\Omega) \textrm{ such that } \int_{\Omega} v \dx= 0 \bigg\}.$$
Let $H^{-1}_{\ast}(\Omega)$ denote the dual space of $H^{1}_{\ast}(\Omega)$. For $s \in (0,1)$ we define the fractional-order Sobolev space as

\begin{equation*}
H^{\s}(\Omega) := \bigg \{ u \in L^{2}(\Omega) \textrm{ such that } \int_{\Omega} \int_{\Omega} \frac{|u(x) - u(y)|^{2}}{|x - y|^{d + 2 \s}} \dx \, \dy < \infty \bigg \},
\end{equation*}
and we endow it with the norm defined by 
\begin{equation*}
\|u\|_{H^{\s}(\Omega)} := \big( \|u\|_{L^{2}(\Omega)}^{2} + |u|^{2}_{H^{\s}(\Omega)} \big)^{\frac{1}{2}},
\end{equation*}
where the Gagliardo--Slobodecki\u{\i} semi-norm $| \cdot |_{H^{\s}(\Omega)}$ is given by 
\begin{equation*}
|u|_{H^{\s}(\Omega)}^{2} := \int_{\Omega} \int_{\Omega} \frac{|u(x) - u(y)|^{2}}{|x - y|^{d + 2 \s}} \dx \dy.
\end{equation*}
For $\sigma = m + s$, where $m$ is a nonnegative integer and $s \in (0,1)$, the Sobolev space of order $\sigma>0$ is defined by 
\begin{equation*}
H^{\sigma}(\Omega) := \big\{ u \in H^{m}(\Omega) \text{ such that } \mathcal{D}^{\alpha} u \in H^{s}(\Omega) \text{ for all } \alpha \in \mathbb{N}^{d},\, |\alpha|=m \big\},
\end{equation*} 
where $H^{m}(\Omega)$ is the Sobolev space of order $m$ and summability 2, $\mathcal{D}^{\alpha}$ stands for the weak derivative of multi-index $\alpha=(\alpha_1,\ldots,\alpha_d)$ and order $|\alpha|:=\alpha_1+ \cdots + \alpha_d$, and $\mathbb{N}$ denotes the set of all nonnegative integers. We will denote by  $H^{-s}(\Omega)$ the dual space of the fractional-order Sobolev space $H^{s}(\Omega)$. Similarly as before, we let $H^{\s}_{\ast}(\Omega)$ be the set of all functions in $H^{\s}(\Omega)$ with zero integral average, and we denote by $H^{-\s}_{\ast}(\Omega)$ its dual space, $s \in (0,1)$.

As was indicated previously, many definitions of the fractional Laplacian are available; for the purposes of our work here we shall consider the spectral definition of the fractional Neumann Laplacian, based on using the eigenvalues and eigenfunctions of the Neumann Laplacian. The Neumann Laplacian $(-\Delta_\mathrm{N}, \mathrm{Dom}(-\Delta_\mathrm{N}))$ on a Lipschitz domain $\Omega$ is defined
as follows: 
\begin{align*}
 \mathrm{Dom}(-\Delta_\mathrm{N}) := \bigg\{ u \in H^1(\Omega)\,:\, H^1(\Omega)\ni & \,v \mapsto \!\int_\Omega \nabla u \cdot \nabla v \dx\quad \mbox{is a continuous linear funcional on } L^2(\Omega)\bigg \},\\
 \int_\Omega (-\Delta_\mathrm{N} u)\, v &:= \int_\Omega \nabla u \cdot \nabla v \dx, \quad u \in \mathrm{Dom}(-\Delta_\mathrm{N}),\, \, v \in H^1(\Omega).
\end{align*}
Thus, $-\Delta_\mathrm{N}$ is a nonnegative and selfadjoint linear operator that is a densely defined in $L^2(\Omega)$. 
There exists therefore an orthonormal basis of $L^2(\Omega)$ consisting of eigenfunctions $\psi_k \in H^1(\Omega)$, $k = 0,1,2,\ldots$, with corresponding eigenvalues $0 = \lambda_0 < \lambda_1 \leq \lambda_2 \leq \cdots \nearrow +\infty$. 
Thus, 
\begin{equation}\label{eigen}
    \left \{ \begin{aligned}
        -\Delta \psi_k &= \lambda_k \psi_k &&\quad \textrm{in } \Omega,  \\
        \partial_{n} \psi_k &= 0 &&\quad\textrm{on } \partial \Omega, 
    \end{aligned} \right.
\end{equation}
with the partial differential equation in the first line of \eqref{eigen} understood as an equality in $L^2(\Omega)$, and the homogeneous Neumann boundary condition appearing in the second line of \eqref{eigen} as an equality in $H^{-\frac{1}{2}}(\partial\Omega)$ in the sense of the trace theorem $\partial_n\,: u \in \{v \in L^2(\Omega)\,: \Delta v \in L^2(\Omega)\} \mapsto \partial_n u \in H^{-\frac{1}{2}}(\partial \Omega)$.
Following the discussion in Sec.~7 of \cite{caffarelli2016fractional},
for $\s \in (0, 1)$ the domain of the spectral fractional Neumann Laplacian $(-\Delta_{\mathrm{N}})^s$ is defined as the Hilbert space
\begin{equation*}
\Hsast(\Omega) := \left \{ u(\cdot) = \sum_{k=1}^{\infty} u_{k} \psi_{k}(\cdot) \in L^{2}_{\ast}(\Omega): \|u\|^{2}_{\Hsast(\Omega)} := \sum_{k=1}^{\infty} \lambda_{k}^{\s} u_{k}^{2} < \infty \right\} \quad \textrm{with } u_{k} := \int_{\Omega} u(x) \psi_{k}(x) \dx
\end{equation*}
for $k=1,2,\ldots$, and, for $u \in \mathbb{H}^s_\ast(\Omega)$, 
\begin{equation}
\label{FracLapNeu}
(-\Delta_{\mathrm{N}})^s u(x) : = \sum_{k=1}^\infty \lambda_k^s u_k \psi_k(x),
\end{equation}
as an element of $\mathbb{H}^{-s}_\ast(\Omega)$. Furthermore,  one has 
\begin{equation*} \|u \|_{\Hsast(\Omega)} = \| (-\Delta_\mathrm{N})^{\s / 2} u\|_{L^{2}(\Omega)}\quad \forall\, u \in \mathbb{H}^s_\ast(\Omega). 
\end{equation*}

\begin{remark} \label{RemFracSobs}
We note that any function $u$ belonging to $\Hsast(\Omega)$, $s \in (0,1)$, satisfies $\int_{\Omega} u \dx = 0$. In addition, thanks to   Lemma 7.1 in \cite{caffarelli2016fractional}, for any $s \in (0,1)$, a function $u$ belongs to $\Hsast(\Omega)$ if, and only if, $u \in H^{\s}(\Omega)$ and $\int_{\Omega} u \dx = 0$; furthermore, the norm $\|\cdot\|_{\Hsast(\Omega)}$ is equivalent to the seminorm $|\cdot|_{H^{\s}(\Omega)}$, which is a norm on $H^{\s}_{\ast}(\Omega)$.
\end{remark}

Having stated the problem (\ref{MainFullProb}) we now propose its weak formulation. Let $V:= H^{1}(\Omega)$. The weak formulation of the initial-boundary-value problem \eqref{MainFullProb} is as follows:
\begin{subequations} \label{BasicWeakFormMain}
\begin{gather}
\text{Find } \rho \in L^2(0,T;V)  \text{ with $\frac{\partial \rho}{\partial t} \in L^\infty(0,T;V')$ such that} \nonumber \\
\Big\langle \frac{\partial \rho}{\partial t}, \phi \Big\rangle = - \int_{\Omega} \nabla \rho \cdot \nabla \phi \dx + \int_{\Omega} \rho \nabla c \cdot \nabla \phi \dx \quad \text{for all } \phi \in V \text{ and a.e. $t \in (0,T]$}, \label{BasicWeakForm}
\end{gather}
subject to the initial condition $\rho(x, 0) = \rho_{0}(x)$, where $\rho_{0} \in L^{\infty}(\Omega)$ and $\rho_{0}(x) \geq 0$ for a.e. $x \in \Omega$ and where
\begin{equation} \label{MyFracPois1}
    -(-\Delta_{\mathrm{N}})^{\s} c = \rho^{\ast} \textrm{ in } \Omega, 
\end{equation}
\end{subequations}
with $\rho^{\ast}$ given by (\ref{ProjZeroAv}). Notice that $c$ is well defined, as the solution to the fractional equation (\ref{MyFracPois1}) exists and it is unique by Theorem \ref{FracPoiExUniq}; we will therefore write $c = -(-\Delta_{\mathrm{N}})^{-\s}\rho^{\ast}$, where $(-\Delta_{\mathrm{N}})^{-s}$ stands for the inverse of the spectral fractional Neumann Laplacian  (\ref{FracLapNeu}). When we speak of a weak solution to the problem (\ref{BasicWeakFormMain}) we mean a couple $\rho$ and $c$ such that (\ref{BasicWeakForm}) and (\ref{MyFracPois1}) hold. 

Now let us assume that a weak solution to problem (\ref{BasicWeakFormMain}) exists and that $\rho$ is nonnegative almost everywhere in $\Omega \times [0,T]$. The existence of a weak solution to our problem and its nonnegativity will be rigorously proved in the next sections. In order to motivate the argument which our proof is based upon, we shall now perform preliminary formal computations; in these formal computations we shall additionally assume that the functions $\rho$ and $c$ have all the regularity that we need in order to ensure that our calculations are meaningful. 

The first remark we make is that, if we take $\phi \equiv 1$ as test function in the weak formulation \eqref{BasicWeakFormMain}, then we have conservation of mass, in the sense that $\frac{1}{|\Omega|} \int_{\Omega} \rho(x, t) \dx = \frac{1}{|\Omega|} \int_{\Omega} \rho(x, 0) \dx = \frac{1}{|\Omega|} \int_{\Omega} \rho_{0}(x) \dx = 1$ for all $t \geq 0$, and because $\rho(x, t)$ is supposed to be nonnegative for all $t\geq 0$, also 
\begin{equation*}
\|\rho(t)\|_{L^{1}(\Omega)} = \|\rho_{0}\|_{L^{1}(\Omega)} \quad \textrm{for } t \geq 0. 
\end{equation*}

A further important property that holds true for our equation is that 
$\|\rho(t)\|_{L^\infty(\Omega)} \leq \|\rho_0\|_{L^\infty(\Omega)}$. This fact is shown in \cite{chen2022analysis} for the integral definition of the fractional Laplacian on $\mathbb{R}^{d}$. We will follow a similar strategy for our argument and show that the property holds true when $\Omega$ is a bounded Lipschitz domain, and using the spectral definition of the fractional Neumann Laplacian. We will also show later on that a similar statement holds true for the discretization of our problem.

\begin{lemma} \label{InfNormBoundClass}
Let $\rho$ be a nonnegative smooth solution of the initial-boundary-value problem \eqref{MainFullProb},  let $F \in C^{2}([0, \infty))$ be a convex function, and suppose that $F(\rho_{0}) \in L^{1}(\Omega)$. Then, 
\[ \frac{\mathrm{d}}{\mathrm{d}t} \int_{\Omega} F(\rho)\dx \leq 0. \]
\begin{proof}
Assume first that $F''$ is bounded. Then, $F'(\rho) - F'(0)$ with $\rho \in H^{1}(\Omega)$ is an admissible test function in the weak formulation (\ref{BasicWeakForm}). Thus we additionally have that 
\begin{align}
& \frac{\mathrm{d}}{\mathrm{d}t} \int_{\Omega} F(\rho) \dx = - \int_{\Omega} \nabla \rho \cdot \nabla F'(\rho) \dx - \int_{\Omega} \rho F''(\rho) \nabla \rho \cdot \nabla (-\Delta_{\mathrm{N}})^{-\s} \rho^{\ast} \dx \nonumber \\
\Leftrightarrow \quad & \frac{\mathrm{d}}{\mathrm{d}t} \int_{\Omega} F(\rho) \dx = - \int_{\Omega} F''(\rho) |\nabla \rho|^{2} \dx - \int_{\Omega} \rho F''(\rho) \nabla \rho \cdot \nabla (-\Delta_{\mathrm{N}})^{-\s} \rho^{\ast} \dx \nonumber \\
\Leftrightarrow \quad & \frac{\mathrm{d}}{\mathrm{d}t} \int_{\Omega} F(\rho) \dx = - \underbrace{\int_{\Omega} F''(\rho) |\nabla \rho|^{2} \dx}_{\circled{1}} - \underbrace{\int_{\Omega} H(\rho) (-\Delta_{\mathrm{N}})^{1-\s} \rho^{\ast} \dx}_{\circled{2}}, \label{IntPartsBound}
\end{align}
where $H(s) := \int_{0}^{s} t F''(t) \dt$. Notice that $H$ is a nondecreasing function by definition. 

Term $\circled{1}$ is clearly nonnegative. For term $\circled{2}$ we notice that letting $1-\s = \frac{\beta}{2}$ we can use the explicit representation of the spectral fractional Laplacian
(cf., for example, eq. (1.11) in \cite{Stinga-Torrea})
\[ (-\Delta_{\mathrm{N}})^{\frac{\beta}{2}} u(x) = \frac{1}{\Gamma(-\beta/2)} \int_{0}^{\infty} (\mathrm{e}^{\tau \Delta_{\mathcal{N}}} u(x) - u(x)) \frac{1}{\tau^{1+\frac{\beta}{2}}} \text{d}\tau, \]
where $\mathrm{e}^{t \Delta_{\mathcal{N}}} u(x)$ is the solution of the problem
\begin{equation} \left\{ \begin{array}{ll} \frac{\partial w}{\partial \tau} = \Delta w & \textrm{in } \Omega \times (0, \infty), \\
w(x, 0) = u(x) & \textrm{in } \Omega, \\ 
\partial_{n}w(x, \tau) = 0 &\textrm{on } \partial \Omega \times [0, \infty). \end{array} \right. 
\label{HeatNeuSemiGroup} \end{equation}
We note that by the definition (\ref{ProjZeroAv}) of $\rho^\ast$ we have that $\mathrm{e}^{\tau \Delta_{\mathcal{N}} }\rho^{\ast} = \mathrm{e}^{\tau \Delta_{\mathcal{N}} }\rho - \bar{\rho}$. Thus,
\begin{align*}
\circled{2} &= \int_{\Omega} H(\rho(x)) \frac{1}{\Gamma(-\beta/2)} \int_{0}^{\infty} (\mathrm{e}^{\tau \Delta_{\mathcal{N}}} \rho^{\ast}(x) - \rho^{\ast}(x)) \frac{1}{\tau^{1+\frac{\beta}{2}}} \,\dtau \dx \\
&= \int_{\Omega} H(\rho(x)) \frac{1}{\Gamma(-\beta/2)} \int_{0}^{\infty} (\mathrm{e}^{\tau\Delta_{\mathcal{N}}} \rho(x) - \rho(x)) \frac{1}{\tau^{1+\frac{\beta}{2}}} \,\dtau \dx \\
&= \frac{1}{\Gamma(-\beta/2)} \int_{0}^{\infty} \bigg( \int_{\Omega} H(\rho(x)) \mathrm{e}^{\tau\Delta_{\mathcal{N}}} \rho(x) - \int_{\Omega} H(\rho(x)) \rho(x)\bigg) \frac{1}{\tau^{1+\frac{\beta}{2}}}\, \dx \dtau\\
&= \frac{1}{\Gamma(-\beta/2)} \int_{0}^{\infty} \big(  (\mathrm{e}^{\tau\Delta_{\mathcal{N}}}\rho, H(\rho) ) - ( \rho, H(\rho)) \big) \frac{1}{\tau^{1+\frac{\beta}{2}}}\, \dtau.
\end{align*}

Let $W_{\tau}(x, y)$ be the distributional heat kernel for $-\Delta$ subject to a homogeneous Neumann boundary condition and $\tau > 0$. We have that $W_{\tau} \geq 0$ for all $\tau>0$ and, for $u, v \in L^{2}(\Omega)$,
\[ (\mathrm{e}^{\tau \Delta_{\mathcal{N}}}u, v ) = \int_{\Omega} \int_{\Omega} W_{\tau}(x, y) u(y) v(x) \dy\dx. \]
Therefore, by plugging the heat kernel in the expression above, we deduce that 
\begin{align*}
\circled{2} &= \frac{1}{\Gamma(-\beta/2)} \int_{0}^{\infty} \int_{\Omega} \bigg( \int_{\Omega} W_{\tau}(x, y) \rho(x) H(\rho(y)) \dx - \rho(y) H(\rho(y)) \bigg) \dy \,\frac{1}{\tau^{1+\frac{\beta}{2}}} \,\dtau \\
&= \frac{1}{\Gamma(-\beta/2)} \int_{0}^{\infty} \int_{\Omega} \bigg( \int_{\Omega} W_{\tau}(x, y)(\rho(x) - \rho(y)) H(\rho(y))\dx \\
&\quad + \rho(y) H(\rho(y)) \bigg( \int_{\Omega} W_{\tau}(x, y)\dx - 1 \bigg) \bigg)\dy\, \frac{1}{\tau^{1+\frac{\beta}{2}}} \,\dtau \\
&= \frac{1}{\Gamma(-\beta/2)} \int_{0}^{\infty} \int_{\Omega} \int_{\Omega} W_{\tau}(x, y) (\rho(x) - \rho(y)) H(\rho(y)) \dy \dx \, \frac{1}{\tau^{1+\frac{\beta}{2}}} \,\dtau \\
&\quad + \frac{1}{\Gamma(-\beta/2)} \int_{0}^{\infty} \int_{\Omega} \rho(y) H(\rho(y))(\mathrm{e}^{\tau\Delta_{\mathcal{N}}}1(y) - 1) \dy \,\frac{1}{\tau^{1+\frac{\beta}{2}}} \,\dtau.
\end{align*}
By exchanging $x$ and $y$ and using the symmetry of the heat kernel we have also 
\begin{align*}
\circled{2} = &-  \frac{1}{\Gamma(-\beta/2)} \int_{0}^{\infty} \int_{\Omega} \int_{\Omega} W_{\tau}(x, y) (\rho(x) - \rho(y)) H(\rho(x)) \dy\dx \, \frac{1}{\tau^{1+\frac{\beta}{2}}} \,\dtau \\
& +\frac{1}{\Gamma(-\beta/2)} \int_{0}^{\infty} \int_{\Omega} \rho(y) H(\rho(y))(\mathrm{e}^{\tau\Delta_{\mathcal{N}}}1(y) - 1) \dy\, \frac{1}{\tau^{1+\frac{\beta}{2}}} \,\dtau. 
\end{align*}
By summing the two expressions for $\circled{2}$ we then find that  
\begin{align*}
\circled{2} = &-\frac{1}{2\Gamma(-\beta/2)} \Bigg( \int_{0}^{\infty} \int_{\Omega} \int_{\Omega} W_{\tau}(x, y) (\rho(x) - \rho(y))(H(\rho(x)) - H(\rho(y))) \dx \dy \,\frac{1}{\tau^{1+\frac{\beta}{2}}} \dtau \\
& + \int_{0}^{\infty} \int_{\Omega} \rho(y) H(\rho(y))(1- \mathrm{e}^{\tau\Delta_{\mathcal{N}}}1(y)) \dy \,\frac{1}{\tau^{1+\frac{\beta}{2}}} \,\dtau \Bigg).
\end{align*}
Thanks to the definition of $\Delta_{\mathcal{N}}$ stated in \eqref{HeatNeuSemiGroup} we have that $\mathrm{e}^{\tau\Delta_{\mathcal{N}}}1(y) \equiv  1$, which results in a convenient simplification of the expression for term $\circled{2}$ appearing in the last equality. Hence, 
from the monotonicity of $H$ we finally deduce that $\circled{2} \geq 0$. In the general case (i.e., when $F''$ is not bounded) we define $F_{k}(u) := F(0) + F'(0)u + \int_{0}^{u} \int_{0}^{v} \min(F''(w), k) \, \text{d}w \, \text{d}v$ for $k>0$; then, $F''_{k}(u)$ is bounded and the result follows for $F$ replaced by $F_{k}$ and taking the limit $k \to \infty$ using monotone convergence.
\end{proof}
\end{lemma}
The previous result allows us to prove the following lemma, which provides a bound on the $L^{\infty}(\Omega)$ norm of $\rho(t)$ in terms of $\|\rho_0\|_{L^\infty(\Omega)}$.
\begin{lemma} \label{InfNormDecayClass}
Let $\rho$ be a nonnegative smooth solution of the initial-boundary-value problem \eqref{MainFullProb}. The following inequality holds for all $t\geq 0$: 
\[ \|\rho(t)\|_{L^{\infty}(\Omega)} \leq \|\rho_{0}\|_{L^{\infty}(\Omega)}. \]
\begin{proof}
By the previous lemma we have that
\[ \sup_{t>0} \int_{\Omega} F(\rho(t)) \dx \leq \int_{\Omega} F(\rho_{0}) \dx. \]
We choose a nonnegative convex function $F \in C^{2}([0, \infty))$ such that $F(u) = 0$ for $u \leq \|\rho_{0}\|_{L^{\infty}(\Omega)}$ and $F(u) > 0$ for $u > \|\rho_{0}\|_{L^{\infty}(\Omega)}$. Then, 
\[ 0 \leq \int_{\Omega} F(\rho(t)) \dx \leq \int_{\Omega} F(\rho_{0}) \dx = 0 \quad \textrm{for } t \geq 0, \] 
and thus $0 \leq \rho(x, t)  \leq \|\rho_{0}\|_{L^{\infty}(\Omega)}$ for $t\geq 0$.
\end{proof}
\end{lemma}

Lemma \ref{InfNormDecayClass} together with the nonnegativity of $\rho$ and mass conservation implies that $\|\rho(t)\|_{L^{p}(\Omega)}$ is bounded for all $t\geq 0$ and $1 \leq p \leq \infty$.

We continue by investigating the energy structure of our problem (\ref{BasicWeakForm}),  (\ref{MyFracPois1}). Let us define the nonnegative convex function $G \in C([0,\infty))$ by
\begin{align}\label{eq:G}
G(s) := s (\log s - 1) + 1\quad \mbox{for $s>0$}\quad \mbox{and}\quad  G(0):=1.
\end{align}

We note that $G(1)=0$ and $\lim_{s \rightarrow +\infty} G(s)/s = + \infty$. We shall assume for the moment that $\rho(x, t) \geq  0$  for all $x \in \Omega$ and all $t>0$ and its spatial integral is normalized as $\frac{1}{|\Omega|} \int_{\Omega} \rho(x, t) \dx = 1$, and consider, as discussed in the Introduction, the free energy functional 
\begin{equation} \label{Energy1} 
E(\rho) := \int_{\Omega} G(\rho)\dx - \frac{1}{2} \int_{\Omega} c \rho \dx = \int_{\Omega} G(\rho)\dx - \frac{1}{2} \int_{\Omega} c \rho^{\ast} \dx, 
\end{equation}
where the second equality follows thanks to the fact that $c$ has zero integral average on $\Omega$ (cf. Theorem \ref{FracPoiExUniq}). Both terms appearing in the free energy functional \eqref{Energy1} are nonnegative. Indeed, by applying Jensen's inequality with the uniform probability measure $\mathrm{d}\mu(x):=\frac{1}{|\Omega|} \dx$ on $\Omega$ to the convex function $G$ yields that 
\[  
\int_\Omega G(\rho) \dd \mu(x) \geq  G\left(\int_\Omega \rho(x) \dd \mu(x)\right) = 0,
\]
thanks to the normalization of $\rho$ and the conservation of  mass. On the other hand, we can express the second term in \eqref{Energy1} in terms of the eigenvalues and eigenfunctions of the Neumann Laplacian as follows:
$$
- \frac{1}{2} \int_{\Omega} c \rho^{\ast} \dx = \frac{1}{2} \sum_{k=1}^\infty \lambda_k^s \rho_k^2\geq 0,
$$
using that $-(-\Delta_{\mathrm{N}})^s c = \rho^\ast$ on $\Omega$, in conjunction with the definition \eqref{FracLapNeu} 
of the fractional Neumann Laplacian and Parseval's identity in $L^2_\ast(\Omega)$.
Therefore, the free energy functional $E(\rho)$ is a suitable quantity for controlling the `distance' between $\rho$ and $\mu$, the uniform density on the domain $\Omega$; this will be rigorously proved in Theorem \ref{th4.7} stated in Section \ref{Sect:4}.
\begin{lemma} \label{FormalDecayFreeEnergy}
Let $\rho$ be a nonnegative smooth solution of the initial-boundary-value problem \eqref{MainFullProb}. The following free energy identity holds for all $t\geq 0$: 
\begin{equation}\label{DecayEnergy1}
    \frac{\mathrm{d}}{\mathrm{d}t} E(\rho) = - \int_{\Omega} \rho |\nabla (\log \rho - c)|^{2}  \dx \leq 0\, . 
\end{equation}
\begin{proof}
By (formally) differentiating $E(\rho)$ with respect to $t$ we have
\begin{align*}
\frac{\mathrm{d}}{\mathrm{d}t} E(\rho) & = \int_{\Omega} \rho_{t} G'(\rho) \dx - \frac{1}{2} \int_{\Omega} c_{t} \rho^{\ast}  \dx - \frac{1}{2} \int_{\Omega} c \rho_{t}  \dx.
\end{align*}
Recalling that we have conservation of mass, $\frac{\mathrm{d}}{\mathrm{d}t} \int_{\Omega} \rho \dx=0$, and so $\partial_{t}\rho^{\ast} = \partial_{t} \rho$, it follows that 
\[  \int_{\Omega} c_{t} \rho^{\ast}  \dx = -
\int_{\Omega} ((-\Delta_{\mathrm{N}})^{-\s} \rho^{\ast}_{t}) \,\rho^{\ast}  \dx= - \int_{\Omega} \rho_{t} \,(-\Delta_{\mathrm{N}})^{-\s} \rho^{\ast}  \dx = \int_{\Omega} \rho_{t} c  \dx. \] 
As a consequence, using the weak formulation (\ref{BasicWeakForm}), we have 
\begin{align*}
\frac{\mathrm{d}}{\mathrm{d}t} E(\rho) & = \int_{\Omega} \rho_{t} \log \rho  \dx - \int_{\Omega} \rho_{t} c  \dx = \int_{\Omega} \rho_{t} (\log \rho - c)  \dx \nonumber \\ &= - \int_{\Omega}  \rho \nabla \log \rho \cdot \nabla(\log \rho - c) \dx + \int_{\Omega} \rho \nabla c \cdot \nabla (\log \rho - c) \dx \nonumber \\ &= - \int_{\Omega} \rho |\nabla (\log \rho - c)|^{2}  \dx \leq 0, 
\end{align*}
and therefore we have the desired decrease of the energy in time.
\end{proof}
\end{lemma}
Our goal is now to formulate a fully discrete numerical scheme based on spatial finite element discretization and implicit Euler time-stepping, which can reproduce all of these properties at the temporally semidiscrete-level upon passing to the limit with the spatial discretization parameter. We shall then pass to the limit with the time-step, while keeping the essential estimates intact in the limit, to establish the existence of a weak solution to the initial-boundary-value problem \eqref{MainFullProb},  and rigorously prove the decrease in time of the energy.

In order to build our numerical scheme for equation (\ref{MainFullProb}) that mimics the properties exhibited by $\rho$, because the assumption that $\rho \geq 0$ played an important role in the derivation of the previous energy estimate, we introduce for $0<\delta < 1< L < \infty$ the positive cut-off function $\beta_{\delta}^{L} \in C^{0,1}(\mathbb{R})$ in the weak formulation (\ref{BasicWeakForm}) via
\begin{equation*} 
\beta_{\delta}^{L}(s) := \left \{ \begin{array}{ll} \delta &  \mbox{if } s \leq \delta, \\ s &  \mbox{if } \delta < s < L, \\ L & \mbox{if }  s \geq L, \end{array} \right. 
\end{equation*}
prior to discretizing the problem. The two cut-off parameters, $\delta$ and $L$, will be present in the statement and the analysis of the fully discrete approximation of our model in Section \ref{Sect:3}. We will then take the limit as the lower cut-off parameter $\delta$ and the spatial discretization parameter tend to zero, and in Section \ref{Sect:4} we will make the semidiscrete-in-time approximation independent of the upper cut-off parameter $L$ by establishing a bound on the $L^{\infty}(\Omega)$ norm, which is uniform with respect to $L$ and analogous to the one in Lemma \ref{InfNormDecayClass}. 

We therefore introduce the following `truncated' weak formulation:
\begin{subequations} \label{RegWeakFormMain}
\begin{gather}
\text{Find } \rhodL  \in L^2(0,T;V) \text{ with $\frac{\partial \rho_{\delta,L}}{\partial t} \in L^2(0,T;V')$ such that } \nonumber \\
    \Big\langle \frac{ \partial \rhodL}{\partial t}, \phi  \Big\rangle = -\int_{\Omega} \nabla \rhodL \cdot \nabla \phi  \dx + \int_{\Omega} \beta_{\delta}^{L}(\rhodL)  \nabla c_{\delta, L} \cdot \nabla \phi  \dx \quad \textrm{for all } \phi \in V \textrm{ and a.e. $t \in (0,T]$}, \label{RegWeakForm}
\end{gather}
subject to the initial condition $\rho_{\delta, L}(\cdot, 0) = \rho_{0}$, where 
\begin{equation}\label{RegWeakFormb}
- (-\Delta_{\mathrm{N}})^{s} c_{\delta, L} = \rhodL^{\ast} \quad \text{ in } \Omega, 
\end{equation}
\end{subequations}
and, as before, $V = H^1(\Omega)$. In this case, similarly as before, we assume that a solution $\rhodL$ to our regularized problem exists and that it has all the regularity properties that we need for the formal computations that follow.

The structure of the problem (\ref{RegWeakForm}) is such that it still preserves decay of the energy in time. To confirm this, we introduce a regularized version $G^L_\delta \in C^{2,1}(\mathbb{R})$ of the function $G(\cdot)$, defined, for $0 < \delta < 1<L < \infty$, by
\begin{equation*} G_{\delta}^{L}(s) := \left \{ \begin{array}{ll} \frac{s^{2} - \delta^{2}}{2 \delta} + (\log \delta - 1)s + 1 & \mbox{if }  s \leq \delta, \\
G(s) & \mbox{if }  \delta < s < L, \\
\frac{s^{2} - L^{2}}{2L} + (\log L - 1) s + 1 & \mbox{if }  s \geq L.
\end{array} \right.
\end{equation*}
Clearly, 
\[ (G_{\delta}^{L})'(s) = \left \{ \begin{array}{ll} \frac{s}{\delta} + \log \delta - 1 & \mbox{if }  s \leq \delta, \\ \log s & \mbox{if }  \delta < s < L, \\ \frac{s}{L} + \log L - 1 & \mbox{if }  s \geq L, \end{array} \right. \qquad (G_{\delta}^{L})''(s) = \left \{ \begin{array}{ll}  \frac{1}{\delta} & \mbox{if }   s \leq \delta, \\ \frac{1}{s} & \mbox{if }  \delta < s < L, \\ \frac{1}{L} & \mbox{if }  s \geq L ,\end{array} \right. \]
and therefore 
\[ \beta_{\delta}^{L}(s) (G_{\delta}^{L})''(s) = 1 \quad \textrm{for all } s \in \mathbb{R}.\]
In addition, for any sufficiently smooth function $\phi$, the following equality holds (see \cite{barrett2012finite}):
\begin{equation} \label{PropBas1} \beta_{\delta}^{L}(\phi) \nabla [(G_{\delta}^{L})'(\phi)] = \nabla \phi. \end{equation}
Moreover, we have that
\begin{equation} \label{ImpProp1}
\min \{ G_{\delta}^{L}(s), s(G_{\delta}^{L})'(s) \} \geq \left \{ \begin{array}{ll} \frac{s^{2}}{2 \delta} &\textrm{if } s \leq 0, \\ \frac{s^{2}}{4L} - C(L) & \textrm{if } s \geq 0, \end{array} \right.
\end{equation}
where $C(L)$ is a positive constant depending only on $L$. Finally we note that $(G_{\delta}^{L})''(s) \geq \frac{1}{L}$ for all $s \in \mathbb{R}$. 
The regularized version of the energy $E(\rho)$ is then defined by
\begin{equation*} 
E_{\delta}^{L}(\rhodL) := \int_{\Omega} G_{\delta}^{L}(\rhodL)  \dx - \frac{1}{2} \int_{\Omega} \rhodL^{\ast} c_{\delta, L}  \dx, 
\end{equation*}
and with similar computations as the ones in (\ref{DecayEnergy1}) it is possible to show that
\[ \frac{\mathrm{d}}{\mathrm{d}t} E_{\delta}^{L}(\rho_{\delta, L}) \leq 0. \]
\begin{remark}
    The relationship between the potential and the density expressed by the fractional Poisson equation $-(-\Delta_{\mathrm{N}})^{\s} \cdL = \rhodL^\ast$ will be crucial throughout and we want to highlight now one direct consequence of it, which is the nonpositivity of the $L^{2}$ inner product in space between the gradients of $\rhodL$ and $c_{\delta, L}$, namely 
    \begin{equation} \label{NegProdGradReg}
    \int_{\Omega} \nabla \rhodL \cdot \nabla \cdL \dx \leq 0.
\end{equation}
    Indeed, by using the definition of the Neumann Laplacian $-\Delta_{\mathrm{N}}$ and recalling that $c_{\delta,L} = (-\Delta_{\mathrm{N}})^{-\s} \rhodL^{\ast}$, we get
    \begin{align*}
        \int_{\Omega} \nabla \rhodL \cdot \nabla c_{\delta, L} \dx= \int_{\Omega} \nabla \rhodL^\ast \cdot \nabla c_{\delta, L} \dx = \int_{\Omega} \rho_{\delta, L}^{\ast}\, (-\Delta_{\mathrm{N}}) c_{\delta, L}  \dx
        = - \int_{\Omega} \rhodL^{\ast} (-\Delta_{\mathrm{N}})^{1-\s} \rhodL^{\ast} \dx \leq 0. 
    \end{align*}   
\end{remark}

We can also derive a formal energy bound. Using the weak formulation (\ref{RegWeakForm}) with $\phi = (G_{\delta}^{L})'(\rhodL)$ and property (\ref{PropBas1}), we have that
\[ \frac{\mathrm{d}}{\mathrm{d}t} E_{\delta}^{L}(\rhodL) =  - \int_{\Omega} (G_{\delta}^{L})''(\rhodL) \nabla \rhodL \cdot \nabla \rhodL \dx + 2 \int_{\Omega} \nabla \rhodL \cdot \nabla \cdL \dx - \int_{\Omega} \beta_{\delta}^{L}(\rhodL) \nabla \cdL \cdot \nabla \cdL \dx. \]
By integrating this equality in time from 0 to $t$ for some $t \in (0, T]$, where $T>0$, and using (\ref{ImpProp1}) with the assumption that the initial datum is $(\rhodL)_{0} =\rho_{0}$, we have that 
\begin{align*} \frac{1}{2\delta} \int_{\Omega} [\rhodL]_{-}^{2} \dx - \frac{1}{2} \int_{\Omega} \rhodL^{\ast} \cdL \dx &\leq C(L) - \int_{0}^{t} \int_{\Omega} (G_{\delta}^{L})''(\rho) |\nabla \rhodL|^{2} \dx \\ & \quad + 2 \int_{0}^{t} \int_{\Omega} \nabla \rhodL \cdot \nabla \cdL \dx - \int_{0}^{t}\int_{\Omega} \beta_{\delta}^{L}(\rhodL) |\nabla \cdL|^{2} \dx.   \end{align*}
The property (\ref{NegProdGradReg}) then implies that 
\[ \frac{1}{2\delta} \int_{\Omega} [\rhodL]_{-}^{2} \dx + \frac{1}{L} \int_{0}^{t} \int_{\Omega}  |\nabla \rhodL|^{2} \dx- \frac{1}{2} \int_{\Omega} \rhodL^{\ast} \cdL \dx + \int_{0}^{t}\int_{\Omega} \beta_{\delta}^{L}(\rhodL) |\nabla \cdL|^{2} \dx \leq C(L), \]
which means, taking the supremum over $t \in (0,T)$ and recalling the equality $-(-\Delta_{\mathrm{N}})^{\s} c_{\delta, L} = \rhodL^{\ast}$, that we have
\begin{align} \label{StabEst}
\begin{aligned}
\frac{1}{\delta} \underset{t \in (0, T)}{\sup} \int_{\Omega} [\rhodL]_{-}^{2} \dx &+ \int_{0}^{T}  \int_{\Omega}  |\nabla \rhodL|^{2} \dx \\
&+ \underset{t \in (0, T)}{\sup} \int_{\Omega} (-\Delta_{\mathrm{N}})^{\s} \cdL \, \cdL \dx + \int_{0}^{T}\int_{\Omega} \beta_{\delta}^{L}(\rhodL) |\nabla \cdL|^{2} \dx
\leq C(L). \end{aligned}\end{align}
Hence, we arrive at the following bound, which is uniform with respect to the parameter $\delta$:
\[ \frac{1}{\delta} \underset{t \in (0, T)}{\sup} \|[\rhodL]_{-}\|_{L^{2}(\Omega)}^{2} + \int_{0}^{T}  \int_{\Omega}  |\nabla \rhodL|^{2} \dx + \underset{t \in (0, T)}{\sup} \|\cdL\|_{\mathbb{H}^{\s}_{\ast}(\Omega)}^{2} + \int_{0}^{T}\int_{\Omega} \beta_{\delta}^{L}(\rhodL) |\nabla \cdL|^{2} \dx \leq C(L). \]
A discrete counterpart of this inequality will be used later on to prove nonnegativity of the numerical solution by passing to the limit $\delta\to 0_+$.

We can also consider the weak formulation (\ref{RegWeakForm}) with $\phi = \rhodL$. We then have that
\[ \int_{\Omega} (\rhodL)_{t}\, \rhodL \dx = -\int_{\Omega} |\nabla \rhodL|^{2} \dx + \int_{\Omega} \beta_{\delta}^{L}(\rhodL) \nabla \cdL \cdot \nabla \rhodL \dx, \]
and the application of the Cauchy--Schwarz inequality and Young's inequality yields, for any $\eta >0$, that
\[
\frac{1}{2} \frac{\mathrm{d}}{\mathrm{d}t} \int_{\Omega} \rhodL^{2} \dx \leq -\int_{\Omega} |\nabla \rhodL|^{2} \dx + \frac{\eta}{2} \int_{\Omega} \beta_{\delta}^{L}(\rhodL) |\nabla \rhodL|^{2} \dx + \frac{1}{2 \eta} \int_{\Omega} \beta_{\delta}^{L}(\rhodL) |\nabla \cdL|^{2} \dx.
\]
Hence, by taking $\eta$ sufficiently small and recalling (\ref{StabEst}), we get the following overall bound on $\rhodL$:
\begin{equation}\label{StabEstPlus}
\underset{t \in (0, T)}{\sup} \int_{\Omega} \rhodL^{2} \dx + \frac{1}{\delta} \underset{t \in (0, T)}{\sup} \int_{\Omega} [\rhodL]_{-}^{2} \dx + \int_{0}^{T}  \int_{\Omega}  |\nabla \rhodL|^{2} \dx \leq C(L).
\end{equation}

Motivated by these formal calculations, we shall design a finite element method whose solution mimics the energy law of (\ref{RegWeakForm}) and the formal bounds (\ref{StabEst}) and (\ref{StabEstPlus}) at a discrete level, and show that, as the spatial and temporal discretization parameters tend to zero, a subsequence of this sequence of approximations converges to a bounded and nonnegative weak solution of our original problem in a suitable sense. In doing so, we shall confine ourselves to the physically relevant cases of $d=2$ and $d=3$ space dimensions. The arguments that follow can be therefore viewed as a constructive proof of the existence of a nonnegative, energy-dissipative, weak solution to the initial-boundary-value problem for the fractional porous medium equation under consideration, based on the Neumann Laplacian.  \section{Finite element approximation}\label{Sect:3}

We now define the fully discrete scheme for the numerical approximation of the problem (\ref{RegWeakForm}),  (\ref{RegWeakFormb}). Our approach is to couple the strategy used in \cite{barrett2011existence, barrett2011finite, barrett2012finite} in terms of the spatial finite element scheme with recent results concerning numerical approximation of the fractional Laplacian \cite{bonito2018numerical, bonito2021approximation, bonito2015numerical, bonito2017numerical}.  The overall aim is to construct a finite element approximation of the  problem (\ref{MainFullProb}) in its weak formulation (\ref{RegWeakForm}), (\ref{RegWeakFormb}), which mimics the estimate (\ref{StabEst}) at a discrete level, and to show that a subsequence of this sequence of approximations converges to a bounded and nonnegative weak solution, first as the spatial discretization parameter $h$ and the lower cut-off parameter $\delta$ tend to zero, and then as the time-step $\Delta t$ tends to zero. It will also transpire that it suffices to take a fixed value of the upper cut-off parameter $L$ such that $L>\|\rho_0\|_{L^\infty(\Omega)}$; passage to the limit $L \to +\infty$ will therefore not be required. We begin by introducing the spatial finite element approximation. 

\subsection{Finite Element Approximation}

We assume that $\Omega$ is a bounded open polygonal domain in $\mathbb{R}^{2}$ or a bounded open Lipschitz polyhedral domain in 
$\mathbb{R}^{3}$. With $h>0$ being our spatial discretization parameter, let $\mathcal{T}_{h}= \{ K_{n} \}_{n=1}^{M_{h}}$ be a quasi-uniform and shape-regular family of triangulations of $\overline\Omega$, where $K_{n}$, $n=1,\ldots, M_h$, are closed simplices with mutually disjoint interiors, such that $\overline{\Omega} = \cup_{n=1}^{M_{h}} K_{n}$; assume further that the triangulation is weakly acute. Let $V_{h}$ be the linear space of continuous piecewise affine functions defined on this triangulation, i.e., 
\begin{equation*}
    V_{h} := \big\{ v_{h} \in C(\overline{\Omega}) \textrm{ such that } v_{h}\big|_{K} \in \mathbb{P}^{1} \textrm{ for all } K \in \mathcal{T}_{h} \big\},
\end{equation*}
and let $N_{h}:=\mbox{dim}(V_h)$. Moreover, let $\pi_{h}: C^{0}(\overline{\Omega}) \to V_{h}$ be the interpolation operator into the linear space $V_{h}$ based on nodal evaluation, i.e., such that, for each $\phi \in C(\overline\Omega)$,
\begin{equation*}
    \pi_{h} \phi(P_{j}) = \phi(P_{j}), \quad j = 1, \dots, N_{h},
\end{equation*} 
where $\{ P_{j} \}_{ j=1}^{N_{h}}$ are the nodes (vertices) of $\mathcal{T}_{h}$ contained in $\overline \Omega$.

In order to construct a discrete version of the fractional spectral Neumann Laplacian we first consider a finite-dimensional approximation, $-\Delta_{h}$, of the standard Neumann Laplacian $-\Delta_{\mathrm{N}}$. Let us consider the bilinear form
\begin{equation} \label{BiForm} a: H^{1}_{\ast}(\Omega) \times H^{1}_{\ast}(\Omega) \rightarrow \mathbb{R} \quad \mbox{where} \quad a(u, v) := \int_{\Omega} \nabla u \cdot \nabla v \dx, \end{equation}
which is symmetric, bounded and coercive; we can restrict $a(\cdot,\cdot)$ to
$(V_{h} \cap L^{2}_{\ast}(\Omega))\times(V_{h} \cap L^{2}_{\ast}(\Omega))$. Clearly, 
$V_{h} \cap L^{2}_{\ast}(\Omega) ~(= V_h \cap H^{1}_{\ast}(\Omega))$ is an $(N_h-1)$-dimensional linear subspace of $H^{1}_{\ast}(\Omega)$.
The bilinear form (\ref{BiForm}) gives rise in $V_{h} \cap L^{2}_{\ast}(\Omega)$ to a discrete orthonormal basis of eigenfunctions that we denote by $\{ \varphi_{k}^{h}\}_{k=1}^{N_{h}-1}$ with corresponding positive eigenvalues $\lambda_{1}^{h}, \dots, \lambda_{N_{h}-1}^{h}$ such that 
\begin{equation*}
    a(\varphi_{k}^{h}, v_{h}) = \lambda_{k}^{h} (\varphi_{k}^{h}, v_{h}) \quad \textrm{for all } v_{h} \in V_{h} \cap L^{2}_{\ast}(\Omega),\quad k=1,\ldots,N_h-1.
\end{equation*}
Using these, we can now propose a similar expression to (\ref{FracLapNeu}) and introduce a finite-dimensional counterpart of the fractional spectral Neumann Laplacian, $(-\Delta_{\mathrm{N}})^{\s}$; namely, for $u_{h} \in V_{h}\cap L^{2}_{\ast}(\Omega)$ and $s \in [0,1]$, we define
\begin{equation} \label{FinElemFracLapDef}
    (-\Delta_{h})^{\s} u_{h} := \sum_{k=1}^{N_{h}-1} (\lambda_{k}^{h})^{\s} u_{k}^{h} \varphi_{k}^{h}  \quad \textrm{with } u_{k}^{h} := \int_{\Omega} u_{h} \varphi_{k}^{h} \dx,\quad k=1,\ldots,N_h-1.
\end{equation} 
\begin{remark}
We notice that this definition is consistent with the definition of $-\Delta_{h}$, in the sense that if $s=1$ then
\[ (-\Delta_{h})u_{h} = \sum_{k=1}^{N_{h}-1} \lambda_{k}^{h} u_{k}^{h} \varphi_{k}^{h}  \quad \textrm{with } u_{k}^{h} = \int_{\Omega} u_{h} \varphi_{k}^{h} \dx,\quad k=1,\ldots,N_h-1, \]
and
\[ ((-\Delta_{h}) u_{h}, v_{h}) = \sum_{k=1}^{N_{h}-1} u_{k}^{h} \lambda_{k}^{h} (\varphi_{k}^{h}, v_{h}) = \sum_{k = 1}^{N_{h}-1} u_{k}^{h} a( \varphi_{k}^{h}, v_{h}) = a(u_{h}, v_{h}). \]

We shall also define the finite element approximation, $T_{h}$, of the inverse spectral Neumann Laplacian $T = (-\Delta_{\mathrm{N}})^{-1}$, which with each $F \in H^{-1}_{\ast}(\Omega)$ (uniquely) associates $T_{h}(F) := u_{h} \in V_{h}\cap L^2_*(\Omega)$ such that 
\[ a(u_{h}, v_{h}) = \langle F, v_{h} \rangle, \quad \text{for all } v_{h} \in V_{h} \cap L^2_*(\Omega). \]
\end{remark}

\begin{remark} \label{RemarkSpecFracCompute}
    In our forthcoming paper \cite{carrillosuli2024} we will show that an accurate approximation to the
    discrete fractional spectral Neumann Laplacian $(-\Delta_h)^s$ can be efficiently computed without explicit knowledge of the eigenpairs $(\lambda_k^h, \varphi_k^h)$, $k=1,\ldots, N_h-1$; only a positive lower bound on the smallest eigenvalue $\lambda_1^h$ and an upper bound on the largest eigenvalue $\lambda_{N_h-1}^{h}$ are required.  Full details and extensive numerical simulations will be given in \cite{carrillosuli2024}.
    \end{remark}

We shall require in the analysis of our fully discrete numerical method a discrete version of the property (\ref{PropBas1}). In order to ensure that this holds, we define a diagonal matrix $\Theta_{\delta}^{L}(\phi_{h}) \in \mathbb{R}^{d \times d}$ in the following way: for an element $K \in \mathcal{T}_h$, let $\{ P_{i} \}_{i = 0}^{d}$ be the vertices of the simplex $K$; then, for $j = 1, \dots, d$ and $x \in K^\circ$ (the interior of $K$), we define
\begin{subequations} \label{ThetaMatrix}
\begin{equation}
\widetilde{\Theta}_{\delta}^{L}(\phi_{h})_{jj}(x) := \left \{ \begin{aligned}& \frac{\phi_{h}(P_{j}) - \phi_{h}(P_{0})}{(G_{\delta}^{L})'(\phi_{h}(P_{j})) - (G_{\delta}^{L})'(\phi_{h}(P_{0}))} & \textrm{if } \phi_{h}(P_{j}) \neq \phi_{h}(P_{0}), \\
&\frac{1}{(G_{\delta}^{L})''(\phi_{h}(P_{j}))} = \beta_{\delta}^{L}(\phi_{h}(P_{j})) & \textrm{if } \phi_{h}(P_{j}) = \phi_{0}(P_{j}). \end{aligned} \right.
\end{equation}
\begin{figure}[H]
\centering
\begin{tikzpicture}[vect/.style={->,
             shorten >=0pt,>=latex'}]
\tkzDefPoint(-0.3, -0.5){A}       
\tkzDefPoint(2.6, 0.2){B}     
\tkzDefPoint(1.2, 1.8){C}
\tkzDefPoint(1.2, 0.8){O}
\tkzDefPoint(-0.8, 1.6){P1}
\tkzDefPoint(3.4, 1.4){P2}
\tkzDefPoint(1.2, -1.5){P3}
\tkzDrawSegment(A, B)
\tkzDrawSegment[line width = 1.3](B, C)
\tkzDrawPoints(B,C)
\tkzDrawSegment(C, A)
\tkzDefPointBy[symmetry= center C](O)
\tkzGetPoint{C1}
\tkzDefPointBy[symmetry= center B](O)
\tkzGetPoint{B1}
\tkzDefPointBy[symmetry= center A](O)
\tkzGetPoint{A1}
\tkzDrawSegment[dashed, add = 0 and 0.2](A, P1)
\tkzDrawSegment[dashed, add = 0 and 0.3](C, P1)
\tkzDrawSegment[dashed, add = 0 and 0.3](B, P2)
\tkzDrawSegment[dashed, add = 0 and 0.2](C, P2)
\tkzDrawSegment[dashed](B, B1)
\tkzDrawSegment[dashed](C, C1)
\tkzDrawSegment[dashed](A, A1)
\tkzDrawSegment[dashed, add = 0 and 0.2](A, P3)
\tkzDrawSegment[dashed, add = 0 and 0.3](B, P3)
\tkzLabelPoint[below = 0.1](B){$P_{0}$}
\tkzLabelPoint[above right](C){$P_{j}$}
\end{tikzpicture}
\caption{Triangulation and the computation of the matrix $\widetilde{\Theta}_{\delta}^{L}$}
\end{figure} 
The matrix $\widetilde{\Theta}_{\delta}^{L}$ is designed to approximate $\beta_{\delta}^{L} I$, where $\beta_{\delta}^{L}$ is the positive cut-off function defined above, and its elements are simply difference quotients approximating the inverse of the derivative of $(G_{\delta}^{L})'$, which is $(\beta_{\delta}^{L})^{-1}$. 

Now,  let $\widehat{K}$ be the reference simplex and $\widehat{x} \in \widehat{K} \mapsto P_{0} + B_{K} \widehat{x}=x \in K$ the affine function that maps $\widehat{K}$ onto $K$, and for $x \in K^\circ$ define
\begin{equation}
\Theta_{\delta}^{L}(\phi_{h})(x) := (B_{K}^{\mathrm{T}})^{-1} \,\widetilde{\Theta}_{\delta}^{L}(\phi_{h})(x) \,B_{K}^{\mathrm{T}}.
\end{equation}
\end{subequations}
With this definition in place we have that 
\begin{equation} \label{PropFE1} \Theta_{\delta}^{L}(\phi_{h})(x)\, \nabla [\pi_{h}(G_{\delta}^{L})'(\phi_{h})(x)] = \nabla \phi_{h} (x) \quad \textrm{for all } x \in K^\circ. \end{equation}

The proof of the following lemma is a simple extension of the proof of Lemma 2.1 in \cite{barrett_convergence_2003} and is therefore omitted.

\begin{lemma} \label{ThetaContLemma}
    Let $\| \cdot \|$ denote the spectral norm on $\mathbb{R}^{d \times d}$. For any $\delta \in (0,1)$ and $L>1$ the function $\Theta_{\delta}^{L}: V_{h} \to \mathbb{R}^{d \times d}$ is continuous, and it satisfies 
    \begin{equation*}
        \delta \,  \xi^{\mathrm{T}} \xi \leq \xi^{\mathrm{T}} \Theta_{\delta}^{L}(\phi_{h}) \xi \leq L \, \xi^{\mathrm{T}} \xi \quad \text{for all } \xi \in \mathbb{R}^{d},\, \phi_{h} \in V_{h}.    
    \end{equation*}
    Furthermore, for all $\phi_{h}^{(1)}, \phi_{h}^{(2)} \in V_{h}$ and all $K \in \mathcal{T}_{h}$, one has that 
        \begin{equation*}
            \big\| \big(\Theta_{\delta}^{L}\big(\phi_{h}^{(1)}\big) - \Theta_{\delta}^{L}\big(\phi_{h}^{(2)}\big)\big)(x) \big\| \leq \frac{L}{\delta} \max_{j=1, \dots, d} (|\phi_{h}^{(1)}(P_{j}) - \phi_{h}^{(2)}(P_{j})| + |\phi_{h}^{(1)}(P_{0}) - \phi_{h}^{(2)}(P_{0})|)\quad \mbox{for all $x \in K^\circ$}.
        \end{equation*}
\end{lemma}

Moreover, it is possible to prove the following lemma (the proof can be found in \cite{barrett2011finite}), which clarifies the sense in which the matrix function $\Theta_{\delta}^{L}(\cdot)$ approximates the cut-off function $\beta_{\delta}^{L}(\cdot) I$, as $\delta$ and $h$ tend to zero.
\begin{lemma} \label{LemmaMatrix}
    For each $K \in \mathcal{T}_{h}$ and for all $\phi_{h} \in V_{h}$ we have that 
    \begin{equation*}
        \int_{K} |\Theta_{\delta}^{L}(\phi_{h}) - \beta^{L}_{\delta}(\phi_{h}) I |^{2} \dx \leq C \bigg(\delta^{2} + h^{2} \int_{K} |\nabla \phi_{h}|^{2} \dx + \int_{K} \pi_{h}([\phi_{h}]_{-}^{2}) \dx \bigg).
    \end{equation*}
\end{lemma}

The following lemma, together with its corollary, were also proved in \cite{barrett2011finite}; they will be crucial in our argument to show convergence of a subsequence as $h$ and $\delta$ tend to zero and we therefore report them here.

\begin{lemma} \label{LemLip}
Let $g \in C^{0,1}(\mathbb{R})$ be monotonically increasing with Lipschitz constant $g_{\mathrm{Lip}}$; then, for all $K \in \mathcal{T}_{h}$ and for all $\phi_{h} \in V_{h}$,
\begin{equation*}
    \int_{K} |\nabla( \pi_{h}(g(\phi_{h}))|^{2} \dx \leq g_{\mathrm{Lip}} \int_{K} \nabla \phi_{h} \cdot \nabla \pi_{h}(g(\phi_{h})) \dx.
\end{equation*}
\begin{cor} \label{CorLip}
Let $g$ be defined and strictly monotonically increasing on $\mathbb{R}$, such that $g^{-1}$, the inverse function of $g$, is Lipschitz continuous on $\mathbb{R}$, with Lipschitz constant $(g^{-1})_{\mathrm{Lip}}$; then, for all $K \in \mathcal{T}_{h}$ and for all $\phi_{h} \in V_{h}$,
\begin{equation*}
    \int_{K} |\nabla \phi_{h}|^{2} \dx \leq (g^{-1})_{\mathrm{Lip}} \int_{K} \nabla \phi_{h} \cdot \nabla \pi_{h}(g(\phi_{h})) \dx.
\end{equation*}
\end{cor}
\end{lemma}
We further note that, for each $K \in \mathcal{T}_{h}$ with vertices $\{ P_{j} \}_{j=0}^{d}$, by virtue of Jensen's inequality we have
\begin{equation} \label{PropFE2}
 (\pi_{h}(\phi)(x))^{2} = \bigg( \sum_{j = 0}^{d} \phi(P_{j}) \phi_{j}(x) \bigg)^{2} \leq \sum_{j = 0}^{d} (\phi(P_{j}))^{2} \phi_{j}(x) = \pi_{h}(\phi^{2})(x) \quad \mbox{for all $x \in K$ and all $\phi \in C(\overline{K})$}. 
 \end{equation}

In order to formulate the discretization of the problem under consideration with respect to the temporal variable, we let $T>0$ denote the final time and we introduce a temporal mesh given by $0 =: t_{0} < t_{1} < \dots < t_{N} := T$, which for the ease of exposition we choose to be uniform, with spacing $\Delta t$; i.e.,  $t_{n} - t_{n-1} = \Delta t$ for all $n=1, \dots, N$. 

Given the initial datum $\rho_{0} \in L^\infty(\Omega)$ such that $\rho_0 \geq 0$ and $\frac{1}{|\Omega|} \int_{\Omega} \rho_{0} \dx = 1$, we first choose $\rho^0_h \in V_{h}$ as the unique solution of the following problem:
\begin{equation} \label{InValDiscProb}
\int_{\Omega} \pi_h(\rho^0_h \phi_{h}) \dx+ \Delta t \int_{\Omega} \nabla \rho^0_h \cdot \nabla \phi_{h} \dx= \int_{\Omega} \rho_{0} \phi_{h} \dx \quad \text{for all } \phi_{h} \in V_{h}.
\end{equation}
The existence of a unique solution $\rho^0_h \in V_h$ to \eqref{InValDiscProb} is a consequence of the fact that \eqref{InValDiscProb} is equivalent to a system of linear algebraic equations with a symmetric positive definite matrix. By taking $\phi_h = \rho^0_h$ in \eqref{InValDiscProb}, it follows using \eqref{PropFE2} and Young's inequality that 
\begin{equation} \label{BoundInVal}
    \int_{\Omega} (\rho^0_h)^{2} \dx + 2\Delta t \int_{\Omega} |\nabla \rho^0_h|^{2} \dx \leq  \int_{\Omega} \rho_{0}^{2} \dx \leq |\Omega|\, \|\rho_0\|^2_{L^\infty(\Omega)}.
\end{equation}
Furthermore, by taking $\phi_h\equiv 1 \in V_h$ in  \eqref{InValDiscProb} and noting that $\pi_h \rho_h^0 = \rho_h^0$, it follows that
\[\frac{1}{|\Omega|} \int_\Omega \rho_h^0(x) \dx = \frac{1}{|\Omega|}\int_\Omega \rho_0(x) \dx = 1.\]
Next, we shall prove that $0 \leq \rho^0_h(x) \leq \|\rho_0\|_{L^\infty(\Omega)}$ for all $x \in \overline\Omega$. To this end, we consider the function $\varphi_p$, defined so that its first derivative is a regularized version of the function $s \in \mathbb{R} \mapsto [s]_{-} \in \mathbb{R}_{\leq 0}$. More precisely, for $p > 0$, we consider
\[ \varphi_p(s) := \begin{cases} 
      \frac{s^2}{2} + \frac{s}{2p} + \frac{1}{6p^2} & \mbox{if } s\leq - \frac{1}{p}, \\
      -\frac{ps^3}{6} & \mbox{if }  -\frac{1}{p}\leq s\leq 0, \\
      0 & \mbox{if }  s \geq 0. 
   \end{cases}
\]
The first and second derivative of $\varphi_p$ are then given by, respectively,
\[
\varphi_p'(s) = \begin{cases} 
      s + \frac{1}{2p}  & \mbox{if }  s\leq - \frac{1}{p}, \\
      -\frac{ps^2}{2} & \mbox{if }  -\frac{1}{p}\leq s\leq 0, \\
      0 & \mbox{if }  s \geq 0; 
   \end{cases}
   \qquad 
\varphi_p''(s) = \begin{cases} 
      1  & \mbox{if }  s\leq - \frac{1}{p}, \\
      -ps & \mbox{if }  -\frac{1}{p}\leq s\leq 0, \\
      0 & \mbox{if }  s \geq 0. 
   \end{cases}
\]
Clearly, $\varphi_p \in C^{2,1}(\mathbb{R})$ and $\varphi_p$ is a nonnegative convex function.  Furthermore, $\varphi_p'(s) \leq 0$ and $s \hspace{0.2mm}\varphi_p'(s) \geq 0$ for all $s \in \mathbb{R}$. 

We then take $\phi_h = \pi_h(\varphi_p'(\rho^0_h))$ as test function in \eqref{InValDiscProb} and apply Lemma \ref{LemLip} to infer that with this choice of $\phi_h$ the second term on the left-hand side of 
\eqref{InValDiscProb} is nonnegative, and therefore
\[ \int_\Omega \pi_h(\rho^0_h \,\pi_h(\varphi_p'(\rho^0_h))) \dx  \leq \int_\Omega \rho_0\, \pi_h(\varphi_p'(\rho^0_h))\dx \leq 0, \]
because, by hypothesis, $\rho_0(x) \geq 0$ for a.e.~$x \in \Omega$ while $\pi_h(\varphi_p'(\rho^0_h(x)))\leq 0$ for all $x \in \overline\Omega$ since $\varphi_p'(\rho^0_h(x)) \leq 0$ for all $x \in \overline\Omega$.
On the other hand, 
\[ 0 \leq \int_\Omega \pi_h(\rho^0_h \,\varphi_p'(\rho^0_h)) \dx = \int_\Omega \pi_h(\rho^0_h\, \pi_h(\varphi_p'(\rho^0_h))) \dx \leq 0,\]
because $\rho^0_h(x) \,\varphi_p'(\rho^0_h(x)) \geq 0$ for all $x \in \overline\Omega$. Hence, 
\[ \int_\Omega \pi_h(\rho^0_h \,\varphi_p'(\rho^0_h))\dx = 0.\]
Passing to the limit $p \rightarrow +\infty$, this implies that
\[ \int_\Omega \pi_h(\rho^0_h \, [\rho^0_h]_{-}) \dx = 0.\]
As $\pi_h(\rho^0_h \, [\rho^0_h]_{-}) = \pi_h(([\rho^0_h]_{-})^2)$ and the latter is a nonnegative function on $\overline\Omega$, we deduce that $\pi_h(([\rho^0_h]_{-})^2)$ is identically zero on $\overline\Omega$; equivalently, $([\rho^0_h]_{-})^2(P_j)=0$ at all vertices $P_j$, $j=1,\ldots, N_h$, of the triangulation. Therefore, $\rho^0_h(P_j) \geq 0$ for all $j=1,\ldots,N_h$, and consequently, because $\rho^0_h$ is piecewise affine, $\rho^0_h(x) \geq 0$ for all $x \in \overline\Omega$.

One can use an analogous argument to show that $\rho^0_h \leq \|\rho_0\|_{L^\infty(\Omega)}$ on $\overline\Omega$. Indeed, if we define $A:=\|\rho_0\|_{L^\infty(\Omega)}$ and the functions $\hat{\rho}^0_h:= A- \rho^0_h$ and $\hat{\rho}_0:= A- \rho_0$, then, because $\pi_h \phi_h = \phi_h$ for all $\phi_h \in V_h$, it follows from 
\eqref{InValDiscProb} that
\begin{equation*} 
\int_{\Omega} \pi_h(\hat{\rho}^0_h \phi_{h}) \dx+ \Delta t \int_{\Omega} \nabla \hat{\rho}^0_h\cdot \nabla \phi_{h} \dx= \int_{\Omega} \hat\rho_{0} \phi_{h} \dx \quad \text{for all } \phi_{h} \in V_{h}.
\end{equation*}
As $\hat\rho_0 \geq 0$ a.e.~on $\Omega$, an identical argument to the one above implies that $\hat{\rho}^0_h \geq 0$ on $\overline\Omega$, whereby $\rho^0_h \leq A = \|\rho_0\|_{L^\infty(\Omega)}$, as required. To summarize, we have shown that
\begin{align}\label{inibound}
0 \leq {\rho}^0_h(x) \leq \|\rho_0\|_{L^\infty(\Omega)}\quad \mbox{for all $x \in \overline\Omega$}\quad \mbox{and}\quad \frac{1}{|\Omega|}\int_\Omega \rho_h^0(x) \dx = 1.
\end{align}

With $\rho^0_h$ thus defined, the fully-discrete scheme for the problem under consideration is then the following: 
\begin{subequations} \label{FullDiscWeakFormMain}
\begin{gather}
\textrm{ Let $\rho^0_{h,\delta,L}:=\rho^0_h \in V_h$.}
\textrm{ For } n = 1, \dots, N, \textrm{ given } \rhodLh^{n-1} \in V_{h}, \textrm{ find } \rhodLh^{n} \in V_
{h} \textrm{ such that } \nonumber \\
\label{FullDiscWeakForm} \int_{\Omega} \pi_h\bigg(\frac{\rhodLh^{n} - \rhodLh^{n-1}}{\Delta t} \phi_{h} \bigg)\dx = - \int_{\Omega} \nabla \rhodLh^{n} \cdot \nabla \phi_{h} \dx + \int_{\Omega} \Theta_{\delta}^{L}(\rhodLh^{n}) \nabla \cdLh^{n} \cdot \nabla \phi_{h} \dx \quad \textrm{for all } \phi_{h} \in V_{h},
\end{gather}
 where $\cdLh^{n} \in V_h \cap L^2_{\ast}(\Omega)$ satisfies
\begin{equation} \label{SpaceDiscFracPois}
- (-\Delta_{h})^{\s} \cdLh^{n} = (\rhodLh^{n})^{\ast}.
\end{equation}
\end{subequations}

The time-stepping scheme is chosen to be implicit and we stress that the numerical approximation relies on the potential $\cdLh$ and the fact that it is repulsive in the sense that \eqref{NegProdGradReg} holds. This will be crucial in deriving global-in-time existence, uniform boundedness and convergence of the sequence of numerical solutions. 

Since our first analytical step is to pass to the limit as $\delta, h \to 0_{+}$, we now introduce the definition of the semidiscrete-in-time approximation of our problem, to which the fully discrete scheme \eqref{FullDiscWeakFormMain} will be shown to converge in this limit. We first need to assign to $\rho_{0}$, for a fixed value of $\Delta t>0$, a certain smoothed initial datum $\rho^{0} = \rho^0({\Delta t}) \in V=H^{1}(\Omega)$, that is the solution to the following problem: 
\begin{equation} \label{InValSemiDiscProb}
    \int_{\Omega} \rho^{0} \phi \, \dx + \Delta t \int_{\Omega} \nabla \rho^{0} \cdot \nabla \phi \dx= \int_{\Omega} \rho_{0} \phi \, \dx \quad \text{for all } \phi \in V.
\end{equation}
It follows by an application of the Lax--Milgram lemma that \eqref{InValSemiDiscProb} has a unique solution $\rho^0= \rho^0(\Delta t) \in V=H^1(\Omega)$. By taking $\phi \equiv 1$ we find that $\frac{1}{|\Omega|} \int_\Omega \rho^0 \dx = \frac{1}{|\Omega|} \int_\Omega \rho_0 \dx = 1 $. Furthermore, one can show that  the following four assertions hold true:
\begin{equation*}
    \rho^{0} \to \rho_{0} \quad \text{weakly in } L^{2}(\Omega) \quad \text{as} \quad  \Delta t \to 0_{+};
\end{equation*} 
thanks to the assumed nonnegativity of $\rho_0$ also $\rho^0 \geq 0$ a.e.~on $\Omega$; the function $\rho^0$ satisfies
\begin{equation} \label{IneqPropIn}
    \int_{\Omega} G(\rho^{0}) \, \dx + 4 \Delta t \int_{\Omega} \Big|\nabla \sqrt{\rho^{0}}\Big|^{2} \, \dx \leq \int_{\Omega} G(\rho_{0})\dx,
\end{equation}
where $G\in C([0,\infty))$ is the nonnegative convex function defined in \eqref{eq:G}; and 
$\rho^{0} \in L^{\infty}(\Omega) \cap V$. The full details of the proofs of these assertions are contained in the Appendix of \cite{barrett2012finite} and in Section 6 of \cite{barrett2011existence}.

With this definition of $\rho^0=\rho^0(\Delta t)$ the semidiscrete-in-time approximation of our problem is then formulated as follows:
\begin{subequations} \label{WeakFormDiscTimeMain}
 \begin{gather}
\textrm{Let $\rho^0_L:= \rho^0 = \rho^0(\Delta t)$. For } n = 1, \dots, N, \textrm{ given } \rho_{L}^{n-1} \in V, \textrm{ find } \rho_{L}^{n} \in V \textrm{ such that } \nonumber \\
 \int_{\Omega} \frac{\rho_{L}^{n} - \rho_{L}^{n-1}}{\Delta t} \phi \dx = - \int_{\Omega} \nabla \rho_{L}^{n} \cdot \nabla \phi \dx+ \int_{\Omega} \beta^{L}(\rho_{L}^{n}) \nabla c_{L}^{n} \cdot \nabla \phi \dx \quad \textrm{for all } \phi \in V, \label{WeakFormDisctime}
\end{gather}
subject to the initial condition $\rhoL^{0} = \rho^{0}$, where
\begin{equation} \label{FinElemFracLap} - (-\Delta_{\mathrm{N}})^{\s} c_{L}^{n} = (\rho_{L}^{n})^{\ast}\quad \mbox{in $\Omega$}. \end{equation}
\end{subequations}

\subsection{Convergence analysis of the discretized fractional porous medium model}

Our aim in the rest of this section is to pass to the limit in (\ref{FullDiscWeakForm}), \eqref{SpaceDiscFracPois} as $\delta$ and $h$ tend to zero, to infer the semidiscrete-in-time approximation (\ref{WeakFormDisctime}), \eqref{FinElemFracLap} as the limiting problem. The analysis of the semidiscrete-in-time approximation will then be the focus of the next section, where we shall pass to the limit $\Delta t \rightarrow 0_+$ in  (\ref{WeakFormDisctime}), \eqref{FinElemFracLap}, and will also show that the sequence of approximations becomes independent of $L$ once $L>\|\rho_0\|_{L^\infty(\Omega)}$, which then makes passage to the limit $L \rightarrow +\infty$ redundant. The first result we prove is an a priori bound on the solution $\rhodLh^{n}$ of (\ref{FullDiscWeakForm}), assuming that it exists, with the existence of a solution to be proved thereafter. 
\begin{lemma} \label{AprioriBoundDisc}
For $n=1, \dots, N$ a solution $\rhodLh^{n} \in V_{h}$ of  \eqref{FullDiscWeakForm}, if it exists, satisfies the following bound:
\begin{equation} \label{APrioriBound1}
    \int_{\Omega} \pi_{h} (G_{\delta}^{L}(\rhodLh^{n})) \dx + \Delta t \int_{\Omega} \nabla \rhodLh^{n} \cdot \nabla \pi_{h}((G_{\delta}^{L})'(\rhodLh^{n}))) \dx \leq \int_{\Omega} \pi_{h} (G_{\delta}^{L}(\rhodLh^{n-1})) \dx.
\end{equation}
\begin{proof}
By taking $\phi_{h} = \pi_{h}((G_{\delta}^{L})'(\rho_{h}^{n}))$ as test function in (\ref{FullDiscWeakForm}), we have that
\begin{align*}
\int_{\Omega} \pi_h\bigg(\frac{\rhodLh^{n} - \rhodLh^{n-1}}{\Delta t} \pi_{h}((G_{\delta}^{L})'(\rhodLh^{n})) \bigg)\dx = &- \int_{\Omega} \nabla \rhodLh^{n} \cdot \nabla \pi_{h}((G_{\delta}^{L})'(\rhodLh^{n})) \dx \\
& + \int_{\Omega} \Theta_{\delta}^{L}(\rhodLh^{n}) \nabla \cdLh^{n} \cdot \nabla \pi_{h}((G_{\delta}^{L})'(\rhodLh^{n})) \dx.
\end{align*}
Hence, by (\ref{PropFE1}), it follows that
\[ 
\int_{\Omega}\! \pi_h \bigg(\!\frac{\rhodLh^{n} - \rhodLh^{n-1}}{\Delta t} \pi_{h}((G_{\delta}^{L})'(\rhodLh^{n})) \!\bigg)\!\dx = - \int_{\Omega} \nabla \rhodLh^{n} \cdot \nabla \pi_{h}((G_{\delta}^{L})'(\rhodLh^{n}))\! \dx\, + \int_{\Omega} \nabla \cdLh^{n} \cdot \nabla \rhodLh^{n}\! \dx.
\]
Next, we bound the left-hand side of this equality from below by using the convexity of $G_{\delta}^{L}$, which implies that
\begin{align*}
   \pi_h \bigg( (\rhodLh^{n} - \rhodLh^{n-1})\,\pi_{h}((G_{\delta}^{L})'(\rhodLh^{n})) \bigg) &=\pi_h \bigg( (\rhodLh^{n} - \rhodLh^{n-1})\,(G_{\delta}^{L})'(\rhodLh^{n}) \bigg)\\&\hspace{-7mm}\geq \pi_{h}\bigg(G_{\delta}^{L}(\rhodLh^{n}) - G_{\delta}^{L}(\rhodLh^{n-1}) \bigg)= \pi_{h}(G_{\delta}^{L}(\rhodLh^{n})) - \pi_{h} (G_{\delta}^{L}(\rhodLh^{n-1})),
\end{align*}
and therefore 
\begin{equation} \label{DiscBound1} \frac{1}{\Delta t}  \int_{\Omega} (\pi_{h}(G_{\delta}^{L}(\rhodLh)) - \pi_{h} (G_{\delta}^{L}(\rhodLh^{n-1}))) \dx + \int_{\Omega} \nabla \rhodLh^{n} \cdot \nabla \pi_{h}((G_{\delta}^{L})'(\rhodLh^{n})) \dx \leq \int_{\Omega} \nabla \cdLh^{n} \cdot \nabla \rhodLh^{n} \dx. \end{equation}

We observe that, in the case of a repulsive potential under consideration here, the right-hand side of inequality \eqref{DiscBound1} contributes `constructively' to the desired bound, in the sense that it satisfies
\begin{equation} \label{NegativeProd}
    \int_{\Omega} \nabla \cdLh^{n} \cdot \nabla \rhodLh^{n} \dx \leq 0.
\end{equation}
Indeed, thanks to the definition of the finite element approximation of the fractional Laplacian stated in equation (\ref{FinElemFracLapDef}), we have that
\begin{align*}
   \int_{\Omega} \nabla c_{h, \delta, L}^{n} \cdot \nabla \rho_{h, \delta, L}^{n} \dx &= \int_{\Omega} \nabla c_{h, \delta, L}^{n} \cdot \nabla (\rho_{h, \delta, L}^{n})^{\ast} \dx = a(c_{h, \delta, L}^{n}, (\rho_{h, \delta, L}^{n})^{\ast}) = a(c_{h, \delta, L}^{n}, -(-\Delta_{h})^{\s} c_{h, \delta, L}^{n}) \\
   &= - \sum_{k=1}^{N_{h}-1} (\lambda_{k}^{h})^{\s} (c_{h, \delta, L}^{n})_{k} a(c_{h, \delta, L}^{n}, \phi_{k}^{h}) =  - \sum_{k=1}^{N_{h}-1} (\lambda_{k}^{h})^{1+\s} (c_{h, \delta, L}^{n})_{k}^{2} \leq 0.
\end{align*}
This then completes the proof of the lemma.
\end{proof}
\end{lemma}

We shall now prove the existence of a solution at a given time $t_{n}$ in order to be able to use the bound (\ref{APrioriBound1}). We will do so by using Brouwer's fixed point theorem. 

\begin{lemma} \label{DiscExistLemma1}
For any $\Delta t > 0$, given $\rhodLh^{n-1} \in V_{h}$, there exists at least one solution $(\rhodLh^{n}, c_{h, \delta, L}^{n}) \in V_{h} \times (V_{h} \cap L^{2}_{\ast}(\Omega))$ to \eqref{FullDiscWeakFormMain}.
\begin{proof}
We begin by equipping $V_h$ with the inner product $((\cdot,\cdot))$, defined by
\[ ((\psi_h, \varphi_h)):= \int_\Omega \pi_h (\psi_h(x) \phi_h(x)) \dx,\quad \phi_h, \varphi_h \in V_h.\]
Next, we define the function $\mathcal{H}: V_{h} \to V_{h}$ such that, for any $\rho_{h} \in V_{h}$,
\[ ((\mathcal{H}(\rho_{h}), \phi_{h})) := \int_{\Omega} \pi_h\bigg(\frac{\rho_{h} - \rhodLh^{n-1}}{\Delta t} \phi_{h} \bigg)\dx + \int_{\Omega} \nabla \rho_{h} \cdot \nabla \phi_{h}\dx - \int_{\Omega} \Theta_{\delta}^{L}(\rho_{h}) \nabla c_{h} \cdot \nabla \phi_{h} \dx \quad \textrm{for all } \phi_{h} \in V_{h},  \] 
where $c_h \in V_h \cap L^2_\ast(\Omega)$ satisfies $-(-\Delta_{h})^{\s} c_{h} = \rho_{h}^{\ast}$. The construction of $\Theta_{\delta}^{L}$, together with  Lemma \ref{ThetaContLemma}, ensures that $\mathcal{H}$ is a continuous mapping. Note that if a solution $\rhodLh^{n}$ to (\ref{FullDiscWeakForm}) exists then it is a zero of $\mathcal{H}$, namely
\[ ((\mathcal{H}(\rhodLh^{n}), \phi_{h})) = 0 \quad \textrm{for all } \phi_{h} \in V_{h}. \] 

We shall therefore prove that $\mathcal{H}$ has a zero. 
For contradiction, let us assume that $\mathcal{H}$ has no zero for any $\gamma \in \mathbb{R}_{>0}$ in the ball 
$B_{\gamma} = \{ \psi_{h} \in V_{h} \textrm{ such that } \|\psi_{h}\|_{L^{2}(\Omega)} \leq \gamma \}$. We define the function $\mathcal{E}_{\gamma}: B_{\gamma} \to B_{\gamma}$ by
\[ \mathcal{E}_{\gamma}(\psi_{h}) = - \gamma \frac{\mathcal{H}(\psi_{h})}{\| \mathcal{H}(\psi_{h})\|_{L^{2}(\Omega)}}. \]

By Brouwer's fixed point theorem $\mathcal{E}_{\gamma}$ has at least one fixed point $\rho_{h}^{\gamma}$ in $B_{\gamma}$, and this means that $\|\rho_{h}^{\gamma}\|_{L^{2}(\Omega)} = \|\mathcal{E}_{\gamma}(\rho_{h}^{\gamma})\|_{L^{2}(\Omega)} = \gamma$. Let $c_{h}^{\gamma}$ be the solution of the equation $-(-\Delta_{\mathrm{N}})^{s} c_{h}^{\gamma} = (\rho_{h}^{\gamma})^{\ast}$ in $\Omega$. Now let us choose $\phi_{h} = \pi_{h}((G_{\delta}^{L})'(\rho_{h}^{\gamma}))$ as test function and notice that, by (\ref{ImpProp1}) and inequality \eqref{PropFE2}, together with the trivial equality $\pi_h \rho^\gamma_h = \rho^\gamma_h$, we have that
\begin{align*}
 ((\mathcal{H}(\rho_{h}^{\gamma}), \pi_{h}((G_{\delta}^{L})'(\rho_{h}^{\gamma})))) 
  &= -\frac{\|\mathcal{H}(\rho_{h}^{\gamma})\|_{L^{2}(\Omega)}}{\gamma} \bigg( \int_{\Omega} \pi_h(\rho_{h}^{\gamma} \, \pi_{h}((G_{\delta}^{L})'(\rho_{h}^{\gamma}))) \dx \bigg) \\
  &\hspace{-7mm}= -\frac{\|\mathcal{H}(\rho_{h}^{\gamma})\|_{L^{2}(\Omega)}}{\gamma} \bigg( \int_{\Omega} \pi_h(\rho_{h}^{\gamma}\,  (G_{\delta}^{L})'(\rho_{h}^{\gamma})) \dx \bigg) \\
  &\hspace{-7mm}\leq -  \frac{\|\mathcal{H}(\rho_{h}^{\gamma})\|_{L^{2}(\Omega)}}{\gamma} \bigg( \int_{\Omega} \bigg( \frac{(\rho_{h}^{\gamma})^{2}}{4L} - C(L) \bigg) \dx \bigg) =- \frac{\|\mathcal{H}(\rho_{h}^{\gamma})\|_{L^{2}(\Omega)}}{\gamma} \Big( \frac{\gamma^{2}}{4L} -  C(L) |\Omega| \Big),
 \end{align*}
 and thus for $\gamma  > [4L\, C(L)\, |\Omega|]^{1/2}$ we also have that 
 \begin{equation*}   
 ((\mathcal{H}(\rho_{h}^{\gamma}), \pi_{h}(G_{\delta}^{L})'(\rho_{h}^{\gamma}))) < 0. \end{equation*}
 
 On the other hand, we can proceed in a similar way as in Lemma \ref{AprioriBoundDisc} to deduce that
 \begin{align}
 ((\mathcal{H}(\rho_{h}^{\gamma}), \pi_{h}(G_{\delta}^{L})'(\rho_{h}^{\gamma}))) &= \frac{1}{\Delta t} \int_{\Omega} \pi_h\big((\rho_{h}^{\gamma} - \rhodLh^{n-1}) \pi_{h}( (G_{\delta}^{L})'(\rho_{h}^{\gamma}))\big) \dx + \int_{\Omega} \nabla \rho_{h}^{\gamma} \cdot \nabla \pi_{h} ((G_{\delta}^{L})'(\rho_{h}^{\gamma})) \dx \nonumber \\
 &  \quad -\int_{\Omega} \Theta_{\delta}^{L}(\rho_{h}^{\gamma}) \nabla c_{h}^{\gamma} \cdot \nabla \pi_{h} ((G_{\delta}^{L})'(\rho_{h}^{\gamma})) \dx \nonumber \\
 &\geq \underbrace{\frac{1}{\Delta t} \int_{\Omega} \pi_h\big((\rho_{h}^{\gamma} - \rhodLh^{n-1}) \pi_{h} ((G_{\delta}^{L})'(\rho_{h}^{\gamma})) \big)\dx}_{\circled{1}}  \underbrace{- \int_{\Omega} \nabla c_{h}^{\gamma} \cdot \nabla \rho_{h}^{\gamma} \dx.}_{\circled{2}} \label{BrouwerEqRef1}
 \end{align}
For the first term in (\ref{BrouwerEqRef1}) we can use the convexity of $G^L_\delta$ to deduce, as in the proof of the previous lemma, that 
 \begin{align*}
 \circled{1} \geq \frac{1}{\Delta t} \int_{\Omega} \Big( \pi_{h} (G_{\delta}^{L}(\rho_{h}^{\gamma})) - \pi_{h} (G_{\delta}^{L}(\rho_{h}^{n-1})) \Big) \dx \geq \frac{1}{\Delta t} \Big( \frac{\gamma^{2}}{4L} - C(L) |\Omega| \Big) - \frac{1}{\Delta t} \int_{\Omega} \pi_{h}( G_{\delta}^{L}(\rho_{h}^{n-1}) )\dx.
 \end{align*}
For the second term in (\ref{BrouwerEqRef1}) we can use the definition of $c_{h}^{\gamma}$ involving the  discrete fractional Neumann Laplacian to deduce an inequality analogous to (\ref{NegativeProd}), which then implies that $\circled{2} \geq 0$.
Therefore,
\begin{align*}
 ((\mathcal{H}(\rho_{h}^{\gamma}), \pi_{h}((G_{\delta}^{L})'(\rho_{h}^{\gamma})))) &\geq \frac{1}{\Delta t} \Big( \frac{\gamma^{2}}{4L} - C(L) |\Omega| \Big) - \frac{1}{\Delta t} \int_{\Omega} \pi_{h} (G_{\delta}^{L}(\rho_{h}^{n-1})) \dx.
\end{align*}
Thus, for $\gamma  > [4L\, (C(L)\, |\Omega| + \| \pi_{h} (G_{\delta}^{L}(\rho_{h}^{n-1}))\|_{L^{1}(\Omega)})]^{1/2}$ we have 
\[  ((\mathcal{H}(\rho_{h}^{\gamma}), \pi_{h}((G_{\delta}^{L})'(\rho_{h}^{\gamma})))) > 0, \]
and hence we have arrived at a contradiction; therefore the function $\mathcal{H}$ has a zero, which means that at time $t_{n}$ there exists a solution $\rhodLh^{n}\in V_h$ to the fully discrete scheme (\ref{FullDiscWeakForm}). Having shown the existence of $\rho_{h, \delta, L}^{n} \in V_{h}$, the existence of an associated (unique) $c_{h, \delta, L}^{n} \in V_{h} \cap L^{2}_{\ast}(\Omega)$ follows from (\ref{SpaceDiscFracPois}).
\end{proof}
\end{lemma}
We need a series of additional results before proceeding to the proof of convergence as $\delta$ and $h$ tend to zero. 
\begin{lemma}\label{lem:uniformdeltah}
A solution $\{\rhodLh^{n}\}_{n=1}^{N} \subset V_h$ of \eqref{FullDiscWeakForm} satisfies the following bound:
\begin{equation} \label{APrioriBound2}
    \underset{n = 1, \dots, N}{\max} \int_{\Omega} (\rhodLh^{n})^{2} \dx + \frac{1}{\delta}\underset{n = 1, \dots, N}{\max} \int_{\Omega} \pi_{h}( [\rhodLh^{n}]_{-}^{2}) \dx + \Delta t \sum_{n =1}^{N} \int_{\Omega} |\nabla \rhodLh^{n}|^{2} \dx \leq C(L). 
\end{equation}
Moreover, let $\cdLh^{n} \in V_h \cap L^2_\ast(\Omega)$ be the solution of \eqref{SpaceDiscFracPois}; then, the following bound holds: 
\begin{equation} \label{AprioriBound3}
    \Delta t \int_{\Omega} \Theta_{\delta}^{L}(\rhodLh^{n}) \nabla \cdLh^{n} \cdot \nabla \cdLh^{n} \dx \leq C(L) + \int_{\Omega} (\cdLh^{n})^{2} \dx.
\end{equation}
\begin{proof}
We sum the bounds (\ref{APrioriBound1}) over $n = 1, \dots, m$ to deduce that, for each $m\in \{1,\ldots,N\}$,
\begin{equation*}
     \int_{\Omega}\pi_{h}(G_{\delta}^{L}(\rhodLh^{m})) \dx + \Delta t \sum_{n = 1}^{m} \int_{\Omega} \nabla \rhodLh^{n} \cdot \nabla \pi_{h}((G_{\delta}^{L})'(\rhodLh^{n})) \dx \leq \int_{\Omega}\pi_{h}(G_{\delta}^{L}(\rho^0_h)) \dx.
\end{equation*}
Using Corollary \ref{CorLip} with $\phi_{h} = \rhodLh^{n}$ and $g = (G_{\delta}^{L})'$, noting that $g$ is strictly monotonically increasing and the inverse $g^{-1}$ is Lipschitz continuous with Lipschitz constant $L$, we then have that,  for each $m\in \{1,\ldots,N\}$,
\[  \int_{\Omega}\pi_{h}(G_{\delta}^{L}(\rhodLh^{m})) \dx + \frac{\Delta t}{L} \sum_{n = 1}^{m} \int_{\Omega} |\nabla \rhodLh^{n}|^{2} \dx \leq \int_{\Omega}\pi_{h}(G_{\delta}^{L}(\rho^0_h)) \dx. \]
Let $\Omega_+:= \{ x \in \overline{\Omega}\,:\, \rho^n_{h,\delta,L}\geq 0\}$ and let $\Omega_{-}:= \overline\Omega \setminus \Omega_+$.
We shall further bound the first term on the left-hand side of this inequality from below by noting that, thanks to (\ref{ImpProp1}) and (\ref{PropFE2}), we have 
\begin{align*}
  \int_{\Omega}\pi_{h}(G_{\delta}^{L}(\rhodLh^{n})) \dx &= \int_{\Omega_{+}}    \pi_{h}(G_{\delta}^{L}(\rhodLh^{n})) \dx + \int_{\Omega_{-}} \pi_{h}(G_{\delta}^{L}(\rhodLh^{n})) \dx \\
  &\hspace{-11mm}\geq \frac{1}{4L} \int_{\Omega_{+}} \pi_{h}((\rhodLh^{n})^{2})\dx - C(L)|\Omega_+| + \frac{1}{2\delta}\int_{\Omega_{-}} \pi_{h}((\rhodLh^{n})^{2}) \dx \\
  &\hspace{-11mm}= \frac{1}{4L} \int_{\Omega} \pi_{h}((\rhodLh^{n})^{2})\dx - C(L)|\Omega_+| + \frac{1}{2\delta}\int_{\Omega_{-}} \pi_{h}((\rhodLh^{n})^{2}) \dx - \frac{1}{4L} \int_{\Omega_{-}} \pi_{h}((\rhodLh^{n})^{2})\dx\\
  &\hspace{-11mm}\geq \frac{1}{4L} \int_{\Omega} (\pi_{h}\rhodLh^{n})^{2} \dx + \frac{1}{4\delta} \int_{\Omega} \pi_{h}([\rhodLh^{n}]_{-}^{2}) \dx - C(L)|\Omega| \\
  &\hspace{-11mm}= \frac{1}{4L} \int_{\Omega} (\rhodLh^{n})^{2} \dx + \frac{1}{4\delta} \int_{\Omega} \pi_{h}([\rhodLh^{n}]_{-}^{2}) \dx - C(L)|\Omega|, 
\end{align*}
where in the transition to the fourth line we have used that $\delta <L$ and that $[\rhodLh^{n}]_{-}^{2}\equiv 0$ on $\Omega_+$. It remains to note that, because for all $\delta$ and $L$ such that $0<\delta<1<L$ we have that
\[ 0 \leq G^L_\delta(s) \leq \max\bigg\{1 - \frac{\delta}{2}, \frac{s^2-L^2}{2L} + (\log L - 1)s + 1\bigg\}, \quad \mbox{for all $s \geq 0$},\]
it follows from \eqref{inibound} that $0 \leq G^L_\delta(\rho^0_{h,\delta, L}) \leq C$, where $C=C(L)$ is a positive constant, independent of $h$ and $\delta$, whereby the same is true of  $\int_\Omega \pi_h(G^L_\delta(\rho^0_{h,\delta, L})) \dx$. As $\Omega$ is a fixed domain, we shall absorb $|\Omega|$ into $C(L)$ and write $C(L)$ instead of $C(L)|\Omega|$. From this,  the inequality (\ref{APrioriBound2}) directly follows. 

To prove (\ref{AprioriBound3}) we take $\phi_{h} = -\cdLh^{n}$ in (\ref{FullDiscWeakForm}); we then have by (\ref{NegativeProd}) that
\[ \Delta t \int_{\Omega} \Theta_{\delta}^{L}(\rhodLh^{n}) \nabla \cdLh^{n} \cdot \nabla \cdLh^{n} \dx \leq \int_{\Omega} (\rhodLh^{n} - \rhodLh^{n-1}) \cdLh^{n} \dx. \]
Young's inequality implies that
\[ \Delta t \int_{\Omega} \Theta_{\delta}^{L}(\rhodLh^{n}) \nabla \cdLh^{n} \cdot \nabla \cdLh^{n} \dx \leq \frac{1}{2} \int_{\Omega}  (\rhodLh^{n})^{2} \dx + \frac{1}{2} \int_{\Omega}  (\rhodLh^{n-1})^{2} \dx +
\int_{\Omega} \cdLh^{2} \dx.\]
Finally, we bound the first two terms on the right-hand side of this inequality from above by using the bound on the first term appearing on the left-hand side of (\ref{APrioriBound2}). That completes the proof of inequality (\ref{AprioriBound3}).
\end{proof}
\end{lemma}
The final piece of information we need before taking the limit concerns the convergence of the discrete fractional Laplacian $(-\Delta_{h})^{s}$ and thus of the solution $\cdLh^{n}$ of the finite-dimensional fractional Poisson problem. We will rely on the results that were proved in \cite{bonito2021approximation, bonito2015numerical, bonito2017numerical}, where fractional powers of elliptic operators and their finite element approximation were studied in the setting we have put in place.

\begin{prop}[Convergence of the spatial discretization of the fractional Laplacian] \label{DiscFracLapConv}
Suppose that $s \in (0,1)$, $\sigma \in [0,1]$ and $f \in H^{2\sigma}_{\ast}(\Omega)$. Then, there exists a positive constant $C$ that is independent of $h$ and $\sigma$ such that 
\begin{equation*} 
    \|(-\Delta_{\mathrm{N}})^{-\s}f - (-\Delta_{h})^{-\s} \pi_{h}f\|_{L^{2}(\Omega)} \leq C \varepsilon(h) \|f\|_{H^{2\sigma}(\Omega)} \quad \textrm{for all } f \in H^{2\sigma}_{\ast}(\Omega),
\end{equation*}
where 
\begin{equation*} 
\varepsilon(h) = \left\{ \begin{array}{ll} 
h^{2} & \text{if } s + \sigma > 1, \\ \log(h^{-1}) h^{2(s+\sigma)} & \textrm{if } \s + \sigma \leq 1.
\end{array} \right.
\end{equation*}
\end{prop}

Now we are ready for passage to the limit with $\delta \rightarrow 0_+$ and $h\rightarrow 0_+$. Passage to the limit $\Delta t \rightarrow 0_+$ will be discussed in the next section. It will be also shown that by taking $L>\|\rho_0\|_{L^\infty(\Omega)}$ the sequence of numerical solutions becomes independent of $L$, so passage to the limit $L \rightarrow +\infty$ is not required.

\begin{theo} \label{convergenceFEMTheo}
The initial data $\{ \rho^0_h \}_{h > 0} $ defined in \eqref{InValDiscProb} are such that, for $\Delta t$ and $L$ fixed, as $h \to 0_{+}$ we have that
\begin{equation} \label{Conv0}
    \rho^0_h \to \rho^{0}_{L} = \rho^0 := \rho^0(\Delta t) \quad \text{strongly in } L^{2}(\Omega). 
\end{equation}
    Furthermore there exists a subsequence of $\{ \rhodLh^{n} \}_{\delta, h >0}$, a nonnegative function $\rho_{L}^{n} \in V=H^1(\Omega)$ and $c_{L}^{n} \in H^{1}_{\ast}(\Omega)$ such that the following convergence results hold, for each $n \in \{1,\ldots,N\}$, as $\delta, h \to 0_{+}$:
    \begin{subequations} \label{Convergence}
        \begin{alignat}{2}
            \label{Conv2}
            \nabla \rhodLh &\to \nabla \rho_{L}^{n} &&\quad\text{weakly in } L^{2}(\Omega;\mathbb{R}^d), \\
            \label{Conv3}
            \rhodLh^{n} &\to \rho_{L}^{n} &&\quad\text{strongly in } L^{2}(\Omega), \\
            \label{Conv4}
            \Theta_{\delta}^{L}(\rhodLh^{n}) &\to \beta^{L}(\rho^{n}_{L}) I&&\quad \text{strongly in } L^{2}(\Omega;\mathbb{R}^{d\times d}), \\
            \label{Conv5}
            \cdLh^{n} &\to c_{L}^{n} &&\quad \text{strongly in } L^{2}_\ast(\Omega;\mathbb{R}^d), \\
            \label{Conv6}
              \nabla \cdLh^{n} &\to  \nabla c_{L}^{n} &&\quad \text{weakly in }  L^{2}_\ast(\Omega;\mathbb{R}^d).
        \end{alignat}
    \end{subequations}
    Moreover $\{ \rho^{n}_{L} \}_{n=1, \dots, N}$ solves the problem \eqref{WeakFormDisctime},  $\rho^{n}_{L} \geq 0$ a.e.~on $\Omega$ for all $n=1, \dots, N$, and, given $\rho^{0}_L=\rho^0$ such that $\frac{1}{|\Omega|} \int_{\Omega} \rho^{0} \dx = 1$, one has $\frac{1}{|\Omega|} \int_{\Omega} \rho_{L}^{n} \dx = 1$ for all $n = 1, \dots, N$.
\begin{proof}
    It follows from (\ref{BoundInVal}) that there exists a subsequence of $\{ \rho^0_h \}_{h > 0}$ (which, for the sake of simplicity of notation, we shall continue to denote by $\{\rho^0_h\}_{h>0}$) that converges weakly in $L^{2}(\Omega)$ to some element $g \in L^{2}(\Omega)$, and $\nabla \rho^0_h$ converges weakly to some element $G \in L^{2}(\Omega;\mathbb{R}^d)$, as $h \to 0_{+}$. To show that $G = \nabla g$ notice that for any $\eta \in C^{\infty}_{0}(\Omega;\mathbb{R}^n)$ we have 
    \[ \int_{\Omega} G \cdot \eta \, \dx \leftarrow \int_{\Omega} \nabla \rho^0_h \cdot \eta \dx = - \int_{\Omega} \rho^0_h \text{div} \, \eta \, \dx \rightarrow - \int_{\Omega} g \,\text{div} \, \eta \, \dx, \]
    which means that $G$ is the distributional gradient of $g$, i.e., $G = \nabla g$ in $\mathcal{D}'(\Omega;\mathbb{R}^d)$, but since $G \in L^{2}(\Omega;\mathbb{R}^d)$ we then have that $G = \nabla g \in L^{2}(\Omega;\mathbb{R}^d)$. The compact embedding $H^{1}(\Omega) \compEmb L^{2}(\Omega)$ implies that (a subsequence of) $\{\rho^0_h\}_{h>0}$ convergences strongly in $L^2(\Omega)$ to $g \in L^2(\Omega)$ as $h \to 0_{+}$, and therefore we can pass to the limit in (\ref{InValDiscProb}) with $\phi_{h} = \pi_{h} \phi$ for any function $\phi \in C^{\infty}(\overline{\Omega})$ to get equation (\ref{InValSemiDiscProb}). This shows that $g = \rho^{0}$, and the fact that $C^{\infty}(\overline{\Omega})$ is dense in $V=H^1(\Omega)$ then implies  (\ref{Conv0}). 
    
    Now, let $n \in \{1,\ldots,N\}$. The weak convergence (\ref{Conv2}) and the strong convergence (\ref{Conv3}) are implied by the inequality  (\ref{APrioriBound2}) and the same considerations as the ones we exposed above in the case of $\nabla \rho_h^0$ and $\rho^h_0$. The nonnegativity of the limit function $\rho_{L}^{n}$ follows from (\ref{Conv3}) and the second bound in (\ref{APrioriBound2}).
    
     Combining the second and third bounds in (\ref{APrioriBound2}) and Lemma \ref{LemmaMatrix} we have the result (\ref{Conv4}).  
        
        The convergence result (\ref{Conv5}) can be proved in the following way. We have, by use of the stability inequality (\ref{StabFracNeu}) satisfied by the solution of the fractional Poisson equation, that
        \begin{align*}
&\|(-\Delta_{\mathrm{N}})^{-\s} (\rho_{L}^{n})^{\ast} - (-\Delta_{h})^{-\s} (\rhodLh^{n})^{\ast}\|_{L^{2}(\Omega)}\\ &\quad \leq
            \|(-\Delta_{\mathrm{N}})^{-\s}((\rho_{L}^{n})^{\ast} - (\rhodLh^{n})^{\ast})\|_{L^{2}(\Omega)}
            + \|((-\Delta_{\mathrm{N}})^{-\s} - (-\Delta_{h})^{-\s})(\rhodLh^{n})^{\ast}\|_{L^{2}(\Omega)} \\
            &\quad \leq C \|(\rho_{L}^{n})^{\ast} - (\rhodLh^{n})^{\ast}\|_{L^{2}(\Omega)} + \|((-\Delta_{\mathrm{N}})^{-\s} - (-\Delta_{h})^{-\s})(\rhodLh^{n})^{\ast}\|_{L^{2}(\Omega)},
        \end{align*}
        and by (\ref{Conv3}) and the convergence result stated in Proposition \ref{DiscFracLapConv} we have (\ref{Conv5}); here and henceforth $C$ signifies a generic positive constant, independent of $\delta$ and $h$. 
        Notice that this inequality, in conjunction with (\ref{StabFracNeu}) and the bound on the first term on the left-hand side of (\ref{APrioriBound2}),  makes the term on the left-hand side of (\ref{AprioriBound3}) uniformly bounded in $\delta$ and $h$.

        To prove (\ref{Conv6}) we need a uniform bound on the gradient of $\cdLh^{n}$; then, the result will follow from (\ref{Conv5}). Thanks to Theorem 4.8.12 in \cite{ern2004theory} we have that 
        \begin{equation} \label{InterpIneq} \|(I - \mathcal{S}_{h})v\|_{L^{2}(\Omega)} + h |(I - \mathcal{S}_{h}) v|_{H^{1}(\Omega)} \leq C h |v|_{H^{1}(\Omega)} \quad \text{for all } v \in H^{1}(\Omega), \end{equation}
        where $\mathcal{S}_{h}$ denotes the Scott--Zhang quasi-interpolation operator.
        For any $v_{h} \in V_{h}$, thanks to the assumed quasi-uniformity of the family of triangulations, we have also
        \begin{equation} \label{GradIneq} \|\nabla v_{h}\|_{L^{2}(\Omega)} \leq C h^{-1} \| v_{h}\|_{L^{2}(\Omega)}.\end{equation}
        Let us denote by $\tilde{c}_{h, \delta, L}^{n} \in H^1_\ast(\Omega)$ the solution of the boundary-value problem
        \[ -(-\Delta_{\mathrm{N}})^{\s} \tilde{c}_{h, \delta, L}^{n} = (\rhodLh^{n})^{\ast}\quad \mbox{in $\Omega$}.\]
        Then, the gradient of $\cdLh^{n}$ satisfies the following chain of inequalities: 
        \begin{align*}
            \|\nabla \cdLh^{n}\|_{L^{2}(\Omega)} &\leq \|\nabla(\cdLh^{n} - \mathcal{S}_{h} \tilde{c}_{h, \delta, L}^{n}) \|_{
            L^{2}(\Omega)} + \|\nabla(\mathcal{S}_{h} \tilde{c}_{h, \delta, L}^{n} - \tilde{c}_{h, \delta, L}^{n})\|_{
            L^{2}(\Omega)} + \| \nabla \tilde{c}_{h, \delta, L}^{n}\|_{L^{2}(\Omega)} \\ &\leq 
           \frac{C}{h} \|\cdLh^{n} - \mathcal{S}_{h} \tilde{c}_{h, \delta, L}^{n} \|_{
            L^{2}(\Omega)} + C \| \nabla \tilde{c}_{h, \delta, L}^{n}\|_{L^{2}(\Omega)} \\
            & \leq \frac{C}{h} \|\cdLh^{n} - \tilde{c}_{h, \delta, L}^{n}\|_{L^{2}(\Omega)} + \frac{C}{h} \|\tilde{c}_{h, \delta, L}^{n} - \mathcal{S}_{h}\tilde{c}_{h, \delta, L}^{n}\|_{L^{2}(\Omega)} + C \| \nabla \tilde{c}_{h, \delta, L}^{n}\|_{L^{2}(\Omega)},
        \end{align*}
        where in the transition from the first to the second line we used (\ref{GradIneq}) to bound the first term on the right-hand side of the first line, and (\ref{InterpIneq}) to bound the second term on the right-hand side of the first line. Because $\cdLh^{n} = -(-\Delta_{h})^{-\s}(\rhodLh^{n})^{\ast}$ and $\tilde{c}_{h, \delta, L}^{n} = -(-\Delta_{\mathrm{N}})^{-\s} (\rhodLh^{n})^{\ast}$, the first term in the last line is now uniformly bounded thanks to Proposition \ref{DiscFracLapConv} with $\sigma = 1/2$ and $\varepsilon(h) = \ln(h^{-1})\,h^{2s+1}$, and the fact that $(\rhodLh^{n})^{\ast}:=\rhodLh^{n} - \int_\Omega \rhodLh^{n}\dx \in V_h$ is uniformly bounded in $H^{2\sigma}_\ast(\Omega)=H^1_\ast(\Omega)$ as $h, \delta \rightarrow 0_+$ because of \eqref{APrioriBound2}.
        The second term in the last line is uniformly bounded by the interpolation inequality (\ref{InterpIneq}), (\ref{H1Stab}) in Lemma \ref{H1stablemma2} and the bound on the third term on the left-hand side of (\ref{APrioriBound2}). For the third term in the last line, we have similarly, thanks to (\ref{H1Stab}), that 
        \[ \| \nabla \tilde{c}_{h, \delta, L}^{n}\|_{L^{2}(\Omega)} \leq C \| \nabla \rhodLh^{n} \|_{L^{2}(\Omega)} \]
        and the bound on the third term on the left-hand side of (\ref{APrioriBound2}) then makes the term uniformly bounded with respect to $h$ and $\delta$.
        Having shown that $\nabla c^n_{h,\delta,L}$ is uniformly bounded in $L^2(\Omega;\mathbb{R}^d)$ as $h, \delta \rightarrow 0_+$, the weak convergence result \eqref{Conv6} follows from \eqref{Conv5} thanks to the uniqueness of the weak limit.
        
        We combine (\ref{Conv2})--(\ref{Conv4}) to pass to the limit as $\delta, h \to 0_{+}$ in (\ref{FullDiscWeakForm}) with $\phi_{h}=\pi_{h} \phi$, for $\phi \in C^{\infty}(\overline{\Omega})$ together with the strong convergence of $\pi_h \phi$ to $\phi$ in the norm of $W^{1,\infty}(\Omega)$ as $h \rightarrow 0_+$
        (cf., for example, inequality (4.4.29) in \cite{BreSco94} with $s=1$, $m=2$, $l=1$, $r=0$ and $p=\infty$ there), to obtain equation (\ref{WeakFormDisctime}). The fact that $C^{\infty}(\overline{\Omega})$ is dense in $V=H^1(\Omega)$ then concludes the argument. The nonnegativity of $\rhodLh^{n}$ on $\overline\Omega$ implies the nonnegativity of $\rhoL^{n}$ a.e.~on $\Omega$ for $n=1, \dots, N$. Finally, we note that if we choose $\phi\equiv 1$ as a test function in (\ref{WeakFormDisctime}) we have conservation of mass, i.e., $\int_\Omega \rho^n_L(x) \dx = \int_\Omega \rho^0 \dx = 1$ for $\rhoL^{n}$, $n=1, \dots, N$. 
    \end{proof}
\end{theo}

\color{black}
Having passed to the limits $h \to 0_{+}$ and $\delta \to 0_{+}$ with the spatial discretization parameter $h$ and the lower cut-off parameter $\delta$, in the next section we shall focus on passage to the limits $\Delta t \to 0_+$ and $L \to +\infty$ with the time step $\Delta t$ and the upper cut-off parameter $L$. Before doing so, however, we close this section by commenting on possible extensions of the analysis pursued in this section to two related models that involve a nonlinear mobility coefficient $\rho \in \mathbb{R}_{\geq 0} \mapsto A(\rho) \in \mathbb{R}_{>0}$.

\begin{remark}[Nonlinear mobility model I] \label{NonLinProbRem}
  Consider the initial-boundary-value problem 
\begin{equation} \label{NonLinearMobProb}
\left \{
\begin{aligned}
&\frac{\partial \rho}{\partial t} = \nabla \cdot( A(\rho) \grad \rho - \rho \nabla c) & \textrm{in } & \Omega \times (0, \infty) , \\
& - (-\Delta)^{\s} c = \rho^{\ast} & \textrm{in } & \Omega \times (0, \infty), \\
&\partial_{n} \rho = 0, \quad \partial_{n} c = 0 & \textrm{on~\!} & \partial \Omega \times (0, \infty), \\
&\rho(x,0) = \rho_0(x)\quad \mbox{for all $x \in \Omega$}, &
\end{aligned}
\right.
\end{equation}
where $A$ is a continuous real-valued function defined on $\mathbb{R}_{\geq 0}$ such that 
\begin{equation} 
\label{Nonlinmobfunc}
   (\exists \sigma_1, \sigma_2 \in \mathbb{R}_{>0}) \, (\forall s \in \mathbb{R}_{\geq 0})\quad \sigma_1 \leq A(s) \leq \sigma_2. 
\end{equation}

For a (nonnegative) continuous piecewise linear function $\phi_h \in V_h$ we then consider the piecewise constant approximation $\Xi(\phi_h)$ of $A(\phi_h)$, defined by
\begin{equation} \label{XiApprox}
    \Xi(\phi_h)(x) = \sum_{K \in \mathcal{T}_h} A(\phi_h(P_K)) \chi_{\mathring{K}}(x), \quad x \in \Omega, 
\end{equation}
where $P_K$ is the barycenter of the symplex $K \in \mathcal{T}_h$ and $\chi_{\mathring{K}}$ is the indicator function of the interior of $K$. 
The fully-discrete approximation of the problem \eqref{NonLinearMobProb} is then defined as follows.  
\begin{gather*}
\textrm{ Let $\rho^0_{h,\delta,L}:=\rho^0_h \in V_h$.}
\textrm{ For } n = 1, \dots, N, \textrm{ given } \rhodLh^{n-1} \in V_{h}, \textrm{ find } \rhodLh^{n} \in V_
{h} \textrm{ such that } \nonumber \\
\label{FullDiscWeakForm} \int_{\Omega} \pi_h\bigg(\frac{\rhodLh^{n} - \rhodLh^{n-1}}{\Delta t} \phi_{h} \bigg)\dx = - \int_{\Omega} \Xi(\rhodLh) \nabla \rhodLh^{n} \cdot \nabla \phi_{h} \dx + \int_{\Omega} \Theta_{\delta}^{L}(\rhodLh^{n}) \nabla \cdLh^{n} \cdot \nabla \phi_{h} \dx \\ \textrm{for all } \phi_{h} \in V_{h},
\end{gather*}
 where $\cdLh^{n} \in V_h \cap L^2_{\ast}(\Omega)$ satisfies \eqref{SpaceDiscFracPois}. 
 Thanks to our assumptions on $A$, it follows that $0<\sigma_1 \leq \Xi(\rhodLh)(x) \leq \sigma_2<+\infty$ for all $x \in \overline\Omega$.
 Under the hypotheses of Lemma \ref{LemLip}, the assumptions on $A$ also imply that
\begin{equation*}
    \int_K \Xi(\phi_h) \nabla\phi_h \cdot \nabla \pi_h(g(\phi_h)) \dx \geq \sigma_1 \int_K \nabla\phi_h \cdot \nabla \pi_h(g(\phi_h)) \dx \qquad \forall\, \phi_h \in V_h.
\end{equation*}
This lower bound is crucial so as to be able to ensure the validity of the bound stated in inequality \eqref{APrioriBound2} in Lemma \ref{lem:uniformdeltah}. In addition, by the Cauchy--Schwarz inequality together with (\ref{H1Stab}), \eqref{Nonlinmobfunc} and the uniform bound on the $H^1(\Omega)$ seminorm, we also have that
\begin{equation*}
    \int_\Omega \Xi(\rhodLh) \nabla c_{h, \delta, L}^n \cdot \nabla \rhodL \dx \leq C \int_\Omega |\nabla \rhodLh|^2 \dx \leq C(L),
\end{equation*}
which then implies that $c_{h, \delta, L}^n$ satisfies the inequality 
\eqref{AprioriBound3}. These considerations enable us to use an identical argument to the one in the proof of Theorem \ref{convergenceFEMTheo} to prove the convergence, as $\delta, h \to 0$,  of a (sub)sequence of the solution sequence of the fully discrete problem to a solution of the following semidiscrete-in-time problem.
\begin{gather*}
\textrm{Let $\rho^0_L:= \rho^0 = \rho^0(\Delta t)$. For } n = 1, \dots, N, \textrm{ given } \rho_{L}^{n-1} \in V, \textrm{ find } \rho_{L}^{n} \in V \textrm{ such that } \nonumber \\
 \int_{\Omega} \frac{\rho_{L}^{n} - \rho_{L}^{n-1}}{\Delta t} \phi \dx = - \int_{\Omega} A(\rho_L^n) \nabla \rho_{L}^{n} \cdot \nabla \phi \dx+ \int_{\Omega} \beta^{L}(\rho_{L}^{n}) \nabla c_{L}^{n} \cdot \nabla \phi \dx \quad \textrm{for all } \phi \in V, \end{gather*}
where $c_L^n$ satisfies \eqref{FinElemFracLap}. 

\end{remark}

\begin{remark}[Nonlinear mobility II]\label{model2}
     Another related problem with nonlinear mobility is given by 
\begin{equation} \label{NonLinProb2}
\left \{
\begin{aligned}
&\frac{\partial \rho}{\partial t} = \nabla \cdot( A(\rho) (\grad \rho - \rho \nabla c)) & \textrm{in } & \Omega \times (0, \infty) , \\
& - (-\Delta)^{\s} c = \rho^{\ast} & \textrm{in } & \Omega \times (0, \infty), \\
&\partial_{n} \rho = 0, \quad \partial_{n} c = 0 & \textrm{on~\!} & \partial \Omega \times (0, \infty), \\
&\rho(x,0) = \rho_0(x)\quad \mbox{for all $x \in \Omega$}, &
\end{aligned}
\right.
\end{equation}
where, again, $A$ is a continuous real-valued function defined on $\mathbb{R}_{\geq 0}$ such that \eqref{Nonlinmobfunc} holds. 
Let $E$ be the free energy functional defined by \eqref{Energy1}. It then follows that the following dissipation law holds: 
\begin{equation*}
    \frac{\dd}{\dd t} E(\rho) = - \int_{\Omega} \rho A(\rho) |\nabla(\log \rho - c)|^2 \dd x \leq 0,
\end{equation*}
which is analogous to that in \eqref{DecayEnergy1}. 
A crucial inequality, which our analysis in the case of constant mobility relied on, was \eqref{NegativeProd}, i.e., that 
\begin{equation*} 
    \int_{\Omega} \nabla \cdLh^{n} \cdot \nabla \rhodLh^{n} \dx \leq 0.
\end{equation*}
While in the case of the nonlinear mobility model I the function $A(\cdot)$ only multiplies $\nabla \rho$, and therefore the above inequality can be guaranteed to hold under the stated assumptions on $A$, just as in the case of constant mobility, in the nonlinear mobility model II the function $A(\cdot)$ now also multiplies $\nabla c$. Identical discretization of $A(\cdot)$ in both terms appearing on the right-hand side of \eqref{NonLinProb2}$_1$ unfortunately results in a numerical method for which the inequality \eqref{NegativeProd} cannot be guaranteed to hold. This means that in the case of the nonlinear mobility model II the coefficient $A(\cdot)$ needs to be approximated differently in the second term on the right-hand side of \eqref{NonLinProb2}$_1$ than in the first term, and this complicates both the definition and the analysis of the resulting scheme. We shall therefore not dwell on this extension here further.  
\end{remark}
\normalcolor
\color{black}

 \section{The semidiscrete-in-time approximation}\label{Sect:4}

Next, we shall pass to the limit as the time-step $\Delta t>0$ tends to zero in the temporally semidiscrete scheme (\ref{WeakFormDisctime}), which we presented in the previous section, and which we obtained as the limit of the finite element approximation when $\delta, h \to 0_+$, as was shown in Theorem \ref{convergenceFEMTheo}.

In order to pass to the limit $\Delta t \rightarrow 0_+$ we shall first show that the approximation is independent of the cut-off parameter $L$, provided that $L>\|\rho_0\|_{L^\infty(\Omega)}$; we shall then derive uniform bounds with respect to $\Delta t$ and proceed by extraction of weakly convergent subsequences as $\Delta t \rightarrow 0_+$, similarly as we did in the previous section in the case of the parameters $\delta$ and $h$.   

Suppose that $N \in \mathbb{N}_{\geq 2}$ and $\Delta t := T/N$. We then define
\begin{subequations}
    \label{LinInterptime}
\begin{equation}
    \rhoL^{\Delta t}(\cdot, t) := \frac{t - t_{n-1}}{\Delta t} \rhoL^{n}(\cdot) + \frac{t_{n} - t}{\Delta t } \rhoL^{n-1}(\cdot), \quad t \in [t_{n-1}, t_{n}], \quad n = 1, \ldots, N,
\end{equation}
that is, the continuous piecewise affine interpolant in time of the sequence of discrete-in-time approxi\-mations $\{ \rhoL^{n} \}_{n =1, \ldots, N}$, in conjunction with the notation 
\begin{equation}
    \rhoL^{\Delta t, +}(\cdot, t) := \rhoL^{n}(\cdot), \quad \rhoL^{\Delta t, -}(\cdot, t) := \rhoL^{n-1}(\cdot), \quad t \in (t_{n-1}, t_{n}], \quad n =1, \ldots, N.
\end{equation}
\end{subequations}

\begin{figure}[H]
\centering
\begin{tikzpicture}[vect/.style={->,
             shorten >=0pt,>=latex'}]
             
\tkzDefPoint(-6, 0){O1}
\tkzDefPoint(-2, 0){O2}
\tkzDefPoint(-6, 2.5){O1v}
\tkzDefPoint(-2, 2.5){O2v}
\tkzDefPoint(-9, 0){O3}
\tkzDefPoint(1, 0){O4}
\tkzDefPoint(-4, 0){Om}
\tkzDefPoint(-8, 0){O3m}
\tkzDefPoint(-8, 2.5){O3mv}
\tkzDefPoint(0, 0){O4m}
\tkzDefPoint(0, 2.5){O4mv}
\tkzDefPoint(-4, 2.5){Omv}
\tkzDefPoint(-9, 0.7){rho11lin}
\tkzDefPoint(-6, 0.3){rho12lin}
\tkzDefPoint(-6, 1.0){rho21lin}
\tkzDefPoint(-2, 1.6){rho22lin}
\tkzDefPoint(-2, 1.0){rho31lin}
\tkzDefPoint(1, 0.6){rho32lin}
\tkzDefPoint(-6, 0.5){rho12linD}
\tkzDefPoint(-8, 0.9){rho12linDm}
\tkzDefPoint(-6, 0.8){rho21linD}
\tkzDefPoint(-4, 1.5){rho21linDm}
\tkzDefPoint(-4, 1.1){rho22linDm}
\tkzDefPoint(-2, 1.8){rho22linD}
\tkzDefPoint(-2, 1.4){rho31linD}
\tkzDefPoint(0, 1.1){rho31linDm}
\tkzDefPoint(2, 0.9){rho32linDm}
\tkzDefPoint(-6, 0.9){F1}
\tkzDefPoint(-2, 1.05){F2}
\tkzDrawSegment(O3, O4)
\tkzDrawSegment(rho12linDm, rho12linD)
\tkzDrawSegment(rho12linD, rho21linDm)
\tkzDrawSegment(rho21linDm, rho22linD)
\tkzDrawSegment(rho22linD, rho31linDm)
\tkzDrawSegment[dashed](O1, O1v)
\tkzDrawSegment[dashed](O2, O2v)
\tkzDrawSegment[dashed](Om, Omv)
\tkzDrawSegment[dashed](O3m, O3mv)
\tkzDrawSegment[dashed](O4m, O4mv)

\tkzDrawPoints[color= black, fill=black](rho12linDm, rho12linD, rho21linDm, rho22linD,  rho31linDm)
\tkzLabelPoint[below](O1){$t_{n-1}$}
\tkzLabelPoint[below](O2){$t_{n+1}$}
\tkzLabelPoint[below](Om){$t_{n}$}
\tkzLabelPoint[below](O3m){$t_{n-2}$}
\tkzLabelPoint[below](O4m){$t_{n+2}$}

\tkzLabelPoint[below right = -0.125](rho12linD){$\rhoL^{n-1}$}
\tkzLabelPoint[above right](rho21linDm){$\rhoL^{n}$}

\tkzLabelSegment[above](rho21linD, rho21linDm){$\rhoLDt$}

\end{tikzpicture}
\caption{Continuous piecewise affine interpolant of the discrete-in-time approximations}
\end{figure}

We shall adopt $\rhoL^{\Delta t (,\pm)}$ as a collective symbol for $\rhoL^{\Delta t}$, $\rhoL^{\Delta t, \pm}$.  For future purposes we observe that 
\begin{equation*} 
\frac{\partial \rhoLDt}{\partial t}(\cdot, t) = \frac{\rhoLDtp(\cdot,t) - \rhoLDtm(\cdot,t)}{\Delta t}, \quad t \in (t_{n-1},t_n],\quad n=1,\ldots, N, 
\end{equation*}
and that, thanks to the final assertion in Theorem \ref{convergenceFEMTheo}, one has that 
\begin{equation*}
    \frac{1}{|\Omega|} \int_{\Omega} \rhoLDtpm(t) \dx =  \frac{1}{|\Omega|} \int_{\Omega} \rho^{0} \dx = 1, \quad t \in (t_{n-1}, t_{n}], \quad n = 1, \ldots, N.
\end{equation*}

For $t \in (0,T]$, let $\cL^{\Delta t}(\cdot,t) \in \mathbb{H}^s(\Omega) \cap H^1_\ast(\Omega)$ be the unique the solution of the equation
\begin{equation} \label{DefcIntTime}
-(-\Delta_{\mathrm{N}})^{\s} \cL^{\Delta t} = (\rhoL^{\Delta t})^{\ast} \quad \mbox{in $\Omega$}.
\end{equation}
We note in passing that, thanks to linearity and because $\frac{1}{|\Omega|} \int_\Omega \rho^n_L \dx = \frac{1}{|\Omega|} \int_\Omega \rho^0_L \dx = \frac{1}{|\Omega|} \int_\Omega \rho_0 \dx = 1$ for all $n=1,\ldots,N$, we have that
\begin{equation*}
  (\rhoL^{\Delta t})^{\ast}(\cdot,t)  := \rhoLDt - \frac{1}{|\Omega|}\int_{\Omega} \rhoL^{\Delta t} \dx = \frac{t - t_{n-1}}{\Delta t} (\rhoL^{n})^{\ast}(\cdot) + \frac{t_{n} - t}{\Delta t } (\rhoL^{n-1})^{\ast}(\cdot), \quad t \in [t_{n-1}, t_{n}], \quad n= 1 \ldots, N,
\end{equation*}
and therefore
\begin{equation} \label{PropPseudoTimeDerAverage}
    \frac{\partial (\rhoLDt)^{\ast}}{\partial t}(\cdot, t) = \frac{(\rhoLDtp)^{\ast}(\cdot,t) - (\rhoLDtm)^{\ast}(\cdot,t)}{\Delta t}= \frac{\partial \rho_L^{\Delta t}}{\partial t}(\cdot,t), \quad t \in (t_{n-1}, t_{n}], \quad n= 1 \ldots, N. 
\end{equation}

Using the above notation, (\ref{WeakFormDisctime}) summed through $n = 1, \ldots, N$, can be restated in a form that is reminiscent of a weak formulation in time, that is:
\begin{subequations} \label{WeakFormSemiDiscTimeLMain}
\begin{gather}
    \text{Find } \rhoL^{\Delta t}(\cdot,t) \in V =H^1(\Omega)\, \mbox{ for } t\in (0, T], \text{ such that } \nonumber \\
    \int_{0}^{T} \int_{\Omega} \frac{\partial \rhoLDt}{\partial t} \phi \dx = - \int_{0}^{T} \int_{\Omega}  \nabla \rhoLDtp \cdot \nabla \phi \dx + \int_{0}^{T} \int_{\Omega} \beta^{L}(\rhoLDtp) \nabla c_{L}^{\Delta t, +} \cdot \nabla \phi \dx \quad \!\!\text{for all } \phi \in L^1(0,T;V), \label{WeakFormSemiDiscTimeL} \end{gather}
    subject to the initial condition $\rho_{L}^{\Delta t}(x, 0) = \rho^{0}(x)$ for a.e.~$x \in \Omega$, where $c_{L}^{\Delta t, +}(\cdot,t) \in H^1_\ast(\Omega)$ for all $t \in (0,T]$ and satisfies
    \begin{equation} \label{FracPoiDiscTime} -(-\Delta_{\mathrm{N}})^{\s} c_{L}^{\Delta t, +}(\cdot,t) = (\rhoLDtp)^{\ast}(\cdot,t)\quad \mbox{in $\Omega$},
\end{equation}
\end{subequations}
for all $t \in (0,T]$, and $\beta^L(s):=\min\{s,L\}$ for $s \in \mathbb{R}$.

We shall now rigorously prove a result similar to the one that we obtained through formal calculations for smooth solutions in Lemma \ref{InfNormBoundClass}, namely a uniform bound in time on the $L^{\infty}(\Omega)$ norm of the solution $\rhoLDt$. This property will allow us to eliminate the parameter $L$. Moreover the bound will be decisive in the convergence argument that follows, and it is there where the structure of the problem and the action of the fractional Laplacian proves advantageous.  As was noted previously, in \cite{chen2022analysis} this result was shown for a smooth solution of the continuous model on the whole of $\mathbb{R}^{d}$ using the integral definition of the fractional Neumann Laplacian, while in Lemma \ref{InfNormBoundClass} above it was asserted for a bounded domain using the spectral definition of the fractional Neumann Laplacian, based on formal calculations. Our starting point is the following lemma.

\begin{lemma}  \label{InftyBoundLemmaDiscTime}
    Let $\rhoLDtp$ be defined as in \eqref{LinInterptime} and let $F \in C^{2}([0, \infty))$ be a convex and nondecreasing function. Then, 
    \begin{equation} \label{BoundFDiscTimeL}
        \sup_{t \in (0, T]} \int_{\Omega} F(\rhoLDtp(t)) \dx \leq \int_{\Omega} F(\rho^{0}) \dx.
    \end{equation}
    \begin{proof}
    We begin by noting that, because $\rho^0 \in L^\infty(\Omega)$, we have that $F(\rho^{0}) \in L^{1}(\Omega)$, and the right-hand side of the inequality \eqref{BoundFDiscTimeL} is therefore finite. 
    Let us assume first that $F''$ is additionally bounded (for the general case one can use a similar argument as the one in the proof of Lemma \ref{InfNormBoundClass}). Let us then pick $\phi=\chi_{[0,t]} F'(\rhoLDtp)$ with $t = t_{n}$ for $n = 1, \ldots, N$, as test function in (\ref{WeakFormSemiDiscTimeL}).  We then  have that
    
    \begin{align} 
    \begin{aligned}\label{Infnorm11}
        \int_{0}^{t_n} \int_{\Omega} \frac{\partial \rhoLDt}{\partial t} F'(\rhoLDtp) \dx\,\dtau = &- \int_{0}^{t_n} \int_{\Omega} F''(\rhoLDtp) |\nabla \rhoLDtp|^{2} \dx\,\dtau \\
        & + \int_{0}^{t_n} \int_{\Omega} \beta^{L}(\rhoLDtp) F''(\rhoLDtp)\nabla c_{L}^{\Delta t, +} \cdot \nabla \rhoLDtp  \dx\,\dtau.   \end{aligned}
        \end{align}

        The first term on the right-hand side is always nonpositive thanks to the convexity of $F$. For the second term on the right-hand side we can write
        \begin{align}
           &\int_{0}^{t_n} \int_{\Omega} \beta^{L}(\rhoLDtp) F''(\rhoLDtp)\nabla c_{L}^{\Delta t, +} \cdot \nabla \rhoLDtp \dx\,\dtau \nonumber \\&\quad= - \int_{0}^{t_n} \int_{\Omega} \beta^{L}(\rhoLDtp) F''(\rhoLDtp) \nabla \rhoLDtp \!\cdot \nabla (- \Delta_{\mathrm{N}})^{-\s} (\rhoLDtp)^{\ast} \dx\,\dtau 
           \label{IntPartsNO}\\
           &\quad= - \int_{0}^{t_n} \int_{\Omega} \nabla H(\rhoLDtp) \cdot \nabla (-\Delta_{\mathrm{N}})^{-\s} (\rhoLDtp)^{\ast} \dx\,\dtau, \nonumber
        \end{align}
        where $H(s) := \int_{0}^{s} \beta^{L}(\tau) F''(\tau) \,\dtau$ is a nondecreasing function by definition. 

We want to employ an argument similar to the one we used in Lemma \ref{InfNormBoundClass} to show that the minus sign on the left-hand side of the fractional Poisson equation (\ref{FracPoiDiscTime}) and the explicit representation of the spectral fractional Neumann Laplacian in terms of the heat kernel guarantee nonpositivity of the expression in the last line of \eqref{IntPartsNO} above, but in order to rigorously replicate the formal argument of Lemma \ref{InfNormBoundClass} we require additional considerations. We would like to perform spatial integration by parts, as we did in (\ref{IntPartsBound}), to prove the nonpositivity of the term (\ref{IntPartsNO}); but now $\rhoLDtp \in L^\infty(0,T;H^{1}(\Omega))$ only, so we cannot proceed directly. We shall therefore consider a sequence of functions $\{ \rho_{L, \epsilon}^{\Delta t, +} \}_{\epsilon > 0} \subset L^\infty(0,T;C^{\infty}(\overline{\Omega}))$, such that $\rho_{L, \epsilon}^{\Delta t, +} \to \rhoLDtp$ strongly in $L^2(0,T;H^{1}(\Omega))$ as $\epsilon \to 0_+$ and $\rho_{L, \epsilon}^{\Delta t, +} \geq 0$ for a.e.~$(x,t) \in \Omega \times (0, T]$. Such a sequence of smooth functions exists thanks to the density of $C^{\infty}(\overline{\Omega})$ in $H^{1}(\Omega)$, and is easily constructed by extending $\rho_{L}^{\Delta t, +}(\cdot,t) \geq 0$ from $\Omega$ to the whole of $\mathbb{R}^d$ by preserving the Sobolev class $H^1(\Omega)$, and convolving the nonnegative part of the extended function with $\theta_\epsilon(x)=\epsilon^{-d} \theta(x/\epsilon)$, where $\theta \in C^\infty_0(\mathbb{R}^d)$ is a nonnegative function such that $\int_{\mathbb{R}^d} \theta(x) \dx = 1$. Instead of the right-hand side of \eqref{IntPartsNO} we shall then consider the term 
\begin{equation*}
    -\int_{0}^{t_n} \int_{\Omega} \nabla H(\rho_{L, \epsilon}^{\Delta t, +}) \cdot \nabla (-\Delta_{\mathrm{N}})^{-\s} (\rho_{L, \epsilon}^{\Delta t, +})^{\ast} \dx\,\dtau. 
\end{equation*}
Proceeding in the same way as we did in Lemma \ref{InfNormBoundClass}, using integration by parts, we have that 
\begin{equation}\label{IntPartsIneqImp1}
   - \int_{0}^{t_n} \int_{\Omega} \nabla H(\rho_{L, \epsilon}^{\Delta t, +}) \cdot \nabla (-\Delta_{\mathrm{N}})^{-\s} (\rho_{L, \epsilon}^{\Delta t, +})^{\ast} \dx\,\dtau = - \int_{0}^{t_n} \int_{\Omega}  H(\rho_{L, \epsilon}^{\Delta t, +})  (-\Delta_{\mathrm{N}})^{1-\s} (\rho_{L, \epsilon}^{\Delta t, +})^{\ast} \dx\,\dtau, 
\end{equation}
and then the nonnegativity of the expression above follows from the same argument as in the case of the second term in \eqref{IntPartsBound} in the proof of Lemma \ref{InfNormBoundClass}.
        
        To prove the nonpositivity of the right-hand side of (\ref{IntPartsNO}) it remains to show that, as $\epsilon \to 0_+$,
        \begin{equation} \label{IntPartsIneqImp2}
           - \int_{0}^{t_n} \int_{\Omega} \nabla H(\rho_{L, \epsilon}^{\Delta t, +}) \cdot \nabla (-\Delta_{\mathrm{N}})^{-\s} (\rho_{L, \epsilon}^{\Delta t, +})^{\ast} \dx\,\dtau \to -\int_{0}^{t_n} \int_{\Omega} \nabla H(\rhoLDtp) \cdot \nabla (-\Delta_{\mathrm{N}})^{-\s} (\rhoLDtp)^{\ast} \dx\,\dtau.
        \end{equation}
To see this, we observe that 
        \begin{align*}
            & \int_{0}^{t_n} \int_{\Omega} (\nabla H(\rho_{L, \epsilon}^{\Delta t, +}) \cdot \nabla (-\Delta_{\mathrm{N}})^{-\s} (\rho_{L, \epsilon}^{\Delta t, +})^{\ast} - \nabla H(\rhoLDtp) \cdot \nabla (-\Delta_{\mathrm{N}})^{-\s} (\rhoLDtp)^{\ast} \dx\,\dtau \\
            &\quad= \underbrace{\int_{0}^{t_n} \int_{\Omega} (\nabla H(\rho_{L, \epsilon}^{\Delta t, +} - \nabla H(\rhoLDtp)) \cdot \nabla (-\Delta_{\mathrm{N}})^{-\s} (\rho_{L, \epsilon}^{\Delta t, +})^{\ast} \dx\,\dtau}_{\circled{1}} \\
            &\qquad  + \underbrace{\int_{0}^{t_n} \int_{\Omega} \nabla H(\rhoLDtp) \cdot \nabla (-\Delta_{\mathrm{N}})^{-\s}((\rho_{L, \epsilon}^{\Delta t, +})^{\ast} - (\rhoLDtp)^{\ast}) \dx \, \dtau }_{\circled{2}}.
        \end{align*}

        Term $\circled{2}$ converges to zero as $\epsilon \to 0_+$, because the function $\nabla H(\rhoLDtp) = \beta^{L}(\rhoLDtp) F''(\rhoLDtp) \nabla \rhoLDtp$ belongs to $L^\infty(0,T; L^{2}(\Omega;\mathbb{R}^d))$ (recall that $\beta^L$ is a bounded function by definition, and $F''$ is a bounded function by assumption), and because for each $\eta \in L^\infty(0,T; C^{\infty}_{0}(\Omega; \mathbb{R}^{d}))$ we have that
        \begin{align*}
            \bigg|\int_0^{t_n} (\eta, \nabla (-\Delta)^{-\s}((\rho_{L, \epsilon}^{\Delta t, +})^{\ast} - (\rhoLDtp)^{\ast})) \dt \bigg| &= \bigg|- \int_{0}^{t_n}(\text{div}\, \eta, (-\Delta_{\mathrm{N}})^{-\s}((\rho_{L, \epsilon}^{\Delta t, +})^{\ast} - (\rhoLDtp)^{\ast})) \dt \bigg|\\
            & \hspace{-17mm} = \bigg| \int_{0}^{t_n} ((-\Delta_{\mathrm{N}})^{-\s}\text{div} \,\eta,(\rho_{L, \epsilon}^{\Delta t, +})^{\ast} - (\rhoLDtp)^{\ast}) \dt \bigg|\\
            &\hspace{-17mm} \leq \| (-\Delta_{\mathrm{N}})^{-\s}\text{div}\, \eta \|_{L^2(0,T;L^{2}(\Omega))} \| (\rho_{L, \epsilon}^{\Delta t, +})^{\ast} - (\rhoLDtp)^{\ast} \|_{L^2(0,T;L^{2}(\Omega))} \to 0
        \end{align*}
        as $\epsilon \to 0_+$.  We note in passing that $\int_\Omega \text{div}\, \eta(x,t)  \dx = \int_{\partial \Omega} \eta(x,t)\cdot n\, \mathrm{d}s = 0$, whereby $\text{div}\, \eta(\cdot,t) \in L^2_\ast(\Omega)$ for a.e.~$t \in (0,T]$ and 
        $(-\Delta_{\mathrm{N}})^{-\s}\text{div}\, \eta(\cdot,t) $ therefore belongs to $L^2_\ast(\Omega)$ for all $s \in (0,1)$, and for a.e.~$t \in (0,T]$. Thanks to the fact that $C^{\infty}_{0}(\Omega; \mathbb{R}^{d})$ is dense in $L^{2}(\Omega;\mathbb{R}^d)$ it follows that 
        \begin{equation*}
            \int_0^{t_n}(\eta, \nabla (-\Delta_{\mathrm{N}})^{-\s}((\rho_{L, \epsilon}^{\Delta t, +})^{\ast} - (\rhoLDtp)^{\ast})) \dt \to 0
        \end{equation*}
        as $\epsilon \to 0_+$ for all $\eta \in L^\infty(0,T; L^{2}(\Omega;\mathbb{R}^d))$, so in particular this is true for $$\eta = \beta^{L}(\rhoLDtp) F''(\rhoLDtp) \nabla \rhoLDtp = \nabla H(\rhoLDtp) \in L^\infty(0,T;L^{2}(\Omega;\mathbb{R}^d)).$$ Thus we have shown that the term \circled{2} converges to 0 as $\epsilon \rightarrow 0_+$.

        For term $\circled{1}$ we have that 
        \begin{align*}
            |\,\circled{1}\,| & = \bigg|\int_{0}^{t_n} \int_{\Omega} ( \beta^{L}(\rho_{L, \epsilon}^{\Delta t, +}) F''(\rho_{L, \epsilon}^{\Delta t, +}) \nabla \rho_{L, \epsilon}^{\Delta t, +} - \beta^{L}(\rhoLDtp) F''( \rhoLDtp) \nabla \rhoLDtp) \cdot \nabla (-\Delta_{\mathrm{N}})^{-\s} (\rhoLDtp)^{\ast} \dx\,\dtau \bigg|\\
            &\quad\leq \| \beta^{L}(\rho_{L, \epsilon}^{\Delta t, +}) F''(\rho_{L, \epsilon}^{\Delta t, +}) \nabla \rho_{L, \epsilon}^{\Delta t, +} - \beta^{L}(\rhoLDtp) F''( \rhoLDtp) \nabla \rhoLDtp \|_{L^{2}(0, T; L^{2}(\Omega))} \\&\hspace{6.3cm} \times \|  \nabla (-\Delta_{\mathrm{N}})^{-\s} (\rho_{L, \epsilon}^{\Delta t, +})^{\ast} \|_{L^{2}(0, T; L^{2}(\Omega))}.
        \end{align*}
Focusing on the first factor on the right-hand side of the last inequality, we have that  
        \begin{align}
        \begin{aligned}\label{IntPartsIneq1}
           & \| \beta^{L}(\rho_{L, \epsilon}^{\Delta t, +}) F''(\rho_{L, \epsilon}^{\Delta t, +}) \nabla \rho_{L, \epsilon}^{\Delta t, +} - \beta^{L}(\rhoLDtp) F''( \rhoLDtp) \nabla \rhoLDtp \|_{L^{2}(0, T; L^{2}(\Omega))} \\
            &\quad \leq   \| \beta^{L}(\rho_{L, \epsilon}^{\Delta t, +}) F''(\rho_{L, \epsilon}^{\Delta t, +})(\nabla \rho_{L, \epsilon}^{\Delta t, +} - \nabla \rhoLDtp) \|_{L^{2}(0, T; L^{2}(\Omega))}  \\
            &\qquad + \| (\beta^{L}(\rho_{L, \epsilon}^{\Delta t, +}) F''(\rho_{L, \epsilon}^{\Delta t, +}) - \beta^{L}(\rhoLDtp) F''( \rhoLDtp)) \nabla \rhoLDtp \|_{L^{2}(0, T; L^{2}(\Omega))}.
        \end{aligned}
        \end{align}
        Let us define the constant $M := \sup_{s \in [0,\infty)} F''(s)\geq 0$ (recall that $F$ has been assumed to be a convex function, with $F''$ bounded on $[0,\infty)$), the first term on the right-hand side of \eqref{IntPartsIneq1} is bounded by $L M \| \nabla \rho_{L, \epsilon}^{\Delta t, +} - \nabla \rhoLDtp \|_{L^{2}(0, T; L^{2}(\Omega))}$, which converges to 0 as $\epsilon \to 0_+$ thanks to the strong convergence of $\rho_{L, \epsilon}^{\Delta t, +}$ to $\rhoLDtp$ in $L^2(0,T;H^{1}(\Omega))$. For the second term on the right-hand side of \eqref{IntPartsIneq1}, because of the strong convergence of $\rho_{L, \epsilon}^{\Delta t, +}$ to $\rhoLDtp$ in $L^2(0,T;H^{1}(\Omega))$ there exists a subsequence (not indicated) such that $\rho_{L, \epsilon}^{\Delta t, +}(x,t) \to \rhoLDtp(x,t)$ as $\epsilon \to 0_+$ for a.e.~$(x,t) \in \Omega\times (0,T)$. As $\beta^{L}$ and $F''$ are continuous functions we have also that $\beta^{L}(\rho_{L, \epsilon}^{\Delta t, +}(x,t)) F''(\rho_{L, \epsilon}^{\Delta t, +}(x,t))$ converges to $\beta^{L}(\rho_{L}^{\Delta t, +}(x,t)) F''(\rho_{L}^{\Delta t, +}(x,t))$ for a.e.~$(x,t) \in \Omega\times (0,T)$ as $\epsilon \to 0_+$. Hence, by defining
\begin{equation*}
  A_\epsilon(x,t) =  | (\beta^{L}(\rho_{L, \epsilon}^{\Delta t, +}(x,t)) F''(\rho_{L, \epsilon}^{\Delta t, +}(x,t)) - \beta^{L}(\rhoLDtp(x,t)) F''( \rhoLDtp(x,t))) \nabla \rhoLDtp(x,t) |^2,
\end{equation*}
we have also that $A_\epsilon(x,t) \to 0_+$ as $\epsilon \to 0_+$ for a.e.~$(x,t) \in \Omega\times(0,T)$. Let us 
note furthermore that 
$0 \leq  A_\epsilon(x,t) \leq 4L^2M^2 | \nabla \rhoLDtp(x,t)|^2$ and that the function $4L^2M^2 |\nabla \rhoLDtp|^2 \in L^1(\Omega\times (0,T))$. Thus, by Lebesgue's dominated convergence theorem $\int_0^T \int_\Omega A_\epsilon(x,t) \dx \dt \rightarrow 0_+$ as $\epsilon \rightarrow 0_+$, and consequently the second term on the right-hand side of \eqref{IntPartsIneq1} also converges to 0
as $\epsilon \rightarrow 0_+$. Thus we have shown that the term $\circled{1}$ converges to 0 as $\epsilon \rightarrow 0_+$, which then completes the proof of \eqref{IntPartsIneqImp2}.
        
By combining \eqref{IntPartsIneqImp1} and \eqref{IntPartsIneqImp2}, we have proved that the entire right-hand side of (\ref{Infnorm11}) is nonpositive. For the term on the left-hand side, by Taylor expansion with a remainder of the function
\[ s \in [0, \infty) \mapsto F(s + \alpha) \in [0, \infty), \]
where $\alpha \in [0,\infty)$, we have that, for any $b \in [0, \infty)$,
\[ (s-b) F'(s + \alpha) = F(s + \alpha) - F(b + \alpha) + \frac{1}{2}(s - b)^{2} F''(\theta s + (1-\theta)b + \alpha),\] with $\theta \in (0,1)$. By noting that $t\in [0, T] \mapsto \rhoLDt(\cdot,t)$ is piecewise affine relative to the temporal partition $0 = t_{0} < t_{1} < \cdots < t_{N} = T$, we have that
        \begin{align*}
            \int_{0}^{t_n} \int_{\Omega} \frac{\partial \rhoLDt}{\partial t} F'(\rhoLDtp) \dx\,\dtau &= \int_{\Omega} F(\rhoLDtp(x,t_n)) \dx - \int_{\Omega} F(\rho^{0}(x)) \dx \\
            & \quad + \frac{1}{2 \Delta t} \int_{0}^{t_n} \int_{\Omega} F''(\theta \rhoLDtp + (1-\theta) \rhoLDtm) (\rhoLDtp - \rhoLDtm)^{2} \dx\,\dtau.
        \end{align*}
Because the second derivative of $F$ is nonnegative thanks to the assumed convexity of $F$, we deduce that
        \begin{equation*}
            \int_{\Omega} F(\rhoLDtp(x,t_n)) \dx \leq \int_{\Omega} F(\rho^{0}(x)) \dx \quad \mbox{for all $n \in \{1,\ldots,N\}$,}
        \end{equation*}
        and therefore,  because $\rhoLDtp(\cdot, t)$ is constant on each interval $(t_{n-1}, t_{n}]$ for $n \in \{1,\ldots,N\}$, the inequality \eqref{BoundFDiscTimeL} directly follows. 

The assumption made at the start of the proof  that $F''$ is bounded can be removed as at the end of the proof of Lemma        \ref{InfNormBoundClass}: i.e., we define $F_k(u) := F(0) + F'(0)u + \int_0^u \int_0^v \min(F''(w),k) \, \mathrm{d}w\, \mathrm{d}v$ for $k>0$;
then, $F''_k (u)$ is bounded and the assertion of the lemma follows from the inequality \eqref{BoundFDiscTimeL} with $F$ replaced by $F_k$; we then 
take the limit $k \rightarrow \infty$ using the monotone convergence theorem. \end{proof}
\end{lemma}

By recalling the definitions (\ref{LinInterptime}) we can write, for $t \in (t_{n-1},t_n]$, $n=2,\ldots,N$,
\begin{equation*}
\int_{\Omega} F(\rhoLDtm(x, t)) \dx = \int_{\Omega} F(\rhoLDtp(x, t - \Delta t)) \dx, \end{equation*}
and for $0=:t_{0} < t \leq t_{1}$ we have $\rhoLDtm(\cdot,t) = \rho^{0}(\cdot)$, and therefore we have the same bound (\ref{BoundFDiscTimeL}) also for $F(\rhoLDtm)$ with $\sup_{t \in (0,T]}$ replaced by $\sup_{t \in (0,T)}$, and consequently, by the convexity of $F$, for $F(\rhoLDt)$, 
with $\sup_{t \in (0,T]}$ replaced by $\sup_{t \in [0,T]}$.
These results can now be used to obtain a result that is analogous to Lemma \ref{InfNormDecayClass}. 

\begin{lemma} \label{InfNormDecayDiscTime}
Let $\rhoLDtpm$ be defined as in \eqref{LinInterptime}. Then, the following bound holds:
\begin{equation} \label{InfNormBoundDiscTime}
  \sup_{t \in (0, T)} \| \rhoLDtpm(t) \|_{L^{\infty}(\Omega)} \leq \| \rho_{0}\|_{L^{\infty}(\Omega)}. 
\end{equation}
\begin{proof}
First notice that if we test (\ref{InValSemiDiscProb}) with $F'(\rho^{0})$, where $F \in C^2([0,\infty))$ is a nonnegative convex function, then we have, thanks to the nonnegativity of $F''$, that 
\[ \int_{\Omega} (\rho^{0} - \rho_{0}) F'(\rho^{0}) \dx \leq 0, \qquad \mbox{and therefore also} \qquad \int_{\Omega} F(\rho^{0}) \dx \leq \int_{\Omega} F(\rho_{0}) \dx, \]
thanks to the convexity of $F$. Let us now choose a nonnegative convex function $F \in C^{2}([0, \infty))$ such that $F(u)=0$ for $u \leq \| \rho_{0}\|_{L^{\infty}(\Omega)}$ and $F(u) > 0$ for $u > \| \rho_{0}\|_{L^{\infty}(\Omega)}$. Then, 
\[ 0 \leq \int_{\Omega} F(\rhoLDtpm(x,t)) \dx \leq \int_{\Omega} F(\rho_{0}(x)) \dx = 0 \quad \text{for all } t \in (0, T).\]
Consequently $\rhoLDtpm \leq \| \rho_{0}\|_{L^{\infty}(\Omega)}$, which then implies the stated $L^{\infty}(\Omega)$ norm bound.
\end{proof}
\end{lemma}

Thanks to the fact that (\ref{InfNormBoundDiscTime}) provides a uniform bound on the $L^{\infty}(\Omega)$ norm of $\rhoDtpm$ for all $t \in (0,T)$, dependent only on the $L^{\infty}(\Omega)$ norm of the initial datum, the temporally semidiscrete weak formulation \eqref{PropPseudoTimeDerAverage}, \eqref{WeakFormSemiDiscTimeLMain}, which we obtained by passing to the limits $h, \delta \rightarrow 0_+$, can be simplified: by assuming without loss of generality that $L>\|\rho_0\|_{L^\infty(\Omega)}$, the (upper) cut-off parameter $L$  can be eliminated from \eqref{WeakFormSemiDiscTimeLMain} by recalling that $\beta^L(s):=\min\{s,L\}$ and noting that $\beta^L(s) = s$ for $L>s$.

Hence, by assuming without loss of generality that $L>\|\rho_0\|_{L^\infty(\Omega)}$, the semidiscrete-in-time formulation of our problem becomes the following:
\begin{subequations} \label{WeakFormSemiDiscTimeMain}
\begin{gather}
    \text{Find } \rho^{\Delta t}(\cdot,t) \in V=H^1(\Omega), \; t \in (0, T], \text{ such that } \nonumber \\
    \int_{0}^{T} \int_{\Omega} \frac{\partial \rhoDt}{\partial t} \phi \dx = - \int_{0}^{T} \int_{\Omega}  \nabla \rhoDtp \cdot \nabla \phi  \dx+ \int_{0}^{T} \int_{\Omega} \rhoDtp \nabla c^{\Delta t, +} \cdot \nabla \phi \dx\quad \text{for all } \phi \in L^1(0,T;V), \label{WeakFormSemiDiscTime} \end{gather}
    subject to the initial condition $\rho^{\Delta t}(x, 0) = \rho^{0}(x)$ for a.e.~$x \in \Omega$, where $c^{\Delta t, +}(\cdot,t) \in H^1_\ast(\Omega)$ for all $t \in (0,T]$ and satisfies
        \begin{equation} -(-\Delta_{\mathrm{N}})^{\s} c^{\Delta t, +}(\cdot,t) = (\rhoDtp)^{\ast}(\cdot,t)\quad \mbox{in $\Omega$}
\end{equation}
\end{subequations}
for all $t \in (0,T]$, and the definitions of $\rhoDtpm$ are analogous to those of $\rhoDtpm_L$ stated in (\ref{LinInterptime}). 

We can now start to derive a bound that is uniform in $\Delta t$. Our candidate test function will be $\phi = \chi_{[0, t]} G'(\rhoDtp)$, where $\chi_{[0,t]}$ denotes the characteristic function of the interval $[0,t]
$. However, since $\lim _{s \rightarrow 0_+}G'(s)=+\infty$, we need to be careful because, although the convergence results (\ref{Convergence}) ensure that $\rhoDtp$ is nonnegative a.e.~on $\Omega \times (0,T]$, there is no reason why $\rhoDtp$ should be strictly positive a.e.~on $\Omega \times (0,T]$. We shall circumvent this problem by choosing as a test function $\chi_{[0, t]} G'(\rhoDtp + \alpha)$ where $\alpha > 0$. In this way  $G'(\rhoDtp + \alpha)$ and $G''(\rhoDtp + \alpha)$ will be well defined. After a series of manipulations we shall arrive at a bound involving $G(\rhoDtp+\alpha)$ only. At that point we will be able to pass to the limit as $\alpha \to 0_{+}$, noting that, unlike $G'(\rhoDtp)$ and $G''(\rhoDtp)$, which are only defined when $\rhoDtp$ is positive, the function $G(\rhoDtp)$ is well defined for nonnegative $\rhoDtp$. 

Thus we take $\alpha \in (0,1)$ and we choose $\phi =\chi_{[0, t]} G'(\rhoDtp+ \alpha)$, with $t = t_{n}$, $n \in \{ 1, \ldots, N \}$. We have 
\begin{align*}
    \underbrace{\int_{0}^{T} \int_{\Omega} \frac{\partial \rhoDt}{\partial s} \chi_{[0,t_n]} G'(\rhoDtp + \alpha) \dx \dt}_{\circled{1}} &+ \underbrace{\int_{0}^{T} \int_{\Omega} \nabla \rhoDtp \cdot \nabla \chi_{[0,t_n]} G'(\rhoDtp + \alpha) \dx \dt}_{\circled{2}} \\
    &= \underbrace{\int_{0}^{T} \int_{\Omega} \rhoDtp \nabla c^{\Delta t, +} \cdot \nabla\chi_{[0,t_n]} G'(\rhoDtp + \alpha) \dx \dt}_{\circled{3}}. 
\end{align*}
Let us now manipulate the three terms one at a time. First, by Taylor expansion with a remainder of the function
\[ s \in [0, \infty) \mapsto G(s + \alpha) \in [0, \infty), \]
we have that, for any $b \in [0, \infty)$,
\[ (s-b) G'(s + \alpha) = G(s + \alpha) - G(b + \alpha) + \frac{1}{2}(s - b)^{2} G''(\theta s + (1-\theta)b + \alpha),\] with $\theta \in (0,1)$. Noting that for $t \in [0,T]$ the function $t \mapsto \rhoDt(\cdot,t)$ is, by definition, piecewise affine relative to the partition $0 = t_{0}<t_{1}< \cdots <t_{N} = T$ of the interval $[0,T]$, it follows that 
\begin{align*}
\circled{1} &= \int_{\Omega} G(\rhoDtp(x,t_n) + \alpha) \dx - \int_{\Omega} G(\rho^{0}(x) + \alpha) \dx \\
& \quad + \frac{1}{2 \Delta t} \int_{0}^{t_n} \int_{\Omega} G''(\theta \rhoDtp + (1-\theta) \rhoDtm + \alpha) (\rhoDtp - \rhoDtm)^{2} \dx\,\dtau.
\end{align*}
As $G''(s + \alpha) = \frac{1}{s+\alpha}$ for all $s \geq 0$ and (\ref{InfNormBoundDiscTime}) implies that $G''(\theta \rhoDtp + (1-\theta) \rhoDtm + \alpha) \geq \frac{1}{\| \rho_{0}\|_{L^{\infty}(\Omega)} + \alpha}$ almost everywhere in $\Omega\times (0,T)$ it follows that
\begin{align*}
    \circled{1} \geq \int_{\Omega} G(\rhoDtp(x,t_n) + \alpha) \dx &- \int_{\Omega} G(\rho^{0}(x) + \alpha) \dx \\
    &+ \frac{1}{2 \Delta t \, (\| \rho_{0}\|_{L^{\infty}(\Omega)} + \alpha)} \int_{0}^{t_n} \int_{\Omega} (\rhoDtp - \rhoDtm)^{2} \dx\,\dtau.
\end{align*}

For the second term we have 
\begin{align*}
    \circled{2} = \int_{0}^{t_n} \int_{\Omega} G''(\rhoDtp + \alpha) |\nabla \rhoDtp|^{2} \dx \, \,\dtau = \int_{0}^{t_n} \int_{\Omega} \frac{|\nabla \rhoDtp|^{2}}{\rhoDtp + \alpha} \dx \, \,\dtau,
\end{align*}
again, since $(G^{L})''(s+\alpha) = \frac{1}{s+\alpha}$ for all $s \geq 0$.

For the third term we have
\begin{align*}
    \circled{3} = \int_{0}^{t_n} \int_{\Omega} \frac{\rhoDtp}{\rhoDtp + \alpha} \nabla c^{\Delta t, +} \cdot \nabla \rhoDtp \dx\,\dtau.
\end{align*}

Putting everything together we therefore have that  
\begin{align}
    \int_{\Omega} G(\rhoDtp(x,t_n) + \alpha) \dx &+ \frac{1}{2 \Delta t \, (\| \rho_{0}\|_{L^{\infty}(\Omega)} + \alpha)} \int_{0}^{t_n} \int_{\Omega} (\rhoDtp - \rhoDtm)^{2} \dx\,\dtau \nonumber \\
    &+\int_{0}^{t_n} \int_{\Omega} \frac{|\nabla \rhoDtp|^{2}}{\rhoDtp + \alpha} \dx \, \,\dtau - \int_{0}^{t_n} \int_{\Omega} \frac{\rhoDtp}{\rhoDtp + \alpha} \nabla c^{\Delta t, +} \cdot \nabla \rhoDtp \dx\,\dtau \label{StabEstAlpha} \\ &\hspace{-7mm} \leq \int_{\Omega} G(\rho^{0}(x) + \alpha) \dx \nonumber.
\end{align}
We can now pass to the limit as $\alpha \to 0_+$ to obtain the following bound for $n=1,\ldots,N$, which is uniform in $\Delta t$ thanks to (\ref{IneqPropIn}):
\begin{align} \int_{\Omega} G(\rhoDtp(x,t_n)) \dx &+ \frac{1}{2 \Delta t \, \| \rho_{0}\|_{L^{\infty}(\Omega)}} \int_{0}^{t_n} \int_{\Omega} (\rhoDtp - \rhoDtm)^{2} \dx\,\dtau  \nonumber \\  & + 4\int_{0}^{t_n} \int_{\Omega} \bigg|\nabla \sqrt{\rhoDtp}\bigg|^{2} \dx\,\dtau - \int_{0}^{t_n} \int_{\Omega} \nabla c^{\Delta t, +} \cdot \nabla \rhoDtp \dx\,\dtau \label{RhoIneqEnergy1} \\ &\hspace{-7mm}\leq \int_{\Omega}G(\rho^{0}(x))  \dx \leq \int_{\Omega} G(\rho_{0})  \dx =: C(\rho_{0}). \nonumber \end{align}  
Finally we notice that, for $t \in (0,T],$
\begin{equation} \label{NegativeProdTime} \int_{\Omega} \nabla c^{\Delta t, +}(x,t) \cdot \nabla \rhoDtp(x,t) \dx \leq 0; \end{equation} 
this can be proved by using the density of $C^{\infty}(\overline{\Omega})$ functions in $H^{1}(\Omega)$ and a similar argument as the one contained in the proof of Lemma \ref{InftyBoundLemmaDiscTime} to show the nonpositivity of the term (\ref{IntPartsNO}). We therefore have, for $n=1,\ldots,N$, that
\begin{align}\label{UniformBoundSpatialDer}
\begin{aligned}
 \int_{\Omega} G(\rhoDtp(x,t_n)) \dx & + \frac{1}{2 \Delta t \, \| \rho_{0}\|_{L^{\infty}(\Omega)} } \int_{0}^{t_n} \int_{\Omega} (\rhoDtp - \rhoDtm)^{2} \dx\,\dtau \\
 & + 4\int_{0}^{t_n} \int_{\Omega} \bigg|\nabla \sqrt{\rhoDtp}\bigg|^{2} \dx\,\dtau  \leq C(\rho_{0}). 
\end{aligned}
\end{align}

Using property (\ref{IneqPropIn}) we can supplement the above inequality with additional bounds. The first of these is arrived at by noticing that 
\begin{align*}
 4 \int_{0}^{T} \int_{\Omega}  \bigg| \nabla \sqrt{\rhoDtm} \bigg|^{2} \dx  \dt 
 &= 4 \Delta t \int_{\Omega} \bigg| \nabla \sqrt{\rho^{0}} \bigg|^{2} \dx + 4 \int_{\Delta t}^{T } \int_{\Omega} \bigg| \nabla \sqrt{\rhoDtm} \bigg|^{2} \dx \dt \\
 &= 4 \Delta t \int_{\Omega} \Big| \nabla \sqrt{\rho^{0}} \bigg|^{2} \dx + 4 \int_{0}^{T - \Delta t} \int_{\Omega} \Big| \nabla \sqrt{\rhoDtp} \bigg|^{2} \dx  \dt \\
 & \leq \int_{\Omega} G(\rho_{0}) \dx  + 4 \int_{0}^{T - \Delta t} \int_{\Omega} \Big| \nabla \sqrt{\rhoDtp} \bigg|^{2} \dx \dt \leq C_{\star},
\end{align*}
where in the two inequalities appearing in the last line we used  \eqref{IneqPropIn} and (\ref{UniformBoundSpatialDer}), respectively; here and henceforth $C_{\star}$ signifies a generic positive constant, independent of $\Delta t$.

A simple calculation then shows that we can also get a similar inequality for $\rhoDt$, and therefore we have the uniform bound 
\begin{equation} \label{UniformBoundSpatialDerPlus}
    4 \int_{0}^{T} \int_{\Omega} \bigg| \nabla \sqrt{\rhoDtpm} \bigg|^{2} \dx  \dt \leq C_{\star}.
\end{equation}

At this stage we can also obtain a uniform bound on the norm of $c^{\Delta t, +}$. By choosing $-\chi_{[0,t_n]} c^{\Delta t, +}$ as a test function in (\ref{WeakFormSemiDiscTime}) and recalling \eqref{PropPseudoTimeDerAverage}  we have, for each $n \in \{1,\ldots,N\}$, that 

\begin{align*}
    - \int_{0}^{t_n} \int_{\Omega} \frac{\partial (\rhoDt)^{\ast}}{\partial t} c^{\Delta t, +} \dx\,\dtau + \int_{0}^{t_n} \int_{\Omega} \rhoDtp \nabla c^{\Delta t, +} \cdot \nabla  c^{\Delta t, +} \dx\,\dtau = \int_{0}^{t_n} \int_\Omega \nabla \rhoDtp \cdot \nabla c^{\Delta t, +} \dx \,\dtau,
\end{align*}
which is equivalent to 
\begin{align*}
    \underbrace{- \int_{0}^{t_n} \int_{\Omega} \frac{\partial (\rhoDt)^{\ast}}{\partial t} c^{\Delta t} \dx\,\dtau}_{\circled{1}} &+ \underbrace{\int_{0}^{t_n} \int_{\Omega} \frac{\partial (\rhoDt)^{\ast}}{\partial t} \big(c^{\Delta t} - c^{\Delta t,+}\big) \dx\,\dtau}_{\circled{2}} \\
    &+ \int_{0}^{t_n} \int_{\Omega} \rhoDtp \Big| \nabla c^{\Delta t, +} \Big|^{2} \dx\,\dtau = \int_{0}^{t_n} \int_\Omega \nabla \rhoDtp \cdot \nabla c^{\Delta t, +} \dx \,\dtau.
\end{align*}

The idea at this point is to extract from term $\circled{1}$ the $\mathbb{H}^{\s}$ norm of $c^{\Delta t, +}$ in a similar way as we did in (\ref{StabEst}). Expanding and using (\ref{DefcIntTime}) and (\ref{PropPseudoTimeDerAverage}), thanks to the fact that the fractional Laplacian is an operator acting only in space, we have for each $n \in \{1,\ldots, N\}$ that
\begin{align*}
    \circled{1} & = -\sum_{m=1}^{n} \int_{t_{m-1}}^{t_{m}} \frac{\partial(\rhoDt)^{\ast}}{\partial t} c^{\Delta t} \dx \,\dtau = - \sum_{m=1}^{n} \int_{t_{m-1}}^{t_{m}} \int_{\Omega} \frac{(\rhoDtp)^{\ast} - (\rhoDtm)^{\ast}}{\Delta t} c^{\Delta t} \dx \,\dtau \\
    &= \sum_{m=1}^{n} \int_{\Omega} \frac{(\rho^{m})^{\ast} - (\rho^{m-1})^{\ast}}{\Delta t} \bigg( \int_{t_{m-1}}^{t_{m}} \!\!\!\frac{t-t_{m-1}}{\Delta t} (-\Delta_{\mathrm{N}})^{-s} (\rho^{m})^{\ast} \,\dtau + \int_{t_{m-1}}^{t_{m}} \!\!\!\frac{t_{m}-t}{\Delta t} (-\Delta_{\mathrm{N}})^{-s}(\rho^{m-1})^{\ast} \dtau \bigg) \dx \\
    & = \frac{1}{2} \sum_{m=1}^{n} \int_{\Omega} ((\rho^{m})^{\ast} - (\rho^{m-1})^{\ast})\,((-\Delta_{\mathrm{N}})^{-s} (\rho^{m})^{\ast} + (-\Delta_{\mathrm{N}})^{-s}(\rho^{m-1})^{\ast}) \dx \\
    &= \frac{1}{2} \int_{\Omega} (\rho^{n})^{\ast}\, (-\Delta_{\mathrm{N}})^{-s}(\rho^{n})^{\ast} \dx - \frac{1}{2} \int_{\Omega} (\rho^{0})^{\ast}\, (-\Delta_{\mathrm{N}})^{-s} (\rho^{0})^{\ast} \dx \\
    &= -\frac{1}{2} \int_{\Omega} (\rhoDtp)^{\ast}(x,t_n)\, c^{\Delta t, +}(x,t_n) \dx + \frac{1}{2} \int_{\Omega} (\rho^{0})^{\ast}(x)\, c^{0}(x) \dx \\
    &= \frac{1}{2} \int_{\Omega} (-\Delta_{\mathrm{N}})^{\s} c^{\Delta t,+}(x,t_n) \, c^{\Delta t,+}(x,t_n) \dx - \frac{1}{2} \int_{\Omega} (-\Delta_{\mathrm{N}})^{\s} c^{0}(x) \, c^{0}(x) \dx.
\end{align*}

We now wish to show that the term $\circled{2}$ is nonnegative. We first notice that on the interval $[t_{n-1}, t_{n}]$, by proceeding similarly as in the case of term \circled{1}, we have by linearity that 
\begin{equation*} 
    c^{\Delta t} - c^{\Delta t, +} = \frac{t - t_{n}}{\Delta t} (c^{\Delta t, +} - c^{\Delta t, -})= \frac{t_{n} - t}{\Delta t} (-\Delta_{\mathrm{N}})^{-\s}((\rhoDtp)^{\ast} - (\rhoDtm)^{\ast}).
\end{equation*}
Therefore, we have that 
\begin{align*}
  \circled{2} &= \int_{0}^{t_n} \int_{\Omega} \frac{(\rhoDtp)^{\ast} - (\rhoDtm)^{\ast}}{\Delta t} (c^{\Delta t} - c^{\Delta t, +} ) \dx\,\dtau \\
  &= \sum_{m=1}^{n} \int_{t_{m-1}}^{t_{m}} \int_{\Omega} \frac{(\rhoDtp)^{\ast} - (\rhoDtm)^{\ast}}{\Delta t} \frac{t_{m} - t}{\Delta t} (-\Delta_{\mathrm{N}})^{-\s} ((\rhoDtp)^{\ast} - (\rhoDtm)^{\ast}) \dx\,\dtau \\
  &= \sum_{m = 1}^{n} \int_{t_{m-1}}^{t_{m}} \frac{t_{m} - t} {\Delta t} \int_{\Omega}((\rhoDtp)^{\ast} - (\rhoDtm)^{\ast}) (-\Delta_{\mathrm{N}})^{-\s} ((\rhoDtp)^{\ast} - (\rhoDtm)^{\ast}) \dx\,\dtau. 
\end{align*}
Since the inverse fractional Laplacian is a positive operator, we have that $\circled{2} \geq 0$, and therefore
\begin{align}\begin{aligned}
    \frac{1}{2} \int_{\Omega} (-\Delta_{\mathrm{N}})^{\s} c^{\Delta t,+}(x,t_n) \, c^{\Delta t,+}(x,t_n) \dx &+ \int_{0}^{t_n} \int_{\Omega} \rhoDtp \Big| \nabla c^{\Delta t, +} \Big|^{2} \dx\dt \\
    &\hspace{-17mm} -\int_{0}^{t_{n}} \int_{\Omega} \nabla \rhoDtp \cdot \nabla c^{\Delta t, +} \dx \dt \leq \frac{1}{2} \int_{\Omega} (-\Delta_{\mathrm{N}})^{\s} c^{0}(x) \, c^{0}(x) \dx, \label{EnrgyIneqPart2}
\end{aligned}
\end{align}
which means, by \eqref{NegativeProdTime}, that
\begin{equation*} 
\frac{1}{2} \int_{\Omega} (-\Delta_{\mathrm{N}})^{\s} c^{\Delta t,+}(x,t_n) \, c^{\Delta t,+}(x,t_n) \dx + \int_{0}^{t_n} \int_{\Omega} \rhoDtp \Big| \nabla c^{\Delta t, +} \Big|^{2} \dx\dt \leq \frac{1}{2} \int_{\Omega} (-\Delta_{\mathrm{N}})^{\s} c^{0}(x) \, c^{0}(x) \dx
\end{equation*}
for all $n=1,\ldots,N$. Notice that, by (\ref{StabFracNeu}),
\begin{equation*}
    \int_{\Omega} (-\Delta_{\mathrm{N}})^{\s} c^{0} \, c^{0} \dx = \|c^{0}\|_{\mathbb{H}^{\s}_{\ast}(\Omega)}^{2} \leq C \|\rho^{0}\|_{L^{2}(\Omega)}^{2} \leq C \| \rho_{0}\|_{L^{2}(\Omega)}^{2} = C_{\ast}.
\end{equation*}

Therefore we have the following overall (uniform in $\Delta t$) bound  on $c^{\Delta t,+}$:
\begin{equation} \label{UniformBoundC}
    \sup_{t \in [0,T]} \int_{\Omega} (-\Delta_{\mathrm{N}})^{\s} c^{\Delta t,+}(t) c^{\Delta t,+}(t) \dx + \int_{0}^{t} \int_{\Omega} \rhoDtp \Big| \nabla c^{\Delta t, +} \Big|^{2} \dx \dt \leq C_{\ast}.
\end{equation}

We observe that \eqref{RhoIneqEnergy1}
in conjunction with \eqref{EnrgyIneqPart2} leads, for any $n\in\{1,\ldots,N\}$,  to the following bound on the energy functional that we defined in \eqref{Energy1}: 
\begin{align}
E(\rhoDtp(t_{n})) - E(\rho^{0})&= \int_{\Omega} G(\rhoDtp(x, t_{n})) \dx - \frac{1}{2} \int_{\Omega} (\rhoDtp)^{\ast}(x, t_{n}) c^{\Delta t, +}(x, t_{n}) \dx  - \int_{\Omega} G(\rho^{0}(x)) \dx \nonumber \\
& \quad + \frac{1}{2} \int_{\Omega} (\rho^{0})^{\ast}(x) c^{0}(x) \dx \leq  -4\int_{0}^{t_n} \int_{\Omega} \bigg|\nabla \sqrt{\rhoDtp}\bigg|^{2} \dx \dt \nonumber \\
& \quad + 2\int_{0}^{t_{n}} \int_{\Omega} \nabla \rhoDtp \cdot \nabla c^{\Delta t, +} \dx \dt  \quad - \int_{0}^{t_{n}} \int_{\Omega} \rhoDtp \Big| \nabla c^{\Delta t, +} \Big|^{2} \dx \dt \label{EnergyIneqSemidisc}.
\end{align}
The inequality  \eqref{EnergyIneqSemidisc} can be rewritten, for any $n\in\{1,\ldots,N\}$, as 
\begin{align} \label{EnergyDissipPerfSquare}
E(\rhoDtp(t_{n})) + \int_{0}^{t_{n}} \int_{\Omega} \bigg| 2 \nabla \sqrt{\rhoDtp} - \sqrt{\rhoDtp} \nabla c^{\Delta t, +} \bigg|^{2} \dx \dt  \leq E(\rho^{0}).
\end{align}
Because $| 2 \nabla \sqrt{\rhoDtp} - \sqrt{\rhoDtp} \nabla c^{\Delta t, +} |^{2} = \rhoDtp |\nabla (\log \rhoDtp - c^{\Delta t, +}) |^{2}$, the inequality \eqref{EnergyDissipPerfSquare} 
can be seen as a temporally semi-discrete counterpart of the time-integral in the (formal) free energy identity \eqref{DecayEnergy1}.

We shall assume throughout the rest of the section that the test functions belong to the function space $L^{2}(0, T, W^{1, \infty}(\Omega))$. Our objective is to derive a (uniform in $\Delta t$) bound on the time-derivative of $\rhoDt$ in terms of the $L^{2}(0, T, W^{1, \infty}(\Omega))$ norm of the test function, followed by the use of Sobolev embedding, which will enable us to apply a variant of the Aubin--Lions compactness lemma, known as Dubinski\u{\i}'s compactness theorem.  From (\ref{WeakFormSemiDiscTime}) we have that, for a test function $\phi \in L^{2}(0, T, W^{1, \infty}(\Omega))$,
\begin{align*}
    \Bigg| \int_{0}^{T} \int_{\Omega} \frac{\partial \rhoDt}{\partial t} \phi \dx\dt \Bigg| \leq \underbrace{\Bigg| \int_{0}^{T} \int_{\Omega} \nabla \rhoDtp \cdot \nabla \phi \dx\dt \Bigg|}_{\circled{1}} + \underbrace{\Bigg| \int_{0}^{T} \int_{\Omega} \rhoDtp \nabla c^{\Delta t, +} \cdot \nabla \phi \dx\dt \Bigg|}_{\circled{2}}.
\end{align*}
We recall that we have 
\begin{subequations}\label{SmallProp1}
    \begin{equation}\label{Positivity}
    \rhoDtpm \geq 0 \quad \textrm{a.e.~on } \Omega \times (0, T), \quad \mbox{and}
    \end{equation}
    \begin{equation} \label{UnitDensityboundDt}
        \frac{1}{|\Omega|} \int_{\Omega} \rhoDtpm(t) = 1 \quad \textrm{for a.e. } t \in (0, T).
    \end{equation}
\end{subequations}
Moreover we notice that the bound on the first term on the left-hand side of (\ref{UniformBoundSpatialDer}) gives us a (uniform in $\Delta t$) bound on the $L^\infty(0,T;L^{1}(\Omega))$ norm of $\rhoDtp$, because the function $s \in [0,\infty) \mapsto G(s) \in [0,\infty)$ has superlinear growth as $s \rightarrow +\infty$. 

For the term $\circled{1}$ we then have by (\ref{UniformBoundSpatialDer}) and (\ref{SmallProp1}) that
\begin{align*}
    \circled{1} &\leq \int_{0}^{T} \|\nabla \phi\|_{L^{\infty}(\Omega)} \int_{\Omega} |\nabla \rhoDtp| \dx \dt= 2 \int_{0}^{T} \|\nabla \phi\|_{L^{\infty}(\Omega)} \int_{\Omega} \sqrt{\rhoDtp}\, \bigg| \nabla \sqrt{\rhoDtp} \bigg| \dx \dt\\
    &\leq 2 \int_{0}^{T} \|\nabla \phi\|_{L^{\infty}(\Omega)} \Bigg( \int_{\Omega} \rhoDtp \dx \Bigg)^{\frac{1}{2}} \Bigg( \int_{\Omega}\bigg| \nabla \sqrt{\rhoDtp} \bigg|^{2} \dx \Bigg)^{\frac{1}{2}} \dt  \leq C_{\ast} \Bigg( \int_{0}^{T} \|\nabla \phi\|_{L^{\infty}(\Omega)}^{2} \dt \Bigg)^{\frac{1}{2}}.
\end{align*}
For the term $\circled{2}$ we recall the bound on the second term on the left-hand side of (\ref{UniformBoundC}) together with (\ref{UniformBoundSpatialDer}); with a similar computation as for term $\circled{1}$ above we get 
\begin{equation*}
    \circled{2} \leq C_{\ast} \Bigg( \int_{0}^{T} \|\nabla \phi\|_{L^{\infty}(\Omega)}^{2} \dt \Bigg)^{\frac{1}{2}}.
\end{equation*}
By combining the bounds on the terms $\circled{1}$ and $\circled{2}$ and noting that by  the Sobolev embedding theorem $H^{d+1}(\Omega)$ is continuously embedded in $W^{1,\infty}(\Omega)$ for $d=2,3$, we then have that 
\begin{equation} \label{UniformBoundTimeDer}
\Bigg| \int_{0}^{T} \int_{\Omega} \frac{\partial \rhoDt}{\partial t} \phi \,\dx  \dt\Bigg| \leq C_{\ast} \Bigg( \int_{0}^{T} \|\nabla \phi\|_{L^{\infty}(\Omega)}^{2} \dt \Bigg)^{\frac{1}{2}}\leq C_{\ast} \Bigg( \int_{0}^{T} \|\phi\|_{H^{d+1}(\Omega)}^{2} \dt \Bigg)^{\frac{1}{2}}.
\end{equation} 

We are now ready to take the limit as $\Delta t \to 0_+$. Our main tool to this end will be Dubinski\u{\i}'s compactness theorem in seminormed sets, which we shall state in the next subsection; cf.  \cite{dubinskii1965weak} and \cite{barrett2012dubinskii}.

\subsection{Dubinski\u{\i}'s compactness theorem}

Let $\mathcal {A}$ be a linear space over the field $\mathbb{R}$ of real numbers, which will be referred as \textit{ambient space}, and suppose that $\mathcal{M}$ is a subset of $\mathcal{A}$ such that, 
\begin{equation} \label{PropSemiNorm}
    \textrm{for all } \varphi \in \mathcal{M} \textrm{ and } c \in \mathbb{R}_{\geq 0}, \quad c\,\varphi \in \mathcal{M}.
\end{equation}
Suppose further that each element $\varphi$ of a set $\mathcal{M}$ with property (\ref{PropSemiNorm}) is assigned a certain real number, denoted by $[\varphi]_{\mathcal{M}}$, such that 

\begin{itemize}
\item[$1.$] $[\varphi]_{\mathcal{M}} \geq 0$, and $[\varphi]_{\mathcal{M}}=0$ if, and only if, $\varphi=0$;
\item[$2.$] $[c \varphi]_{\mathcal{M}} = c[\varphi]_{\mathcal{M}}$ for all $c \in \mathbb{R}_{\geq 0}$.
\end{itemize}
A subset $\mathcal{M}$ satisfying \eqref{PropSemiNorm}, equipped with $[\cdot]$ satisfying these two properties is referred to as a \textit{seminormed set}. 

\begin{theo} [c.f. \cite{barrett2012dubinskii}] \label{Dubinsky}
    Suppose that $\mathcal{A}_{0}$ and $\mathcal{A}_{1}$ are Banach spaces, $\mathcal{A}_{0} \hookrightarrow \mathcal{A}_{1}$ (i.e., $\mathcal{A}_{0}$ is continuously embedded in $\mathcal{A}_{1}$), and $\mathcal{M}$ is a seminormed set contained in $\mathcal{A}_{0}$ such that $\mathcal{M}$ is compactly embedded in $\mathcal{A}_{0}$. Consider the set 
    \[ \mathcal{Y} := \Bigg \{ \varphi: [0, T] \to \mathcal{M}: [\varphi]_{L^{p}(0, T; \mathcal{M})} + \bigg\| \frac{\mathrm{d} \varphi}{\dt}\bigg\|_{L^{p_{1}}(0, T; \mathcal{A}_{1})} < \infty \Bigg \}, \]
    where $1\leq p \leq \infty$, $1\leq p_{1} \leq \infty$, $\| \cdot \|_{\mathcal{A}_{1}}$ is the norm of $\mathcal{A}_{1}$ and $\frac{\mathrm{d}\varphi}{\dt}$ is understood in the sense of $\mathcal{A}_{1}$-valued distributions on the open interval $(0, T)$. Then $\mathcal{Y}$, with 
    \[ [\varphi]_{\mathcal{Y}} := [\varphi]_{L^{p}(0, T; \mathcal{M})} + \bigg\| \frac{\mathrm{d} \varphi}{\dt} \bigg\|_{L^{p_{1}}(0, T; \mathcal{A}_{1})}, \]
    is a seminormed set in $L^{p}(0, T; \mathcal{A}_{0}) \cap W^{1, p_{1}}(0, T; \mathcal{A}_{1})$, and $\mathcal{Y}$ is compactly embedded in $L^{p}(0, T; \mathcal{A}_{0})$ if either $1 \leq p \leq \infty$ and $1< p_{1} < \infty$, or if $1 \leq p < \infty$ and $p_{1}=1$.
\end{theo}

In our current setting we shall apply  Dubinski\u{\i}'s theorem, as in \cite{barrett2011existence}, by selecting 
\[ \mathcal{A}_{0} = L^{1}(\Omega), \textrm{ with the usual Lebesgue norm } \|\varphi\|_{\mathcal{A}_{0}} := \int_{\Omega} |\varphi| \dx \]
and 
\[ \mathcal{M} = \Bigg \{ \phi \in \mathcal{A}_{0}: \varphi \geq 0 \textrm{ with } \int_{\Omega}\big| \nabla \sqrt{\varphi} \big|^{2} \dx < \infty \Bigg \};\]
and, for $\varphi \in \mathcal{M}$, we define 
\[ [\varphi]_{\mathcal{M}} := \|\varphi\|_{\mathcal{A}_{0}} + \int_{\Omega}\big| \nabla \sqrt{\varphi} \big|^{2} \dx.\]
It is a straightforward matter to confirm that (\ref{PropSemiNorm}) and properties 1. and 2. stated above hold, rendering $\mathcal{M}$ a seminormed subset of the ambient space $\mathcal{A}_{0}=L^1(\Omega)$. Finally we put 
\[ \mathcal{A}_{1} = H^{-\beta}(\Omega):= [H^{\beta}(\Omega)]'\quad \mbox{with $\beta = d + 1$}, \]
equipped with the dual norm $\|\varphi\|_{\mathcal{A}_{1}} := \|\varphi\|_{H^{-\beta}(\Omega)}$. The choice of $\mathcal{A}_{1}$ is motivated by the (continuous) Sobolev embedding 
\[ H^{\beta}(\Omega) \hookrightarrow W^{1, \infty}(\Omega).\]
In this way the right-hand side of (\ref{UniformBoundTimeDer}) can be further bounded by a constant multiple of $\|\phi\|_{L^{2}(0,T:H^{\beta}(\Omega))}$. The Banach spaces $\mathcal{A}_{0}$ and $\mathcal{A}_{1}$ and the seminormed set $\mathcal{M}$ then satisfy the conditions of Theorem \ref{Dubinsky}, with $\mathcal{M} \hookrightarrow \mathcal{A}_{0}$ compactly and $\mathcal{A}_{0} \hookrightarrow \mathcal{A}_{1}$
(cf. \cite{barrett2011existence}).\footnote{The continous embedding of $\mathcal{A}_0$ into $\mathcal{A}_1$ follows from the following abstract result, whose proof of straightforward. \textit{Let $X$ and $Y$ be normed linear spaces with (continuous) dual spaces $X'$ and $Y'$, respectively, and let $i : X \rightarrow Y'$ be a continuous embedding. Then the mapping $j : Y \rightarrow X'$, defined by $\langle j(y), x \rangle := \langle y , i(x) \rangle$, is continuous. Moreover, if $i(X)$ is weak-$\ast$ dense in $Y'$, then $j$ is injective, i.e., $j$ is a continuous embedding.} We take $X=H^\beta(\Omega)$ with $\beta=d+1$, $Y=\mathcal{A}_0=L^1(\Omega)$, and note that $Y'=L^\infty(\Omega)$.}

We shall also make use of the next lemma, whose proof via H\"older's inequality is straightforward and is therefore omitted.  

\begin{lemma} \label{DubLemmaGen}
    Suppose that a sequence $\{\varphi_{n}\}_{n\geq1}$ converges strongly in $L^{1}(0,T; L^{1}(\Omega))$ to a function $\varphi \in L^{1}(0,T; L^{1}(\Omega))$, and is bounded in $L^{\infty}(0,T; L^{\infty}(\Omega))$; namely, there exists a $K>0$ such that $\|\varphi_{n}\|_{L^{\infty}(0, T; L^{\infty}(\Omega))} \leq K$ for all $n\geq1$. Then $\varphi \in L^{p}(0,T; L^{q}(\Omega))$ for all $p, q \in [1, \infty)$ and the sequence $\{ \varphi_{n} \}_{n \geq 1}$ converges to $\varphi$ strongly  in $L^{p}(0, T; L^{q}(\Omega))$ for all $p, q \in [1, \infty)$.
\end{lemma}

\subsection{Passage to the limit $\Delta t \rightarrow 0_{+}$}

We begin by collecting the $\Delta t$-independent bounds (\ref{UniformBoundSpatialDer}), (\ref{UniformBoundSpatialDerPlus}), (\ref{UniformBoundTimeDer});  we have shown that there exists a constant $C_{\star}$, independent on $\Delta t$, such that 
\begin{align}\label{UniformBoundDer}
\begin{aligned}
    \sup_{t \in (0,T]} \int_{\Omega} G(\rhoDtp(t)) \dx & + \int_{0}^{T} \int_{\Omega} (\rhoDtp - \rhoDtm)^{2} \dx\dt \\
    &+ \int_{0}^{T} \int_{\Omega} \bigg|\nabla \sqrt{\rhoDtpm}\bigg|^{2} \dx\dt  + \int_{0}^{T}\bigg\| \frac{\partial \rhoDt}{\partial t} \bigg\|^{2}_{H^{-\beta}(\Omega)} \dt \leq C_{\star}.
\end{aligned}
\end{align}

In the proof of the next theorem we shall require a technical result concerning weakly continuous and weak-$\ast$ continuous functions. Suppose to this end that $X$ is a Banach space and let $X'$ be the (continuous) dual space of $X$. Let us denote by $C_{w}([0, T]; X)$ the set of all weakly continuous functions that map $[0, T]$ into $X$, namely the set of all functions $\eta \in L^{\infty}(0, T; X)$ such that $t \in [0,T] \mapsto \langle \varphi, \eta(t) \rangle_{X} \in \mathbb{R}$ is continuous on $[0, T]$ for all $\varphi \in X'$. Whenever $X$ has a predual $E$, we shall denote by $C_{w^{\ast}}([0, T]; X)$ the set of all functions $\eta \in L^{\infty}(0, T; X)$ such that $t \in [0, T] \mapsto \langle \eta(t), \zeta \rangle \in \mathbb{R}$ is continuous on $[0, T]$ for all $\zeta \in E$. We will use Lemma 3.1 of \cite{barrett2016existence}, that we cite here. 

\begin{lemma} \label{WeakContLemma}
    Suppose that $X$ and $E$ are Banach spaces. 
    \begin{itemize}
\item[(a)] Assume that the space $X$ is reflexive and is continuously embedded in a Banach space $Y$; then $L^{\infty}(0, T; X) \cap C_{w}([0, T]; Y) = C_{w}([0, T]; X)$.
\item[(b)] Assume that $X$ has a separable predual $E$ and the Banach space $Y$ has a predual $F$ such that $F$
 is continuously embedded in $E$; then, $L^{\infty}(0, T; X) \cap C_{w^{\ast}}([0, T]; Y) = C_{w^{\ast}}([0, T]; X)$.    \end{itemize}
\end{lemma}

We are now ready to pass to the limit $\Delta t \rightarrow 0_+$.
\begin{theo} \label{TheoConvLDt}
 There exists a subsequence of $\{ \rhoDtpm \}_{\Delta t>0}$ (not indicated) and a function $\widehat{\rho}$ such that 
\begin{equation*}
    \widehat{\rho} \in L^{\infty}(0, T, L^{\infty}(\Omega)) \cap H^{1}(0, T, H^{-\beta}(\Omega)), \quad \beta = d+1, 
\end{equation*}
with $\widehat{\rho} \geq 0$ almost everywhere on $\Omega\times (0,T)$ and $\frac{1}{|\Omega|} \int_{\Omega}\widehat{\rho}(x.t) \dx = 1$ for a.e.~$t \in [0, T]$, and $\widehat{c}(\cdot,t) \in \mathbb{H}^{s}(\Omega) \cap H^1_\ast(\Omega)$ defined as $-(-\Delta_{\mathrm{N}})^{s} \widehat{c}(\cdot,t) = \widehat{\rho}^{\ast}(\cdot,t)$ in $\Omega$  for $t \in (0,T]$ such that, for all $p \in [1, \infty)$, 
\begin{subequations}
    \begin{alignat}{2} \label{ConvL1}
\rho^{\Delta t(,\pm)} & \to \widehat{\rho} &&\qquad\text{strongly in } L^{p}(0,T;L^{1}(\Omega)), \\ \label{ConvL2}
\nabla \sqrt{\rho^{\Delta t(,\pm)}} & \to \nabla \sqrt{\widehat{\rho}} &&\qquad\text{weakly in } L^{2}(0,T;L^{2}(\Omega;\mathbb{R}^d)), \\
\label{ConvL3}
\frac{\partial \rhoDt}{\partial t} & \to \frac{\partial \widehat{\rho}}{\partial t} &&\qquad \text{weakly in } L^{2}(0,T; H^{-\beta}(\Omega)), \\
\label{ConvL3a}
\rho^{\Delta t(,\pm)} & \to \widehat{\rho} &&\qquad \text{weakly-$\ast$ in } L^{\infty}(0,T;L^{\infty}(\Omega)),\\
\label{ConvL5}
\rho^{\Delta t(,\pm)} & \to \widehat{\rho} &&\qquad \text{strongly in } L^{p}(0,T;L^{p}(\Omega)), \\
\label{ConvL6}
c^{\Delta t(,\pm)} & \to \widehat{c} &&\qquad \text{strongly in }  L^{p
}(0, T; \mathbb{H}^{s}(\Omega)), \\
\label{ConvL7}
\nabla c^{\Delta t, +} & \to \nabla \widehat{c} &&\qquad  \text{weakly in }  L^{2}(0, T; L^{2}(\Omega;\mathbb{R}^d)).
    \end{alignat}
The function $\widehat{\rho}$ is a global weak solution to the problem 
\begin{gather} \label{WeakSolutionEq}
    -\int_{0}^{T} \int_{\Omega} \widehat{\rho} \frac{\partial \phi}{\partial t}  \dx \dt+ \int_{0}^{T} \int_{\Omega}  \nabla \widehat{\rho} \cdot \nabla \phi  \dx \dt+ \int_{0}^{T} \int_{\Omega} \widehat{\rho} \, \nabla \widehat{c} \cdot \nabla \phi  \dx \dt = \int_{\Omega} \rho_{0} \;\phi|_{t=0} \dx \\
    \nonumber \text{for all } \phi \in W^{1,1}(0, T; H^{\beta}(\Omega)) \text{ such that } \phi(\cdot, T) = 0 . \end{gather}
In addition, the function $\widehat{\rho}$ is weak-$\ast$ continuous as a mapping from $[0, T]$ to $L^{\infty}(\Omega)$
and it is weakly continuous as a mapping from $[0,T]$ to $L^1(\Omega)$.  The energy functional $E(\cdot)$ defined in \eqref{Energy1} satisfies the inequality
    \begin{equation} \label{DecayEnergyWeakSol}
    E(\widehat\rho(t)) + \int_{0}^{t} \int_{\Omega} \bigg| 2 \nabla \sqrt{\widehat{\rho}} - \sqrt{\widehat{\rho}}\, \nabla \widehat{c} \bigg|^{2} \dx\, \dtau \leq E(\rho_{0}), \end{equation}
for a.e. $t \in [0, T]$.
    
\end{subequations}
 \begin{proof}
    We choose $\mathcal{A}_{0}, \mathcal{A}_{1}$ and $\mathcal{M}$ as in the discussion following the statement of  Dubinski\u{\i}'s theorem (Theorem \ref{Dubinsky}). Then, by taking $p=1$ and $p_{1}=2$ we have that 
    \begin{equation*}
        \Bigg \{ \varphi: [0, T] \to \mathcal{M}: [\varphi]_{L^{p}(0, T; \mathcal{M})} + \bigg\| \frac{\mathrm{d}\varphi}{\dt}\bigg\|_{L^{p_{1}}(0, T; \mathcal{A}_{1})} < \infty \Bigg \}
    \end{equation*}
    is compactly embedded in $L^{1}(0, T; \mathcal{A}_{0}) = L^{1}(0, T; L^{1}(\Omega))$. By using the bounds on the last two terms on the left-hand side of (\ref{UniformBoundDer}), we have that there exists a subsequence of $\{\rhoDt \}_{\Delta t >0}$ (not indicated), which converges strongly in $L^{1}(0, T; L^{1}(\Omega))$ to an element $\widehat{\rho} \in L^{1}(0, T; L^{1}(\Omega))$ as $\Delta t \to 0$. The strong convergence in $L^{1}(0, T; L^{1}(\Omega))$ ensures almost everywhere convergence of a subsequence (not indicated) to $\widehat{\rho}$.
    We have  by  (\ref{LinInterptime}) and the bound on the second term on the left-hand side of (\ref{UniformBoundDer}) that 
    \begin{equation} \label{IneqDeltaT}
        \Bigg( \int_{0}^{T} \int_{\Omega} | \rhoDt - \rhoDtp| \dx\dt \Bigg)^{2} \leq 
        T|\Omega| \int_{0}^{T} \int_{\Omega} (\rhoDtp - \rhoDtm)^{2} \dx\dt \leq C_{\star} T |\Omega|.
    \end{equation}
    We apply the triangle inequality in $L^{1}(0, T; L^{1}(\Omega))$, together with the strong convergence of the sequence $\{ \rhoDt \}_{\Delta t>0}$ in $L^{1}(0,T; L^{1}(\Omega))$ to deduce strong convergence in $L^{1}(0,T; L^{1}(\Omega))$ to the same element $\widehat{\rho} \in L^{1}(0,T; L^{1}(\Omega))$. ~The inequality  (\ref{IneqDeltaT}) also implies strong convergence of $\{ \rhoDtm \}_{\Delta t>0}$ in $L^{1}(0,T; L^{1}(\Omega))$ to $\widehat{\rho}$. This completes the proof of (\ref{ConvL1}) for $p=1$. 
    
    From (\ref{SmallProp1}) we have that 
    \begin{equation*}
        \|\rho^{\Delta t(,\pm)}(t)\|_{L^{1}(\Omega)} \leq |\Omega|
    \end{equation*}
    for a.e.~$t \in (0,T]$. This means that the sequences $\{ \rho^{\Delta t (,\pm)} \}_{\Delta t>0}$ are bounded in $L^{\infty}(0,T;L^{1}(\Omega))$. We can use Lemma \ref{DubLemmaGen} together with the strong convergence of the sequences to $\widehat{\rho}$ in $L^{1}(0,T;L^{1}(\Omega))$ to deduce that we have strong convergence in $L^{p}(0,T;L^{1}(\Omega))$ to the same limit for all values of $p \in [1, \infty)$. We have then completed the proof of (\ref{ConvL1}). 
    
    Because strong convergence in $L^{p}(0,T;L^{1}(\Omega))$ for $p \in [1, \infty)$ implies convergence almost everywhere on $\Omega \times (0,T])$ of a subsequence (not indicated), it follows from (\ref{Positivity}) that $\widehat{\rho} \geq 0$ a.e.~on $\Omega \times (0,T)$. By applying Fubini's theorem we have that 
    \begin{equation*}
        \int_{\Omega} |\rho^{\Delta t(,\pm)} - \widehat{\rho} | \dx \to 0 \quad \text{as } \Delta t \to 0_+ \text{ for 
 a.e. } t \in (0, T].
    \end{equation*}
Hence we have by (\ref{UnitDensityboundDt}) that
    \begin{equation*} 
       \frac{1}{|\Omega|} \int_{\Omega} \widehat{\rho}(t) \dx = 1 \quad \text{for a.e. } t \in (0, T].
    \end{equation*}
    
    We shall now prove (\ref{ConvL2}). We notice that since $|\sqrt{c_{1}} - \sqrt{c_{2}}| \leq \sqrt{|c_{1}-c_{2}|}$ for any two nonnegative real numbers $c_{1}$ and $c_{2}$, the strong convergence (\ref{ConvL1}) directly implies that, as $\Delta t \to 0_+$,
    \begin{equation*}
        \sqrt{\rho^{\Delta t(,\pm)}} \to \sqrt{\widehat{\rho}} \qquad  \text{strongly in } L^{p}(0, T; L^{2}(\Omega)) \quad \text{for all } p \in [1, \infty).
    \end{equation*}
The bound on the third term on the left-hand side of (\ref{UniformBoundDer}) implies the existence of a subsequence (not indicated) and an element $G \in L^{2}(0, T; L^{2}(\Omega;\mathbb{R}^d))$ such that 
    \begin{equation*}
        \nabla \sqrt{\rho^{\Delta t(,\pm)}} \to G \qquad \text{weakly in } L^{2}(0, T; L^{2}(\Omega;\mathbb{R}^d)).
    \end{equation*}
We have therefore that, for a test function $\eta \in C([0,T]; C^{\infty}_{0}(\Omega;\mathbb{R}^d))$,
    \begin{equation*}
       \int_{0}^{T} \int_{\Omega} G \cdot \eta \,\dx \dt \leftarrow -\int_{0}^{T} \int_{\Omega} \sqrt{\rho^{\Delta t(,\pm)}} \, \text{div} \, \eta \dx \dt \rightarrow - \int_{0}^{T} \int_{\Omega} \sqrt{\widehat{\rho}} \,\, \text{div} \, \eta \dx \dt.
    \end{equation*}
Therefore we have the equality 
    \begin{equation*}
        \int_{0}^{T}\int_{\Omega} G \cdot \eta\, \dx \dt =  - \int_{0}^{T}\int_{\Omega} \sqrt{\widehat{\rho}} \,\, \text{div} \, \eta \dx \dt,
    \end{equation*}
    and this means that $G$ is the distributional gradient of $\sqrt{\widehat{\rho}}$. As $G \in L^{2}(0, T; L^{2}(\Omega;\mathbb{R}^d))$ it follows that 
    \begin{equation*}
        G = \nabla \sqrt{\widehat{\rho}} \in L^{2}(0, T; L^{2}(\Omega;\mathbb{R}^d))
    \end{equation*}
    and hence (\ref{ConvL2}) has been proved. 
    
    The weak convergence result (\ref{ConvL3}) follows from the uniform bound on the last term on the left-hand side of (\ref{UniformBoundDer}) and the weak compactness of bounded balls in the Hilbert space $L^{2}(0,T;H^{-\beta}(\Omega))$. 

    The weak-$\ast$ convergence result (\ref{ConvL3a}) follows from the uniform bound \eqref{InfNormBoundDiscTime} stated in Lemma \ref{InfNormDecayDiscTime}, the weak-$\ast$ compactness of bounded balls in the Banach space $L^\infty(0,T;L^\infty(\Omega))$, and the uniqueness of the weak limit. 
    
    The strong convergence result (\ref{ConvL5}) is a direct consequence of (\ref{ConvL1}) and the uniform bound \eqref{InfNormBoundDiscTime} stated in Lemma \ref{InfNormDecayDiscTime}, using Lemma \ref{DubLemmaGen}.
    
    With the convergence (\ref{ConvL5}) in place, (\ref{ConvL6}) follows directly from (\ref{StabFracNeu}) and (\ref{FracPoincareUse}); that is,
    \begin{align*}
\int_0^T \| c^{\Delta t (,\pm)}(t) - \widehat{c}(t) \|_{{\mathbb{H}}^s_\ast(\Omega)}^p \dt \leq C \int_0^T \| \rhoDtpm(t) - \widehat{\rho}(t) \|_{L^{2}_\ast(\Omega)}^p \dt
    \end{align*}
    for all $p \in [1,\infty)$.
    
    Let us now consider the $L^{2}(0, T; L^{2}(\Omega;\mathbb{R}^d))$ norm of $\nabla c^{\Delta t, +}$. By the stability result (\ref{H1Stab})
    we have 
    \begin{align*} 
    \int_{0}^{T} \int_{\Omega} \big|\nabla c^{\Delta t, +} \big|^{2} \dx\dt &\leq C \int_{0}^{T} \int_{\Omega} \big| \nabla \rhoDtp \big|^{2} \dx\dt = 4C \int_{0}^{T} \int_{\Omega} \rhoDtp \bigg| \nabla \sqrt{\rhoDtp} \bigg|^{2} \dx\dt  \\
    &\leq 4C \sup_{t \in (0,T]} \|\rhoDtp\|_{L^{\infty}(\Omega)} \int_{0}^{T} \int_{\Omega}  \bigg| \nabla \sqrt{\rhoDtp} \bigg|^{2} \dx\dt, 
    \end{align*}
    and thanks to the uniform bound on the $L^{\infty}(0, T; L^{\infty}(\Omega))$ norm of $\rhoDtp$ from Lemma \ref{InfNormDecayDiscTime} and the bound on the third term on the left-hand side of (\ref{UniformBoundDer}) we have that the $L^{2}(0, T; L^{2}(\Omega;\mathbb{R}^d))$ norm of the gradient of $c^{\Delta t, +}$ is uniformly bounded and therefore $\nabla c^{\Delta t, +}$ converges weakly to an element, that we call $\mathcal{C}$, as $\Delta t \rightarrow 0_+$. By (\ref{ConvL6}) with $s=0$ we have that $\nabla c^{\Delta t, \pm}$ converges strongly to $\nabla \widehat{c}$ in $L^{2}(0, T, H^{-1}(\Omega;\mathbb{R}^d))$, but then $\nabla \widehat{c}$ must coincide with $\mathcal{C}$ 
    because of the uniqueness of the weak limit. 
    
    We now want to use the convergence results we have just proved to pass to the limit as $\Delta t \to 0_{+}$ in equation (\ref{WeakFormDisctime}). Throughout the argument we will consider test functions $\phi \in C^{1}([0, T]; C^{\infty}(\overline{\Omega}))$ such that $\phi(\cdot, T) = 0$. Note that the set of all such test functions is dense in the set of functions belonging to $W^{1,1}(0, T; H^{\beta}(\Omega))$ and vanishing at $t=T$ (in the sense of the trace theorem in $W^{1,1}(0,T)$), which is continuously embedded in $L^{2}(0, T; H^{\beta}(\Omega))$; therefore the use of such test functions is fully justified for the purposes of our argument. 
    
    We begin by considering the first term in (\ref{WeakFormSemiDiscTime}) and use integration by parts with respect to $t$ to deduce that
    \begin{equation*}
        \int_{0}^{T} \int_{\Omega} \frac{\partial \rhoDt}{\partial t} \phi \dx\dt = -\int_{0}^{T} \int_{\Omega} \rhoDt \frac{\partial \phi}{\partial t} \dx\dt - \int_{\Omega} \rho^{0} \; \phi|_{t=0} \dx
    \end{equation*}
    for all $\phi \in C^{1}([0,T], C^{\infty}(\overline{\Omega}))$ such that $\phi(\cdot , T)=0$. Moreover, 
\begin{equation*}
        \lim_{\Delta t \to 0_{+}} \rho^{0} = \rho_{0} \quad \text{weakly in } L^{2}(\Omega).
    \end{equation*}
    Therefore using (\ref{ConvL1}) we immediately have that, as $\Delta t \to 0_{+}$,
    \begin{equation*}
        \int_{0}^{T} \int_{\Omega} \frac{\partial \rhoLDt}{\partial t} \phi \dx\dt \to -\int_{0}^{T} \int_{\Omega} \widehat{\rho} \frac{\partial \phi}{\partial t} \dx\dt - \int_{\Omega} \rho_{0} \;\phi|_{t=0} \dx.
    \end{equation*}
    The second term in (\ref{WeakFormSemiDiscTime}) will be dealt with by decomposing it as follows:
\begin{align*} \int_{0}^{T} \int_{\Omega}  \nabla \rhoDtp \cdot \nabla \phi \dx\dt &= \underbrace{2\int_{0}^{T} \int_{\Omega} \bigg( \sqrt{\rhoDtp} - \sqrt{\widehat{\rho}} \bigg) \nabla \sqrt{\rhoDtp} \cdot \nabla \phi \dx\dt}_{\circled{1}} \\
    & \quad + \underbrace{2\int_{0}^{T} \int_{\Omega} \sqrt{\widehat{\rho}} \,\nabla \sqrt{\rhoDtp} \cdot \nabla \phi \dx\dt}_{\circled{2}}. 
    \end{align*}
For term $\circled{1}$ we have by H\"older's inequality that 
    \begin{align*}
        |\circled{1}| & \leq 2\int_{0}^{T} \bigg( \int_{\Omega}  \Big| \sqrt{\rhoDtp} - \sqrt{\widehat{\rho}} \Big|^{2} \dx \bigg)^{\frac{1}{2}} \bigg( \int_{\Omega} \Big| \nabla \sqrt{\rhoDtp} \Big|^{2} \dx \bigg)^{\frac{1}{2}} \|\nabla \phi\|_{L^{\infty}(\Omega)} \dt \\
        &\leq 2\,\Bigg( \int_{0}^{T} \int_{\Omega} \Big| \nabla \sqrt{\rhoDtp} \Big|^{2} \dx\dt  \Bigg)^{\frac{1}{2}} \Bigg( \int_{0}^{T} \Big\| \sqrt{\rhoDtp} - \sqrt{\widehat{\rho}} \Big\|_{L^{2}(\Omega)}^{r}  \dt\Bigg)^{\frac{1}{r}} \Bigg( \int_{0}^{T} \|\nabla \phi\|_{L^{\infty}(\Omega)}^{\frac{2r}{r-2}} \dt \Bigg)^{\frac{r-2}{2r}},
    \end{align*}
    where $r \in (2, \infty)$. Using the bound on the third term in the inequality stated in (\ref{UniformBoundDer}) we have 
    \begin{align*} |\circled{1}| &\leq 2C_{\star}^{\frac{1}{2}} \Big\| \sqrt{\rhoDtp} - \sqrt{\widehat{\rho}} \Big\|_{L^{r}(0, T; L^{2}(\Omega))} \|\nabla \phi\|_{L^{\frac{2r}{r-2}}(0, T; L^{\infty}(\Omega))} \\
    & \leq 2C_{\star}^{\frac{1}{2}} \|\rhoDtp - \widehat{\rho}\|_{L^{\frac{r}{2}}(0, T; L^{1}(\Omega))}^{\frac{1}{2}} \|\nabla \phi\|_{L^{\frac{2r}{r-2}}(0, T; L^{\infty}(\Omega))},   \end{align*}
    where we have used the elementary inequality $|\sqrt{c_{1}} - \sqrt{c_{2}}| \leq \sqrt{|c_{1} - c_{2}|}$ with $c_{1}, c_{2} \geq 0$. The first factor in the last line converges to 0 as $\Delta t \to 0_+$ thanks to  (\ref{ConvL1}), and therefore the term $\circled{1}$ converges to 0
    for every $\phi \in C^1([0,T];C^\infty(\overline\Omega))$ as $\Delta t \to 0_+$.
    
    Concerning the term $\circled{2}$, as $\widehat\rho \in L^\infty(0,T;L^\infty(\Omega))$ and $\widehat\rho \geq 0$ a.e.~on $\Omega \times (0,T)$, we have that  $\sqrt{\widehat{\rho}}\, \nabla \phi$ belongs to $L^{2}(0, T; L^{2}(\Omega))$, and the weak convergence result (\ref{ConvL2}) then directly implies that
    \[ \lim_{\Delta t \rightarrow 0_+} 2\int_{0}^{T} \int_{\Omega} \sqrt{\widehat{\rho}}\, \nabla \sqrt{\rhoDtp} \cdot \nabla \phi \,\dx\dt 
    = 2\int_{0}^{T} \int_{\Omega} \sqrt{\widehat{\rho}}\, \nabla \sqrt{\widehat{\rho}} \cdot \nabla \phi \,\dx\dt 
    = \int_{0}^{T} \int_{\Omega} \nabla \widehat{\rho} \cdot \nabla \phi \,\dx\dt. \]

    Let us consider the final term in (\ref{WeakFormSemiDiscTime}). By \eqref{ConvL5}, $\rho^{\Delta t,+} \rightarrow \widehat\rho$ strongly in $L^2(0,T;L^2(\Omega))$ and by \eqref{ConvL7} $\nabla c^{\Delta t, +} \rightharpoonup \nabla\widehat{c}$ weakly in $L^2(0,T;L^2(\Omega;\mathbb{R}^d))$; it follows that $\rho^{\Delta t,+} \nabla c^{\Delta t, +}\rightharpoonup \widehat\rho\,\nabla\widehat{c}$ weakly in $L^1(0,T;L^1(\Omega;\mathbb{R}^d))$, and therefore
    \begin{align*}
        \lim_{\Delta t \to 0_+} \int_{0}^{T} \int_{\Omega}\rho^{\Delta t, +} \nabla c^{\Delta t, +} \cdot \nabla \phi \dx\dt = \int_{0}^{T} \int_{\Omega} \widehat{\rho} \, \nabla \widehat{c} \cdot \nabla \phi \dx\dt
    \end{align*}
    for all $\phi \in C^1([0,T];C^\infty(\overline\Omega))$.

Putting these convergence results together, and noting that $C^1([0,T];C^\infty(\overline\Omega))$ is dense in the function space $W^{1,1}(0.T;H^\beta(\Omega))$ with $\beta=d+1$, we get (\ref{WeakSolutionEq}).

Let us now prove the weak-$\ast$ continuity of the function $\widehat{\rho}$ as a map from $[0, T]$ to $L^{\infty}(\Omega)$. 
As $\widehat\rho \in L^\infty(0,T;L^\infty(\Omega))$ and $\widehat\rho \in H^1(0,T; \mathbb{H}^{-\beta}(\Omega)) \subset C([0,T];\mathbb{H}^{-\beta}(\Omega))$, $\beta=d+1$, and $L^\infty(\Omega) = L^1(\Omega)'$, it follows from part (b) of Lemma \ref{WeakContLemma} with $X= L^\infty(\Omega)$, $Y = \mathbb{H}^{-\beta}(\Omega)$ and $E=L^1(\Omega)$ that $\widehat \rho \in C_{w^\ast}([0,T]; L^\infty(\Omega))$. Because $\Omega$ is bounded, $L^\infty(\Omega) \subset L^1(\Omega)$, and it then directly follows from this weak-$\ast$ continuity property of $\widehat\rho$ that $\widehat \rho \in C_w([0,T]; L^1(\Omega))$, i.e.,  $\widehat{\rho}$ is weakly continuous as a mapping from $[0,T]$ to $L^1(\Omega)$.

It remains to prove the energy inequality \eqref{DecayEnergyWeakSol}.
Thanks to \eqref{EnergyIneqSemidisc}, in the form \eqref{EnergyDissipPerfSquare}, and (\ref{IneqPropIn}) we have that, for all $t \in (0, T]$,
\begin{align}
\int_{\Omega} G(\rho^{\Delta t, +}(x, t)) \dx &-\frac{1}{2} \int_{\Omega}  (\rho^{\Delta t,+})^\ast(x,t) \, c^{\Delta t,+}(x,t) \dx + \int_{0}^{t} \int_{\Omega} \bigg| 2 \nabla \sqrt{\rhoDtp} - \sqrt{\rhoDtp} \nabla c^{\Delta t, +} \bigg|^{2} \dx \,\dtau \nonumber \\
&\leq \int_{\Omega}G(\rho^{0}(x))\dx  -\frac{1}{2} \int_{\Omega} 
(\rho^{0})^\ast(x) \, c^0(x) \dx  \label{energy2}\\
&\leq \int_{\Omega} G(\rho_{0}(x))  \dx - \frac{1}{2} \int_{\Omega} 
(\rho^{0})^\ast(x) \, c^0(x) \dx.\nonumber
\end{align}

Since $G$ is a nonnegative, continuous and convex function, by Proposition 7.7 in \cite{Santambrogio2015} and the strong convergence \eqref{ConvL5}
it follows that, for a.e.~$t \in (0,T]$,  
\begin{align} \label{DecayEnergyDt1}
    \int_{\Omega} G(\widehat{\rho}(x,t)) \dx  \leq \liminf_{\Delta t \to 0_+} \int_{\Omega} G(\rhoDtp(x,t)) \dx.
\end{align}

As $\rho^0=\rho^0(\Delta t)$ is uniformly bounded in $L^2(\Omega)$ with respect to $\Delta t>0$, it follows that $(\rho^0)^\ast$ is uniformly bounded in $L^2_\ast(\Omega)$, and therefore by the continuous embedding of $L^2_\ast(\Omega) = \mathbb{H}^0_\ast(\Omega)$ into $\mathbb{H}^{-s}_\ast(\Omega)$, $(\rho^0)^\ast$ is uniformly bounded in $\mathbb{H}^{-s}_\ast(\Omega)$. By Theorem \ref{FracPoiExUniq}, $c^0=-(-\Delta_{\mathrm{N}})^{-s}(\rho^0)^\ast$ is then uniformly bounded in $\mathbb{H}^s(\Omega)$, and has therefore a weakly convergent subsequence (not indicated) in $\mathbb{H}^s(\Omega)$ as $\Delta t \rightarrow 0_+$; because $\mathbb{H}^s(\Omega)$ is compactly embedded in $\mathbb{H}^0(\Omega)=L^2(\Omega)$, we can extract a strongly convergent subsequence (not indicated) in $L^2(\Omega)$ as $\Delta t \rightarrow 0_+$. It remains to identify the limit of this subsequence. As $(\rho^0)^\ast \rightharpoonup (\rho_0)^\ast$ weakly in $L^2_\ast(\Omega)$ as $\Delta t \rightarrow 0_+$, thanks to the linearity of the operator $(-\Delta_{\mathrm{N}})^{-s}$, we have that $c^0 \rightharpoonup c_0 := -(-\Delta_{\mathrm{N}})^{-s} (\rho_0)^\ast$ weakly in $\mathbb{H}^s(\Omega)$, and therefore, thanks to the uniqueness of the weak limit $c^0 \rightarrow c_0= -(-\Delta_{\mathrm{N}})^{-s}(\rho_0)^\ast$ strongly in $L^2(\Omega)$ as $\Delta t\rightarrow 0_+$. Hence, $(\rho^0)^\ast c^0 \rightharpoonup (\rho_0)^\ast c_0$ weakly in $L^1(\Omega)$ as $\Delta t \rightarrow 0_+$. Using this to pass to the limit in the second term on the right-hand side of the inequality \eqref{energy2} we have
that
\begin{align}\label{energy3} 
\lim_{\Delta t \to 0_+}\frac{1}{2} \int_{\Omega}  (\rho^{0})^\ast(x) \, c^{0}(x) \dx = \frac{1}{2} \int_{\Omega} 
(\rho_{0})^\ast(x) \, c_0(x) \dx.
\end{align}

Concerning the second term on the left-hand side of \eqref{energy2}, thanks to the strong convergence results \eqref{ConvL5} and \eqref{ConvL6}  we have that $(\rho^{\Delta t,+})^\ast c^{\Delta t,+} \rightarrow (\widehat\rho)^{\ast}\, \widehat{c}$, strongly in $L^1(0,T;L^1(\Omega))$, and therefore also, 
\[ \frac{1}{2}\int_\Omega (\rho^{\Delta t,+})^\ast(x,\cdot)\, c^{\Delta t,+}(x,\cdot) \dx \rightarrow \frac{1}{2}\int_\Omega (\widehat\rho)^{\ast}(x,\cdot)\, \widehat{c}(x,\cdot) \dx\]
in $L^1(0,T)$ as $\Delta t \rightarrow 0_+$. Thus, for a subsequence (not indicated), 
\begin{align}\label{energy4}
\frac{1}{2}\int_\Omega (\rho^{\Delta t,+})^\ast(x,t)\, c^{\Delta t,+}(x,t) \dx \rightarrow \frac{1}{2}\int_\Omega (\widehat\rho)^\ast(x,t)\, \widehat{c}(x,t) \dx
\end{align} 
for a.e.~$t \in (0,T]$ as $\Delta t \rightarrow 0_+$.

Finally, for the third term on the left-hand side of \eqref{energy2} we use Theorem 7.5 in \cite{FonsecaLeoni} for the functional 
\begin{equation*}
u \mapsto \mathcal{J}(u) = \int_{0}^{t} \int_{\Omega} |u|^{2} \dx\, \dtau,
\end{equation*}
with respect to the sequence $2 \nabla \sqrt{\rho^{\Delta t, +}} - \sqrt{\rho^{\Delta t, +}} \nabla c^{\Delta t, +} = 2 \nabla \sqrt{\rho^{\Delta t, +}} + \sqrt{\rho^{\Delta t, +}} \nabla (-\Delta_{\mathrm{N}})^{-s}(\rho^{\Delta t, +})^{\ast}$, which is weakly convergent in $L^{2}(0, T; L^{2}(\Omega; \mathbb{R}^{d}))$ (identified with $L^{2}((0, T) \times \Omega; \mathbb{R}^{2})$) to $2 \nabla \sqrt{\widehat{\rho}} - \sqrt{\widehat{\rho}}\, \nabla \widehat{c}$ by \eqref{ConvL2}, \eqref{ConvL5} and \eqref{ConvL7}. The functional $\mathcal{J}(u)$ is clearly convex and therefore we have
\begin{equation} \label{energy5}
 \int_{0}^{t} \int_{\Omega} \bigg| 2 \nabla \sqrt{\widehat{\rho}} - \sqrt{\widehat{\rho}} \nabla \widehat{c} \bigg|^{2} \dx\, \dtau \leq \liminf_{\Delta t \to 0_{+}} \int_{0}^{t} \int_{\Omega} \bigg| 2 \nabla \sqrt{\rhoDtp} - \sqrt{\rhoDtp} \nabla c^{\Delta t, +} \bigg|^{2} \dx\, \dtau.
\end{equation}

Therefore, by inserting \eqref{DecayEnergyDt1}, \eqref{energy3}, \eqref{energy4} and \eqref{energy5} into \eqref{energy2}, we deduce that
\begin{align*}
    \int_{\Omega} G(\widehat{\rho}(x, t)) \dx & - \frac{1}{2} \int_{\Omega} (\widehat{\rho})^{\ast}(x, t)\, \widehat{c}(x, t)  \dx + \int_{0}^{t} \int_{\Omega} \bigg| 2 \nabla \sqrt{\widehat{\rho}} - \sqrt{\widehat{\rho}}\, \nabla \widehat{c} \bigg|^{2} \dx\, \dtau\\
    &\leq \int_{\Omega} G(\rho_{0}(x))  \dx - \frac{1}{2} \int_{\Omega} 
(\rho_{0})^\ast(x) \, c_0(x) \dx\label{energytot}
\end{align*}
for a.e.~$t \in (0,T]$. Thus we have shown that 
\begin{equation*}
E(\widehat\rho(t)) + \int_{0}^{t} \int_{\Omega} \bigg| 2 \nabla \sqrt{\widehat{\rho}} - \sqrt{\widehat{\rho}} \nabla \widehat{c} \bigg|^{2} \dx \,\dtau \leq E(\rho_{0})
\end{equation*}
for a.e.~$t \in (0,T]$ and the final result (\ref{DecayEnergyWeakSol}) then directly follows.
\end{proof}
\end{theo}

\color{black}
\begin{remark}
    We note in passing that, unlike the nonlinear mobility model II mentioned in Remark \ref{model2}, model I stated in Remark \ref{NonLinProbRem}
    does not have a gradient flow structure  since the free energy \eqref{Energy1} does not have a discrete counterpart similar to \eqref{NonLinearMobProb}. Passage to the limit as $\Delta t \to 0$ in the nonlinear mobility model I therefore remains an open problem. 
\end{remark}
\color{black}

As a conclusion of our argument in this section, we shall prove that the function 
$\widehat{\rho}$, which we defined via the limiting procedure detailed in Theorem \ref{TheoConvLDt}, decays exponentially in time to the trivial solution of problem (\ref{MainFullProb}), and it does so at a rate that depends on the $L^{\infty}(\Omega)$ norm of the initial datum $\rho_{0}$. As a consequence of this we shall also be able to strengthen the bound asserted in (\ref{DecayEnergyWeakSol}); in particular we shall show that, provided that $\Omega$ is a bounded convex (and therefore Lipschitz) polytopal domain in $\mathbb{R}^d$, $d \in \{2,3\}$,
and $s \in (1/2,1)$, the function $t \in [0,T] \mapsto E(\widehat{\rho}(t)) \in \mathbb{R}_{\geq 0}$ exhibits exponential decay.

\begin{theo}\label{th4.7}
   Suppose that the assumptions of Theorem \ref{TheoConvLDt} hold. Then, the weak solution $\widehat\rho$ we have constructed exhibits the following exponential decay properties.  
   \begin{itemize}
       \item[(a)] For any  $T>0$, the functional $G(\cdot)$, defined in \eqref{eq:G}, exhibits exponential decay in the sense that
   \begin{equation}\label{eq:decay-a}
        \int_{\Omega} G(\widehat{\rho}(T)) \dx \leq \mathrm{e}^{-\frac{2 T}{C_{\Omega} \|\rho_{0}\|_{L^{\infty}(\Omega)}}}  \int_{\Omega} G(\rho_{0}) \dx,
   \end{equation}
  where $C_{\Omega}=1/\lambda_1$, and $\lambda_1$ is the smallest positive eigenvalue of the Neumann Laplacian on the domain $\Omega$.
  \item[(b)] The function $\widehat\rho$ satisfies the inequality
\begin{equation}\label{eq:decay}
       \frac{1}{2|\Omega|} \| \widehat{\rho}(\cdot,T) - 1 \|_{L^{1}(\Omega)}^{2} \leq \mathrm{e}^{-\frac{2 T}{C_{\Omega} \|\rho_{0}\|_{L^{\infty}(\Omega)}}}  \int_{\Omega} E(\rho_{0}) \dx.
   \end{equation}
\item[(c)] Assuming in addition that $\Omega$ is convex and that $s \in (1/2, 1)$, the energy functional $E(\cdot)$, defined in \eqref{Energy1}, exhibits exponential decay for any $T>0$, that is,
\begin{equation}\label{eq:decay-c}
        \int_{\Omega} E(\widehat{\rho}(T)) \dx \leq \mathrm{e}^{-\frac{8 T}{3C_{\Omega} \|\rho_{0}\|_{L^{\infty}(\Omega)}}}  \int_{\Omega} E(\rho_{0}) \dx.
   \end{equation}
\end{itemize}

\begin{proof}
\textit{(a)} We start by taking $t=t_{1}=\Delta t$ and $0=t_{0}$ in \eqref{RhoIneqEnergy1}. We have that, by definition, $\rho^{0} = \rhoDtm(t_{1})$. We adopt the notation $t_{-1} = -\infty$ in order to write $\rhoDtp(t_{0})$ instead of $\rhoDtm(t_{1})$. Therefore we have that 
    \begin{align}
            \int_{\Omega} G(\rhoDtp(t_{1})) \dx &+ \frac{1}{2 \Delta t \| \rho_{0}\|_{L^{\infty}(\Omega)}} \int_{t_{0}}^{t_{1}} \int_{\Omega} (\rhoDtp - \rhoDtm)^{2} \dx\,\dtau + 4\int_{t_{0}}^{t_{1}} \int_{\Omega} \bigg|\nabla \sqrt{\rhoDtp}\bigg|^{2} \dx \,\dtau \nonumber \\
    &\quad - \int_{t_{0}}^{t_{1}} \int_{\Omega} \nabla c^{\Delta t, +} \cdot \nabla \rhoDtp \dx\,\dtau \nonumber \\ &\leq \int_{\Omega} G(\rhoDtp(t_{0})) \dx \label{StabEstAlpha1}.
    \end{align}
We can now repeat the argument that gave us (\ref{StabEstAlpha1}) on each interval $[t_{n-1}, t_{n}]$ for $n=1, \ldots, N$, resulting in 
    \begin{align}
            \int_{\Omega} \!G(\rhoDtp(t_{n})) \dx &+ \frac{1}{2 \Delta t \| \rho_{0}\|_{L^{\infty}(\Omega)}} \int_{t_{n-1}}^{t_{n}} \!\int_{\Omega} (\rhoDtp \! - \rhoDtm)^{2} \dx\,\dtau \!+4\!\int_{t_{n-1}}^{t_{n}}\!\int_{\Omega} \bigg|\nabla \sqrt{\rhoDtp}\bigg|^{2} \dx \,\dtau \nonumber \\
    &\quad - \int_{t_{n-1}}^{t_{n}} \int_{\Omega} \nabla c^{\Delta t, +} \cdot \nabla \rhoDtp \dx \,\dtau  
    \leq \int_{\Omega} G(\rhoDtp(t_{n-1}) + \alpha) \dx \label{StabEstAlphaGen}.
    \end{align}
We discard the nonnegative second term on the left-hand side of (\ref{StabEstAlphaGen}) to find that
\begin{align}
            \int_{\Omega} G(\rhoDtp(t_{n})) \dx & + 4\int_{t_{n-1}}^{t_{n}} \int_{\Omega} \bigg|\nabla \sqrt{\rhoDtp}\bigg|^{2} \dx \, \dtau \nonumber - \int_{t_{n-1}}^{t_{n}} \int_{\Omega} \nabla c^{\Delta t, +} \cdot \nabla \rhoDtp \dx\,\dtau \nonumber \\ &\leq \int_{\Omega} G(\rhoDtp(t_{n-1})) \dx \label{StabEstAlphaGenPlus}.
    \end{align}
    
    For $s \in [t_{n-1}, t_{n}]$, $n = 1, \ldots, N$, by Poincar\'{e}'s inequality we have that 
\begin{equation*}
        \int_{\Omega} \rhoDtp(s) \log \rhoDtp(s) \dx \leq \int_{\Omega} (\rhoDtp(s))^{2} \dx \leq C_{\Omega} \int_{\Omega} |\nabla \rhoDtp(s)|^{2} \dx,  
    \end{equation*}
    with a positive constant $C_{\Omega}>0$ depending only on the domain $\Omega$, and the bound on the $L^{\infty}(\Omega)$ norm of $\rhoLDtp$ in Lemma \ref{InfNormDecayDiscTime} gives us
    \begin{align*}
        \int_{\Omega} \rhoDtp(s) \log \rhoDtp(s) \dx &\leq C_{\Omega} \int_{\Omega} \rhoDtp(s) \bigg|\nabla \sqrt{\rhoDtp(s)} \bigg|^{2} \dx \nonumber  \leq C_{\Omega} \| \rho_{0} \|_{L^{\infty}(\Omega)} \int_{\Omega} \bigg|\nabla \sqrt{\rhoDtp(s)} \bigg|^{2} \dx.
    \end{align*}
We note in passing that $C_\Omega=1/\lambda_1$, the reciprocal of the smallest positive eigenvalue $\lambda_1$ of the Neumann Laplacian on $\Omega$.
As $s \log s = G(s) - (1-s)$, we can rewrite the above inequality as
    \begin{align*}
        \int_{\Omega} G(\rhoDtp(s)) \dx &\leq C_{\Omega} \| \rho_{0} \|_{L^{\infty}(\Omega)} \int_{\Omega} \bigg|\nabla \sqrt{\rhoDtp(s)} \bigg|^{2} \dx + \int_{\Omega} (1 - \rhoDtp(s)) \dx
    \end{align*}
    for $s \in (t_{n-1}, t_{n}]$, $n = 1, \ldots, N$. 
    Thanks to the fact that, by definition, $\rhoDtp$ is constant on the interval $(t_{n-1}, t_{n}]$ and \eqref{UnitDensityboundDt} holds, this implies that 
\begin{align}
         \Delta t \int_{\Omega} G(\rhoDtp(t_{n}))  \dx & \leq C_{\Omega} \| \rho_{0} \|_{L^{\infty}(\Omega)} \int_{t_{n-1}}^{t_{n}} \int_{\Omega} \bigg|\nabla \sqrt{\rhoDtp(s)} \bigg|^{2} \dx \, \textrm{d}s\label{StabEstAlphaHelp1}
    \end{align}
for $n = 1, \ldots, N$.  Using \eqref{StabEstAlphaGenPlus} and \eqref{StabEstAlphaHelp1} we then have the inequality
\begin{align*}
    & \bigg(1 + \frac{ 2 \Delta t}{C_{\Omega} \| \rho_{0} \|_{L^{\infty}(\Omega)}} \bigg) \int_{\Omega} G(\rhoDtp(t_{n})) \dx - \int_{t_{n-1}}^{t_{n}} \int_{\Omega} \nabla c^{\Delta t, +} \cdot \nabla \rhoDtp \dx\,\dtau  \leq  \int_{\Omega} G(\rhoDtp(t_{n-1})) \dx. 
\end{align*}
For the sake of simplicity of the exposition we introduce the following notation:
\begin{align*}
    \gamma &:= \frac{2}{C_{\Omega}\| \rho_{0} \|_{L^{\infty}(\Omega)}}, \qquad
    A_{n} :=  \int_{\Omega} G(\rhoDtp(t_{n})) \dx, \qquad B_{n}:= \int_{t_{n-1}}^{t_{n}} \int_{\Omega} \nabla c^{\Delta t, +} \cdot \nabla \rhoDtp \dx\,\dtau. 
\end{align*}
We have that $A_{n}$ is nonnegative. In terms of the new notation \eqref{StabEstAlphaComb} can be rewritten as 
\begin{equation*}
(1+\gamma \Delta t) A_{n} - B_{n} \leq A_{n-1}, \qquad n = 1, \ldots, N, 
\end{equation*}
and therefore by induction we have that
\begin{equation*}
    A_{n} - (1 + \gamma \Delta t)^{-n} \sum_{j=1}^{n} B_{j} \leq (1 + \gamma \Delta t)^{-n} A_{0}, \qquad n= 1, \ldots, N. 
\end{equation*}
In particular, with $n = N$ we have 
\begin{align*}
    &\int_{\Omega} G(\rhoDtp(T)) \dx - (1+ \gamma \Delta t)^{-N} \int_{0}^{T} \int_{\Omega} \nabla c^{\Delta t, +} \cdot \nabla \rhoDtp \dx\,\dtau \leq (1+ \gamma \Delta t)^{-N} \int_{\Omega} G(\rhoDtp(0) ) \dx . 
\end{align*}
We recall that $\rhoDtp(0) = \rho^{0}$ and thus we arrive at the following inequality:  
\begin{align}
     \int_{\Omega} G(\rhoDtp(T)) \dx &- \bigg( 1+ \frac{2\Delta t}{C_\Omega\|\rho_{0}\|_{L^{\infty}(\Omega)}} \bigg)^{-N} \int_{0}^{T} \int_{\Omega}  \nabla c^{\Delta t, +} \cdot \nabla \rhoDtp \dx\,\dtau \nonumber \\
    &\leq \bigg( 1+ \frac{2\Delta t}{C_\Omega\|\rho_{0}\|_{L^{\infty}(\Omega)}} \bigg)^{-N} \int_{\Omega} G(\rho^{0}) \dx.
    \label{StabEstHelp1}
\end{align}

We first discard the (nonnegative) second term on the left-hand side of (\ref{StabEstHelp1}) recalling (\ref{NegativeProdTime}). It then remains to pass to the limit in the resulting inequality as $\Delta t \to 0_{+}$. Before doing so, we need make a few additional technical observations.

Thanks to \eqref{ConvL1} we have that, on passing to a subsequence, $\rhoDt(T)$ converges weakly to $\widehat{\rho}(T)$ in $H^{-\beta}(\Omega)$ with $\beta = d+1$ as $\Delta t \to 0_{+}$. Noting the bound on the first term of \eqref{UniformBoundDer} and that $G(s)/s \to +\infty$ as $s \to +\infty$, we deduce from de la Vall\'{e}e-Poussin's theorem that the family $\{ \rhoDt(T) \}_{\Delta t > 0}$ is uniformly integrable in $L^{1}(\Omega)$. Hence, by the Dunford--Pettis theorem, $\{ \rhoDt(T) \}_{\Delta t > 0}$ is weakly relatively compact in $L^{1}(\Omega)$. Consequently, upon extraction of a further subsequence (not indicated), $\rhoDt(T)$ converges weakly in $L^{1}(\Omega)$ to some limit as $\Delta t \rightarrow 0_+$; however the uniqueness of the weak limit together with the weak convergence of the sequence in $H^{-\beta}(\Omega)$, $\beta = d+1$, implies that the sequence $\rhoDt(T)$ converges to $\widehat{\rho}(T)$ weakly in $L^{1}(\Omega)$, because $L^1(\Omega)$ is continuously embedded in $H^{-\beta}(\Omega)$ for $\beta=d+1$. Finally, since the function $G$ is continuous and convex, by applying Tonelli's semicontinuity theorem in $L^{1}(\Omega)$,
\begin{equation} \label{TonelliHelp}
    \int_{\Omega} G(\widehat{\rho}(x,T)) \dx \leq \liminf_{\Delta t \to 0_+} \int_{\Omega} G(\rhoDtp(x,T)) \dx, 
\end{equation}
by noticing that $\rhoDtp(x,T) = \rhoDt(x, T)$ for a.e.~$x \in \Omega$. We can now pass to the limit as $\Delta t \to 0_+$ (with $N = T/\Delta t \to +\infty$) in \eqref{StabEstHelp1}, using \eqref{ConvL1} and \eqref{TonelliHelp}; hence we find that

\begin{equation} \label{ExpDecay1}
    \int_{\Omega} G(\widehat{\rho}(x,T)) \dx \leq \mathrm{e}^{-\frac{2 T}{C_{\Omega} \|\rho_{0}\|_{L^{\infty}(\Omega)}}}  \int_{\Omega} G(\rho_{0}(x)) \dx,
\end{equation}
which proves part (a).

\textit{(b)} The Csisz\'{a}r--Kullback inequality \cite{unterreiter2000}, with the probability measure $\mathrm{d}\mu(x):=\frac{1}{|\Omega|} \dx$ on $\Omega$, yields that 
\[ \int_\Omega |\widehat{\rho}(x,T) - 1| \dd \mu(x) \leq  \left(2 \int_\Omega G(\widehat{\rho}(x,T)) \dd \mu(x) \right)^{\frac{1}{2}},\]
thanks to the fact that $\widehat{\rho}(x,t)\geq 0$ for a.e.~$x \in \Omega$ and $\int_\Omega \widehat\rho(x,T)\,\mathrm{d}\mu(x) = 1$, and therefore\begin{equation*}
    \| \widehat{\rho}(\cdot,T) - 1 \|_{L^{1}(\Omega)}^{2} \leq 2 |\Omega| \int_{\Omega} G(\widehat{\rho}(x,T)) \dx. 
\end{equation*}
Combining this with \eqref{ExpDecay1} finally gives the desired result \eqref{eq:decay} of part (b), by the fact that $G(\rho^{0}) \leq E(\rho^{0})$ for the energy functional $E(\cdot)$ defined in \eqref{Energy1}.

\textit{(c)} We conclude by proving part (c) with the additional assumptions that $\Omega$ is convex and that the fractional order $s$ lies in between $\frac{1}{2}$ and 1. To this end, we shall focus on the second term in the expression \eqref{Energy1} of the energy functional $E(\rho^{\Delta t, +})$, that is, on the interaction term between $\rho^{\Delta t, +}$ and $c^{\Delta t, +}$.

  We notice that, from \eqref{EnrgyIneqPart2}, taking $t = t_{1}$ and $0 = t_{0}$ as we did previously, we have 
    \begin{align}\begin{aligned}
    \frac{1}{2} \int_{\Omega} (-\Delta_{\mathrm{N}})^{\s} c^{\Delta t,+}(t_1) \, c^{\Delta t,+}(t_1) \dx &+ \int_{t_0}^{t_1} \int_{\Omega} \rhoDtp \Big| \nabla c^{\Delta t, +} \Big|^{2} \dx \, \dtau \\
    &\hspace{-17mm} -\int_{t_0}^{t_{1}} \int_{\Omega} \nabla \rhoDtp \cdot \nabla c^{\Delta t, +} \dx \, \dtau \leq \frac{1}{2} \int_{\Omega} (-\Delta_{\mathrm{N}})^{\s}c^{\Delta t,+}(t_0) \, c^{\Delta t,+}(t_0) \dx, \label{EnrgyIneqC1}
    \end{aligned}
    \end{align}
    and then, on each interval $[t_{n-1}, t_{n}]$ for $n=1,\ldots, 
 N$, 
 \begin{align}\begin{aligned}
    \frac{1}{2} \int_{\Omega} (-\Delta_{\mathrm{N}})^{\s} c^{\Delta t,+}(t_n) \, c^{\Delta t,+}(t_n) \dx &+ \int_{t_{n-1}}^{t_n} \int_{\Omega} \rhoDtp \Big| \nabla c^{\Delta t, +} \Big|^{2} \dx \, \dtau \\
    &\hspace{-17mm} -\int_{t_{n-1}}^{t_{n}} \int_{\Omega} \nabla \rhoDtp \cdot \nabla c^{\Delta t, +} \dx \, \dtau \leq \frac{1}{2} \int_{\Omega} (-\Delta_{\mathrm{N}})^{\s}c^{\Delta t,+}(t_{n-1}) \, c^{\Delta t,+}(t_{n-1}) \dx. \label{EnergyIneqC2}
    \end{aligned}
    \end{align}
For $s \in [t_{n-1}, t_{n}]$ with $n =1, \ldots, N$, by \eqref{StabFracNeu} and Poincar\'{e}'s inequality we have that
\begin{align*}
    \| c^{\Delta t, +}(s) \|_{\mathbb{H}^{s}_{\ast}(\Omega)}^{2} = \int_{\Omega} (- \Delta_{\mathrm{N}})^{\s} c^{\Delta t,+}(s) \, c^{\Delta t,+}(s) \dx &\leq \int_{\Omega} ((\rhoDtp)^{\ast}(s))^{2} \dx \leq C_{\Omega} \int_{\Omega} |\nabla (\rhoDtp)^{\ast}(s)|^{2} \dx \\ &\leq C_{\Omega} \| \rho_{0} \|_{L^{\infty}(\Omega)} \int_{\Omega} \bigg|\nabla \sqrt{\rhoDtp(s)} \bigg|^{2} \dx,
\end{align*}
and hence
\begin{equation} \label{EnergyIneqC3}
    \Delta t  \int_{\Omega} (-\Delta_{\mathrm{N}})^{\s} c^{\Delta t,+}(t_n) \, c^{\Delta t,+}(t_n) \dx \leq C_{\Omega} \| \rho_{0} \|_{L^{\infty}(\Omega)} \int_{t_{n-1}}^{t_{n}} \int_{\Omega} \bigg|\nabla \sqrt{\rhoDtp(s)} \bigg|^{2} \dx \, \textrm{d}s.
\end{equation}
We combine \eqref{StabEstAlphaHelp1} and  \eqref{EnergyIneqC3} to find the following inequality for the energy functional $E(\cdot)$:
\begin{align}
\Delta t  \int_{\Omega} E(\rhoDtp(t_{n})) \dx &= \Delta t \int_{\Omega} G(\rhoDtp(t_{n})) \dx + \frac{1}{2} \Delta t    \int_{\Omega} (-\Delta_{\mathrm{N}})^{\s} c^{\Delta t,+}(t_n) \, c^{\Delta t,+}(t_n) \dx \nonumber \\ 
&\leq \frac{3}{2} C_{\Omega} \| \rho_{0} \|_{L^{\infty}(\Omega)} \int_{t_{n-1}}^{t_{n}} \int_{\Omega} \bigg|\nabla \sqrt{\rhoDtp(s)} \bigg|^{2} \dx \, \textrm{d}s.\label{EnergyIneqDecay1}
\end{align}
By summing \eqref{StabEstAlphaGenPlus}, \eqref{EnergyIneqC2} we also have that 
\begin{align*}
     \int_{\Omega} E(\rhoDtp(t_{n})) \dx &+4\int_{t_{n-1}}^{t_{n}} \int_{\Omega} \bigg|\nabla \sqrt{\rhoDtp}\bigg|^{2} \dx \, \dtau + \int_{t_{n-1}}^{t_n} \int_{\Omega} \rhoDtp \Big| \nabla c^{\Delta t, +} \Big|^{2} \dx\,\dtau \nonumber \\ &- \int_{t_{n-1}}^{t_{n}} \int_{\Omega} \nabla c^{\Delta t, +} \cdot \nabla \rhoDtp \dx\,\dtau \leq  \int_{\Omega} E(\rhoDtp(t_{n-1})) \dx,
\end{align*}
and by \eqref{EnergyIneqDecay1} we deduce the following inequality:
\begin{align}
    \bigg(1 + \frac{ 8 \Delta t}{3C_{\Omega} \| \rho_{0} \|_{L^{\infty}(\Omega)}} \bigg) \int_{\Omega} E(\rhoDtp(t_{n})) \dx &+ \int_{t_{n-1}}^{t_n} \int_{\Omega} \rhoDtp \Big| \nabla c^{\Delta t, +} \Big|^{2} \dx \, \dtau \nonumber \\ &- \int_{t_{n-1}}^{t_{n}} \int_{\Omega} \nabla c^{\Delta t, +} \cdot \nabla \rhoDtp \dx\,\dtau \leq  \int_{\Omega} E(\rhoDtp(t_{n-1})) \dx. \label{StabEstAlphaComb}
\end{align}
Similarly as we did previously, we introduce the following auxiliary notation:
\begin{align*}
    \delta &:= \frac{8}{3C_{\Omega}\| \rho_{0} \|_{L^{\infty}(\Omega)}}, \qquad
    C_{n} :=  \int_{\Omega} E(\rhoDtp(t_{n})) \dx, \\
    D_{n}&:= \int_{t_{n-1}}^{t_n} \int_{\Omega} \rhoDtp \Big| \nabla c^{\Delta t, +} \Big|^{2} \dx \, \dtau - \int_{t_{n-1}}^{t_{n}} \int_{\Omega} \nabla c^{\Delta t, +} \cdot \nabla \rhoDtp \dx\,\dtau. 
\end{align*}
We have that $C_{n}$ is nonnegative. In terms of this new notation \eqref{StabEstAlphaComb} can be rewritten as 
\begin{equation*}
(1+\delta \Delta t) C_{n} + D_{n} \leq C_{n-1}, \qquad n = 1, \ldots, N, 
\end{equation*}
and therefore by induction we have that
\begin{equation*}
    C_{n} + (1 + \delta \Delta t)^{-n} \sum_{j=1}^{n} D_{j} \leq (1 + \delta \Delta t)^{-n} C_{0}, \qquad n= 1, \ldots, N. 
\end{equation*}
In particular, with $n = N$ we have 
\begin{align*}
    \int_{\Omega} E(\rhoDtp(T)) \dx &+(1+ \delta \Delta t)^{-N} \int_{0}^{T} \int_{\Omega} \rhoDtp \Big| \nabla c^{\Delta t, +} \Big|^{2} \dx \, \dtau  \\ &- (1+ \delta \Delta t)^{-N} \int_{0}^{T} \int_{\Omega} \nabla c^{\Delta t, +} \cdot \nabla \rhoDtp \dx\,\dtau \leq (1+ \delta \Delta t)^{-N} \int_{\Omega} E(\rhoDtp(0)) \dx . 
\end{align*}
We recall that $\rhoDtp(0) = \rho^{0}$ and thus we arrive at the following inequality:  
\begin{align}
    \int_{\Omega} E(\rhoDtp(T)) \dx &+\bigg( 1+ \frac{8\Delta t}{3C_\Omega\|\rho_{0}\|_{L^{\infty}(\Omega)}} \bigg)^{-N}\int_{0}^{T} \int_{\Omega} \rhoDtp \Big| \nabla c^{\Delta t, +} \Big|^{2} \dx \, \dtau  \nonumber \\  &- \bigg( 1+ \frac{8\Delta t}{3C_\Omega\|\rho_{0}\|_{L^{\infty}(\Omega)}} \bigg)^{-N} \int_{0}^{T} \int_{\Omega}  \nabla c^{\Delta t, +} \cdot \nabla \rhoDtp \dx\,\dtau  \nonumber \\ &\leq \bigg( 1+ \frac{8\Delta t}{3C_\Omega\|\rho_{0}\|_{L^{\infty}(\Omega)}} \bigg)^{-N} \int_{\Omega} E(\rho^{0}) \dx.
    \label{StabEstHelp1a}
\end{align}

We first discard the (nonnegative) second and third terms on the left-hand side of (\ref{StabEstHelp1a}) recalling (\ref{NegativeProdTime}). It then remains to pass to the limit in the resulting inequality as $\Delta t \to 0_{+}$. We use the argument that allowed us to pass to the limit before and prove \eqref{ExpDecay1}, together with additional considerations for the interaction term involved in the energy functional.

We begin by noting that 
\begin{align*}
(\rho^{\Delta t})^\ast(x,t) &:= \rho^{\Delta t}(x,t) - \frac{1}{|\Omega|}\int_\Omega \rho^{\Delta t}(\xi,t) \dd \xi = \rho^{\Delta t}(x,t) - \frac{1}{|\Omega|}\int_\Omega \rho^{\Delta t}(\xi,0) \dd \xi \\&= \rho^{\Delta t}(x,t) - \frac{1}{|\Omega|}\int_\Omega \rho^0(\xi) \dd \xi =  \rho^{\Delta t}(x,t) - \frac{1}{|\Omega|}\int_\Omega \rho_0(\xi) \dd \xi = \rho^{\Delta t}(x,t) - 1.
\end{align*}
Hence,  
\[ \frac{\partial (\rho^{\Delta t})^\ast}{\partial t}= \frac{\partial \rho^{\Delta t}}{\partial t}.\]
Because $c^{\Delta t} = -(-\Delta_{\mathrm{N}})^{-s} (\rho^{\Delta t})^\ast$ and $c^{\Delta t}$ and $(\rho^{\Delta t})^\ast$ are continuous piecewise linear functions of $t \in [0,T]$
on the partition $0=t_0 < t_1 < \cdots < t_N=T$, it follows from the linearity of the differential operator $-(-\Delta_{\mathrm{N}})^{-s}$ and its independence of $t$ that
\begin{equation*}
    \frac{\partial c^{\Delta t}}{\partial t} = -(-\Delta_{\mathrm{N}})^{-s} \frac{\partial (\rho^{\Delta t})^\ast}{\partial t}  = -(-\Delta_{\mathrm{N}})^{-s} \frac{\partial \rho^{\Delta t}}{\partial t}.
\end{equation*}
Thereby, for any function $\phi \in W^{1, \infty}(\Omega)$, we have that
\begin{align}
\bigg| \int_{0}^{T} \int_{\Omega}  \frac{\partial c^{\Delta t}}{\partial t} \phi \dx\dt \bigg| &=  \bigg| \int_{0}^{T} \int_{\Omega}  \frac{\partial c^{\Delta t}}{\partial t} \phi^{\ast} \dx\dt \bigg|= \bigg| \int_{0}^{T} \int_{\Omega}  (-\Delta_{\mathrm{N}})^{-s} \frac{\partial \rho^{\Delta t}}{\partial t}  \phi^{\ast} \dx\dt \bigg| \nonumber \\ 
&= \bigg| \int_{0}^{T} \int_{\Omega}  \frac{\partial \rho^{\Delta t}}{\partial t} (-\Delta_{\mathrm{N}})^{-s}  \phi^{\ast} \dx\dt \bigg| \leq C_{\ast} \bigg( \int_{0}^{T} \| (-\Delta_{\mathrm{N}})^{-s}  \phi^{\ast} \|^{2}_{W^{1, \infty}(\Omega)} \dt \bigg)^{\frac{1}{2}} \nonumber \\ 
&\leq C_{\ast} \bigg( \int_{0}^{T} \|  \phi^{\ast} \|^{2}_{L^{\infty}(\Omega)} \dt \bigg)^{\frac{1}{2}} \leq C_{\ast} \bigg( \int_{0}^{T} \|  \phi^{\ast} \|^{2}_{W^{1, \infty}(\Omega)} \dt \bigg)^{\frac{1}{2}}, \label{UniformBoundTimeDerHelp}
\end{align}
where in the first equality of the first line we used the fact that, by definition, $c^{\Delta t}$ has zero integral average over $\Omega$
and therefore the same is true of its time-derivative; in the equality appearing in the second line we used the selfadjointness of $(-\Delta_{\mathrm{N}})^{-s}$ on $L^2_\ast(\Omega)$, and in the inequality appearing in the second line we made use of the first inequality in \eqref{UniformBoundTimeDer};
and in the passage from the second to the third line we used Lemma \ref{W1InftyLemma}, which requires the additional assumption on the convexity of $\Omega$ and the restriction on the fractional order $s$ to the range $s \in (1/2, 1)$.

Hence, thanks to the fact that $H^2(\Omega)$ is continuously embedded in $L^{\infty}(\Omega)$ for $d=2, 3$,  we deduce from \eqref{UniformBoundTimeDerHelp} the following uniform bound on the time-derivative of  $ c^{\Delta t}$ in $L^{2}(0, T; H^{-2}(\Omega))$, where $H^{-2}(\Omega):=[H^2(\Omega)]'$:
\begin{equation*} \Bigg| \int_{0}^{T} \int_{\Omega} \frac{\partial c^{\Delta t}}{\partial t} \phi \,\dx  \dt\,\Bigg| \leq C_{\ast} \Bigg( \int_{0}^{T} \|\phi\|_{H^{2}(\Omega)}^{2} \dt \Bigg)^{\frac{1}{2}},
\end{equation*} 
for all $\phi \in H^{2}(\Omega)$, where $C_\ast$ is a positive constant, independent of $\phi$ and $\Delta t$. Using the weak compactness of bounded balls in the Hilbert space $L^{2}(0, T; H^{-2}(\Omega))$ we then have that 
\begin{equation*} 
\frac{\partial c^{\Delta t}}{\partial t} \to \frac{\partial \widehat{c}}{\partial t} \quad \text{weakly in } L^{2}(0, T; H^{-2}(\Omega)).
\end{equation*}

By Lemma \ref{W1InftyLemma}, combined with the fact that $\rho^{\Delta t} \in L^{\infty}(0, T; L^{\infty}(\Omega))$  and is uniformly bounded with respect to $\Delta t$ in the norm of $L^{\infty}(0, T; L^{\infty}(\Omega))$ and the fact that $\widehat{\rho} \in L^{\infty}(0, T; L^{\infty}(\Omega))$ we have, respectively, that $c^{\Delta t} \in L^{\infty}(0, T; W^{1,\infty}(\Omega))$ and is uniformly bounded with respect to $\Delta t$ in the norm of $L^{\infty}(0, T; W^{1,\infty}(\Omega))$ and $\widehat{c} \in L^{\infty}(0, T; W^{1,\infty}(\Omega))$. Hence by the Aubin--Lions lemma $\{c^{\Delta t}\}_{\Delta t>0}$ is compact in $C([0,T];X)$ for any Banach space $X$ such that $W^{1,\infty}(\Omega) \hookrightarrow \hookrightarrow X \hookrightarrow H^{-2}(\Omega)$; in particular, (upon extraction of a subsequence, not indicated) $c^{\Delta t}$ converges to $\widehat{c}$ strongly in $C([0,T];L^\infty(\Omega))$. As $\widehat{\rho} \in C_{w^\ast}([0,T];L^\infty(\Omega))$, it then follows that the product $\widehat{\rho}\, \widehat{c}$ belongs to $ C_{w^\ast}([0,T];L^\infty(\Omega))$, and therefore also to $C_w([0,T];L^1(\Omega))$; in particular, 
\[ t \in [0,T] \mapsto \int_\Omega \widehat{\rho}(x,t) \, \widehat{c}(x,t) \dx \in \mathbb{R}\]
is a continuous function. Furthermore, because $\rho^{\Delta t}(T) \rightharpoonup \widehat{\rho}(T)$ weakly in $L^1(\Omega)$ and $c^{\Delta t}(T) \rightarrow \widehat{c}(T)$ strongly in $L^\infty(\Omega)$ as $\Delta t \rightarrow 0_+$, it follows that $\rho^{\Delta t}(T)\, c^{\Delta t}(T)  \rightharpoonup \widehat{\rho}(T)\, \widehat{c}(T)$ weakly in $L^1(\Omega)$ as $\Delta t \rightarrow 0_+$. Thus we deduce that
\begin{align}\label{CHelp1}\lim_{\Delta t \rightarrow 0_+} \int_\Omega \rho^{\Delta t}(x,T) \, c^{\Delta t}(x,T)\dx  = \int_\Omega \widehat{\rho}(x,T)\, \widehat{c}(x,T)\dx.
\end{align}

We can now pass to the limit as $\Delta t \to 0_+$ (with $N = T/\Delta t \to +\infty$) in \eqref{StabEstHelp1}, using \eqref{TonelliHelp}, and \eqref{CHelp1}. Hence we find that
\begin{equation*}
    \int_{\Omega} E(\widehat{\rho}(x,T)) \dx \leq \mathrm{e}^{-\frac{8 T}{3C_{\Omega} \|\rho_{0}\|_{L^{\infty}(\Omega)}}}  \int_{\Omega} E(\rho_{0}(x)) \dx. 
\end{equation*}
This completes the proof.     
\end{proof}
\end{theo}
\color{black}

 \section{ \color{black}Numerical experiments\normalcolor} \label{Sect:5}

In the preceding sections we have confined ourselves to theoretical considerations: the construction and convergence analysis of a `structure-preserving' numerical scheme for the fractional porous medium equation. 
\color{black}
The implementation of the numerical method introduced here, as well as associated computational aspects, are the subject of the paper \cite{carrillosuli2024}, where we describe an efficient technique based on rational approximation, which is capable of providing accurate approximations to the finite element counterpart of the spectral fractional Neumann Laplacian featuring in our numerical method, without explicit knowledge of its eigenvalues and eigenfunctions. 

For details on aspects of our work related to scientific computing and a variety of numerical experiments, we refer the reader to \cite{carrillosuli2024}. We shall confine ourselves here to numerical simulations that validate the proposed method and indicate some open problems associated with it.
\color{black}

\subsection{Numerical validation}

\color{black} In \cite{BilerGrzegorz2015} the authors identified an explicit steady state solution to equation \eqref{MainProb} in self-similar variables, namely the function $\Phi: \mathbb{R}^d \to \mathbb{R}$ given by 
\begin{equation} \label{BarenblattProfile}
    \Phi(y) = k_{s,d}(1-|y|^2)^s_+, \quad k_{s,d} = \frac{d \Gamma(d/2)}{(d+2s)4^s \Gamma(2-s) \Gamma(d/2 + 1-s)},
\end{equation}
which is a solution of the elliptic equation
\begin{equation*}
    -\lambda \nabla \cdot(y \Phi) = \nabla \cdot (\Phi \nabla (-\Delta)^{-s}\Phi), \quad y= \frac{t}{x^\lambda}, \quad \lambda = \frac{1}{d+2-2s}.
\end{equation*}
The equation lacks the parabolic regularization of the problem \eqref{FullProbWSpace} that we have considered, which is also used in the proof of existence of a solution to our finite element method. Instead of the elliptic problem stated above, we shall therefore consider the following parabolic problem:
\begin{equation} \label{FullProbViscosity}
\left \{
\begin{aligned}
&\frac{\partial \rho_\varepsilon}{\partial t} = \varepsilon \Delta \rho_\varepsilon - \nabla \cdot (\rho_\varepsilon \nabla c_\varepsilon) + \lambda \nabla \cdot (y\rho_\varepsilon), \\
& - (-\Delta)^{\s} c_\varepsilon = \rho_\varepsilon.
\end{aligned}
\right.
\end{equation}
for small positive values of $\varepsilon \ll 1$.

\begin{figure} 
\begin{subfigure}{0.48 \textwidth}
\centering
\includegraphics[width = \textwidth]{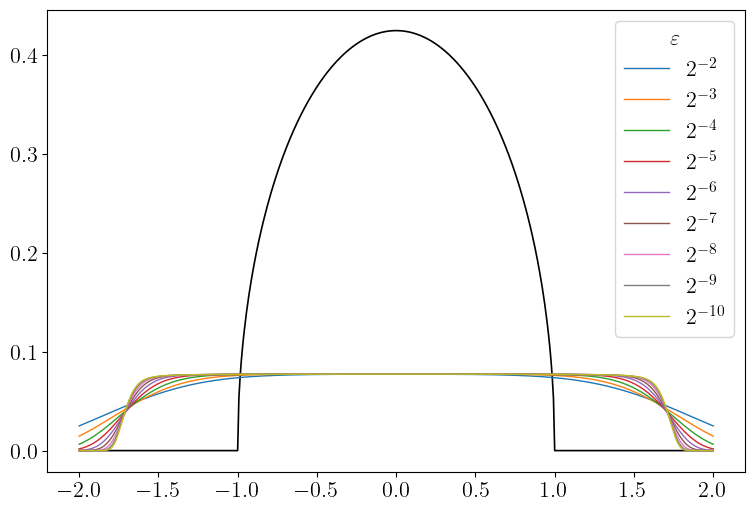}
\caption{$t=0.50$}
\end{subfigure}
\hfill
\begin{subfigure}{0.48 \textwidth}
\centering
\includegraphics[width = \textwidth]{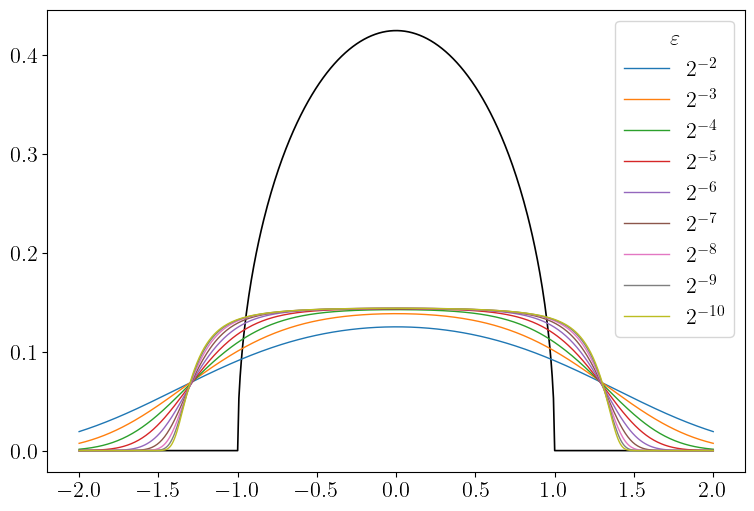}
\caption{$t=1.50$}
\end{subfigure}
\hfill \begin{subfigure}{0.48 \textwidth}
\centering
\includegraphics[width = \textwidth]{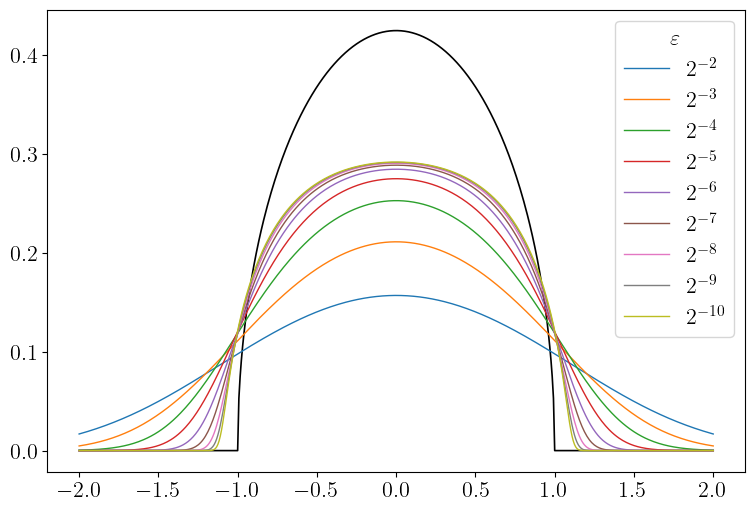}
\caption{$t=3.0$}
\end{subfigure}
\hfill
\begin{subfigure}{0.48 \textwidth}
\centering
\includegraphics[width = \textwidth]{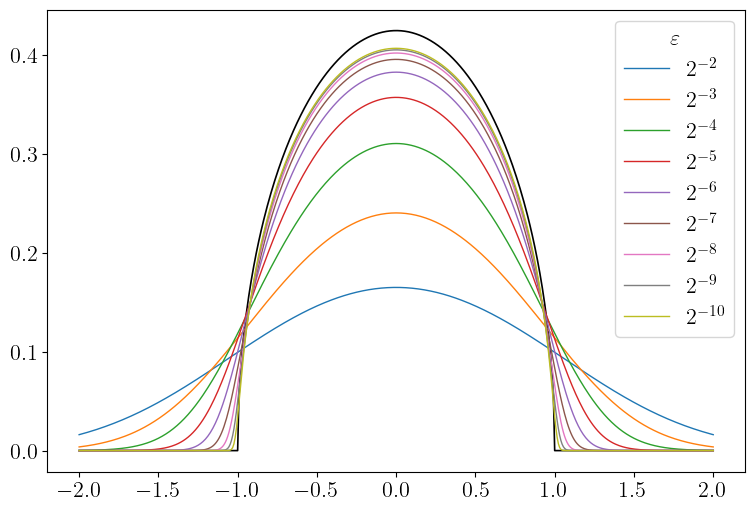}
\caption{$t=5.0$}
\end{subfigure}
\hfill
\centering
\begin{subfigure}{0.75 \textwidth}
\centering
\includegraphics[width = \textwidth]{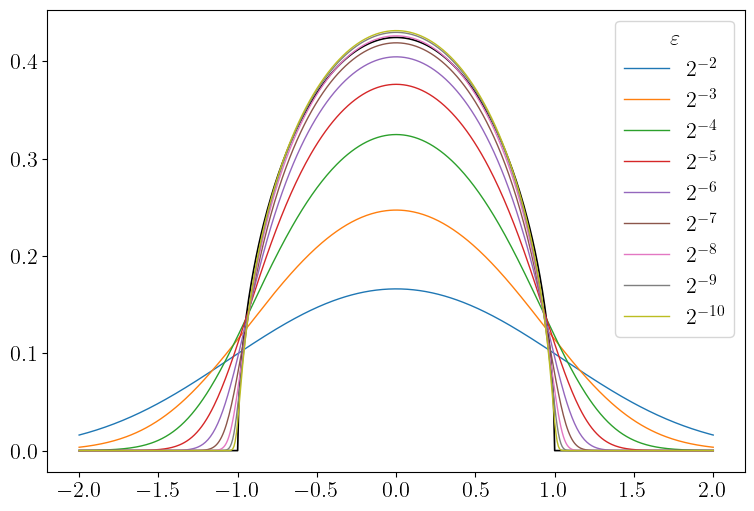}
\caption{$t=13.0$}
\end{subfigure}

\caption{Numerical solutions to the fractional porous medium equation in self-similar variables \eqref{FullProbViscosity}. The simulations were performed for the fractional order $s=0.5$ with different values of the parameter $\varepsilon = 2^{-m}$, $m=2, \dots, 10$, controlling the parabolic regularization. The computational domain is $\Omega = (-2,2)^2$, $\rho_0(x) \propto \frac{1}{|\Omega|}$; a uniform triangulation based on $2^9 \times 2^9$ equally 
 spaced nodes was used with $\Delta t = 0.01$. The black reference curve is the cross-section of the function \eqref{BarenblattProfile} at $x_2=0$.}
\label{FigsimulationsTime}
\end{figure}

\begin{figure} 
\begin{subfigure}{0.48 \textwidth}
\centering
\includegraphics[width = \textwidth]{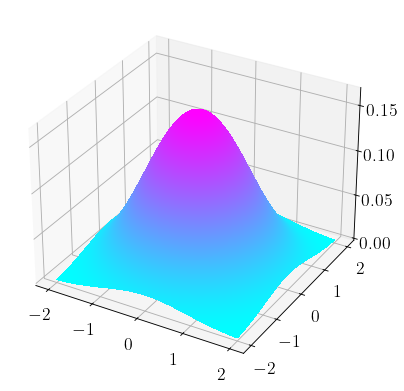}
\caption{$\varepsilon = 2^{-2}$}
\end{subfigure}
\hfill
\begin{subfigure}{0.48 \textwidth}
\centering
\includegraphics[width = \textwidth]{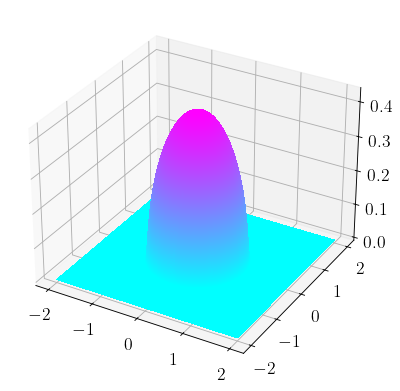}
\caption{$\varepsilon = 2^{-10}$}
\end{subfigure}
\caption{Numerical solutions to the fractional porous medium equation in self-similar variables \eqref{FullProbViscosity}. The simulations were performed for the fractional order $s=0.5$ with different values of the parameter $\varepsilon = 2^{-m}$, $m=2,\ldots, 10$, controlling the parabolic regularization. The computational domain is $\Omega = (-2,2)^2$, $\rho_0(x) \propto \frac{1}{|\Omega|}$; a uniform triangulation based on $2^9 \times 2^9$ equally 
 spaced nodes was used with $\Delta t = 0.01$.}
\label{3Dplots}
\end{figure}

In Figure \ref{FigsimulationsTime} and Figure \ref{3Dplots} we report the simulations performed for equation \eqref{FullProbViscosity} for different values of $\varepsilon$. 
We observe the effect of the parabolic regularization, which introduces Brownian diffusion of the density, and we also see that our numerical solution approaches the profile stated in \eqref{BarenblattProfile} as $\varepsilon$ decreases. 

\normalcolor
\subsection{Future work, challenges and open problems}
We close by highlighting several open problems that concern the implementation of the fully discrete scheme developed here. These, too, will be addressed in \cite{carrillosuli2024}. 
The presence of the upper cut-off parameter $L$ was crucial in the fixed point argument on which the proof,  in Lemma \ref{DiscExistLemma1},  of the existence of a solution to the fully-discrete scheme was based; and its presence was also essential in the derivation of the uniform bounds with respect to $\delta$ and the spatial mesh size stated in Lemma \ref{lem:uniformdeltah}.
While, as was noted above, the precise choice of the upper cut-off parameter $L$ becomes irrelevant in the convergence analysis (provided that $L$ is chosen so that it exceeds the supremum of the initial datum), this independence on $L$ was only shown to hold for the temporally semidiscrete scheme, i.e., \textit{after} passage to the limit with $\delta$ and the spatial discretization parameter has taken place; it is unclear whether a similar independence on $L$ of the solution to the fully discrete scheme holds true. This is a question that will need to be explored through numerical simulations. Similarly, the presence of the lower cut-off parameter $\delta$ in the fully-discrete scheme was essential in our proof of the nonnegativity of the solution of the temporally semidiscrete scheme (cf. the second paragraph of the proof of Theorem \ref{convergenceFEMTheo}), from which the nonnegativity of the limiting weak solution ultimately follows. The impact of the size of $\delta$ on the (sign of the) solution of the \textit{fully discrete scheme} is however another open question that will need to be assessed through numerical experiments. We do know however, thanks to Lemma \ref{lem:uniformdeltah}, that
$\underset{n = 1, \dots, N}{\max} \int_{\Omega} \pi_{h}( [\rhodLh^{n}]_{-}^{2}) \dx
\leq C(L)\delta$ for all $\delta \in (0,1)$.
\color{black}
In connection with this an interesting open problem is in fact whether the limits for the spatial discretization parameter $h \to 0$ and the time step $\Delta t \to 0$ commute. As was previously noted, from the theoretical point of view the dependence on $L$ in the energy estimate \eqref{APrioriBound2} in Lemma \ref{lem:uniformdeltah} prevents simultaneous passage to the limit in space and time. The proof of independence from $L$ of the semidiscrete-in-time solution, through Lemma \ref{InftyBoundLemmaDiscTime} and Lemma \ref{InfNormDecayDiscTime}, is grounded in the properties of the spectral fractional Laplacian $(-\Delta)^s$ and its representation formula. Having said this, we believe that these limits do commute. 

\begin{figure} 
\begin{subfigure}{0.48 \textwidth}
\centering
\includegraphics[width = \textwidth]{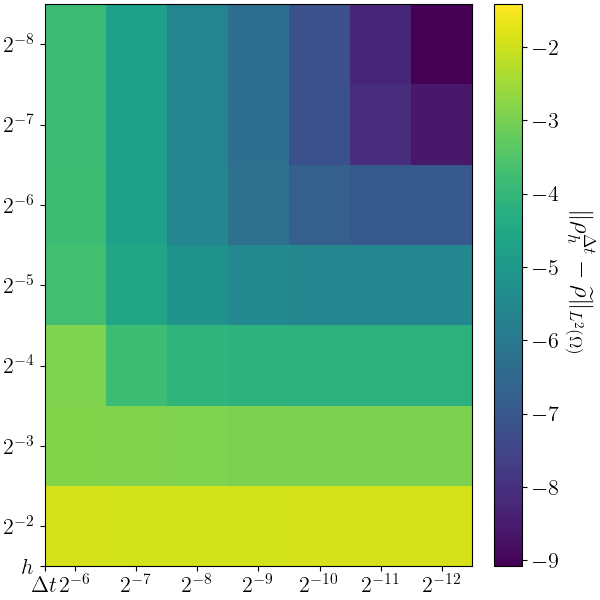}
\caption{$t=0.125$}
\end{subfigure}
\hfill
\begin{subfigure}{0.48 \textwidth}
\centering
\includegraphics[width = \textwidth]{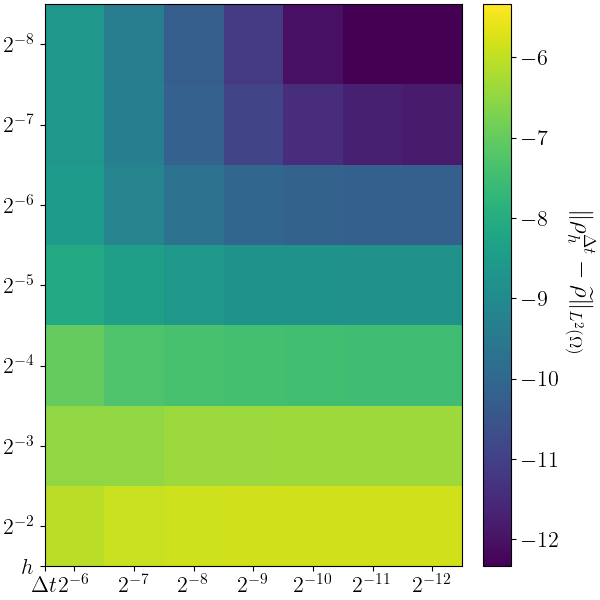}
\caption{$t=0.625$}
\end{subfigure}
\caption{Numerical solutions to the fractional porous medium equation \eqref{FullProbWSpace}, with passage to the limit for $h\to 0$ and $\Delta t \to 0$. The simulations were performed for the fractional order $s=0.5$ with different values of the spatial discretization parameter $h = 2^{-n}$, $n=2, \dots, 8$, and different time steps $\Delta t = 2^{-k}$, $k=6, \dots, 12$. The computational domain is $\Omega = (0,1)^2$, $\rho_0(x) \propto \exp(-|x|^2/(2 \pi \sigma))$. The parabolic regularization coefficient was $\varepsilon = 0.01$. }
\label{FigConvhDtHeatMap}
\end{figure}
To support this hypothesis, we conducted numerical experiments using the implementation described in \cite{carrillosuli2024}, and our results empirically confirm that the order of taking these limits does not affect the final numerical solution. More precisely, denoting by $\rho_h^{\Delta t}$ the finite element approximation of the solution to problem \eqref{MainFullProb} for fixed $h$ and $\Delta t$, we computed a series of approximations for $h=2^{-n}$, $n=2, \dots, 8$, and $\Delta t = 2^{-k}$, $k=6, \dots, 12$. We then evaluated the  $L^2(\Omega)$ norm of the difference between each approximation $\rho_h^{\Delta t}$ and a reference solution that we denote by $\widetilde\rho$, defined as the approximation corresponding to the smallest values of the discretization parameters $h$ and $\Delta t$, namely $\widetilde{\rho}=\rho_{h=2^{-8}}^{\Delta t = 2^{-12}}$.

\begin{figure} 
\begin{subfigure}{\textwidth}
\centering
\includegraphics[width = 0.7 \textwidth]{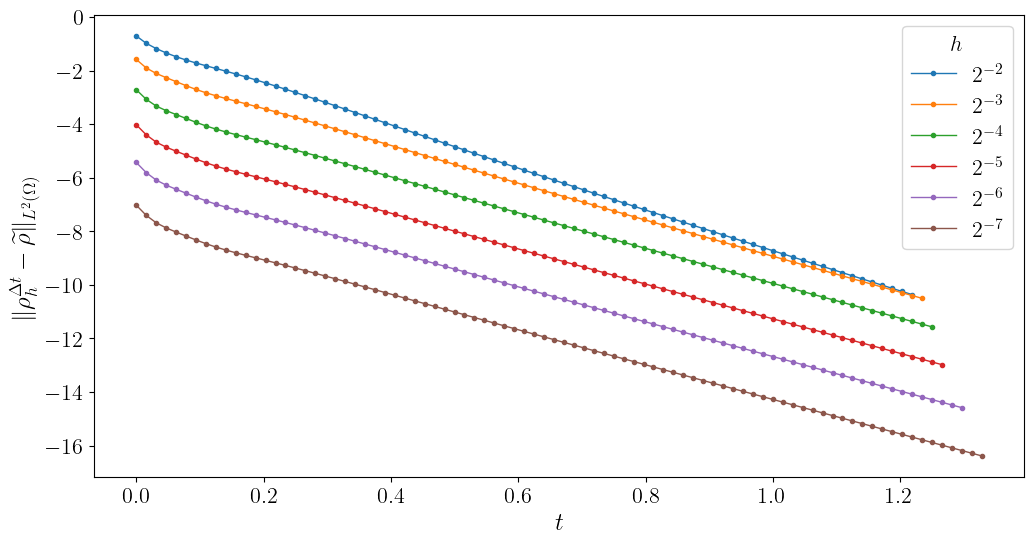}
\caption{$\Delta t = 2^{-12}$}
\end{subfigure}
\hfill
\begin{subfigure}{\textwidth}
\centering
\includegraphics[width = 0.7 \textwidth]{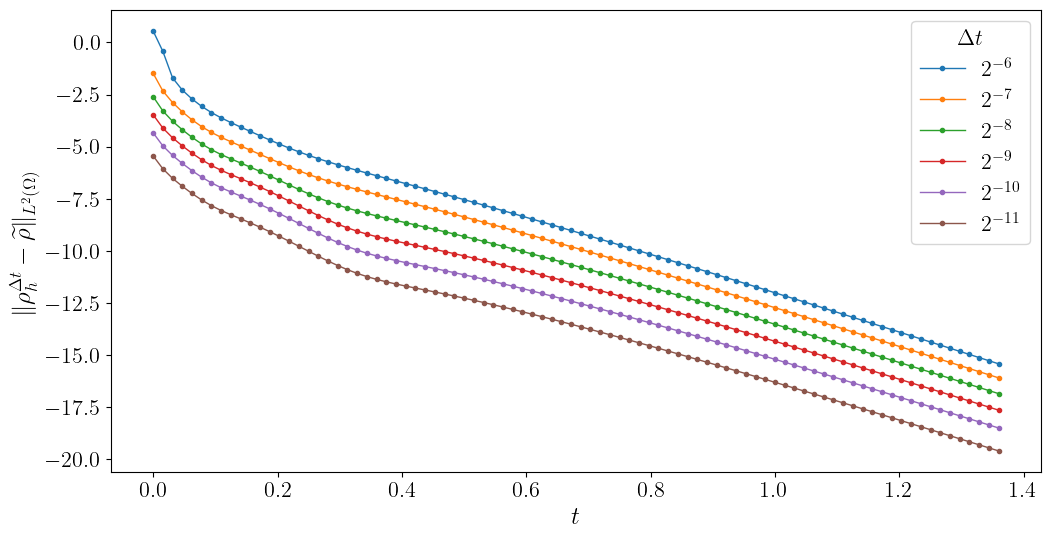}
\caption{$h = 2^{-8}$}
\end{subfigure}
\caption{Numerical solutions to the fractional porous medium equation \eqref{FullProbWSpace}, with passage to the limit for $h\to 0$ and $\Delta t \to 0$. The simulations were performed for the fractional order $s=0.5$ with different values of the spatial discretization parameter $h = 2^{-n}$, $n=2, \dots, 7$, keeping $\Delta t = 2^{-12}$ fixed, and for different values of time step $\Delta t = 2^{-k}$, $k=6, \dots, 11$, keeping $h=2^{-8}$ fixed. The computational domain is $\Omega = (0,1)^2$, $\rho_0(x) \propto \exp(-|x|^2/(2 \pi \sigma))$. The parabolic regularization parameter was $\varepsilon = 0.01$. }
\label{FigConvDtScatter}
\end{figure}

By examining the plots in Figures \ref{FigConvhDtHeatMap} and \ref{FigConvDtScatter},
we observe that $\| \rho_h^{\Delta t} - \widetilde{\rho}\|_{L^{2}(\Omega)} \to 0$ as $h \to 0$ and  $\Delta t \to 0$, regardless of whether the limits are taken simultaneously or sequentially (first in space then time, or vice versa). This supports our intuition that the limits commute in practice, even if the convergence analysis presented above was based on passing to the limit $h\to0$, before taking the limit $\Delta t \to 0$.

\begin{figure} 
\begin{subfigure}{0.48 \textwidth}
\centering
\includegraphics[width = \textwidth]{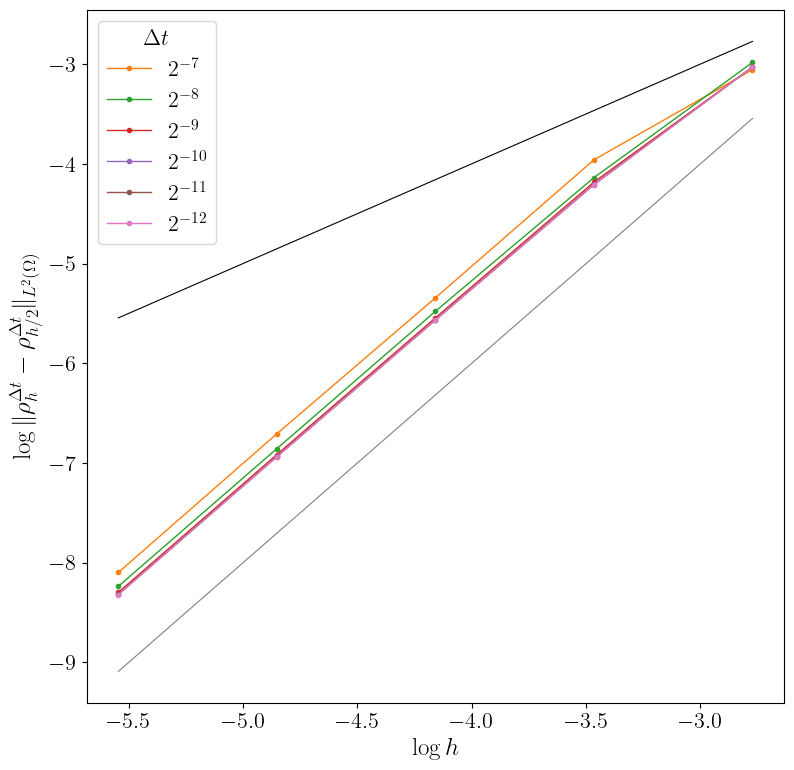}
\caption{$t=0.125$}
\end{subfigure}
\hfill
\begin{subfigure}{0.48 \textwidth}
\centering
\includegraphics[width = \textwidth]{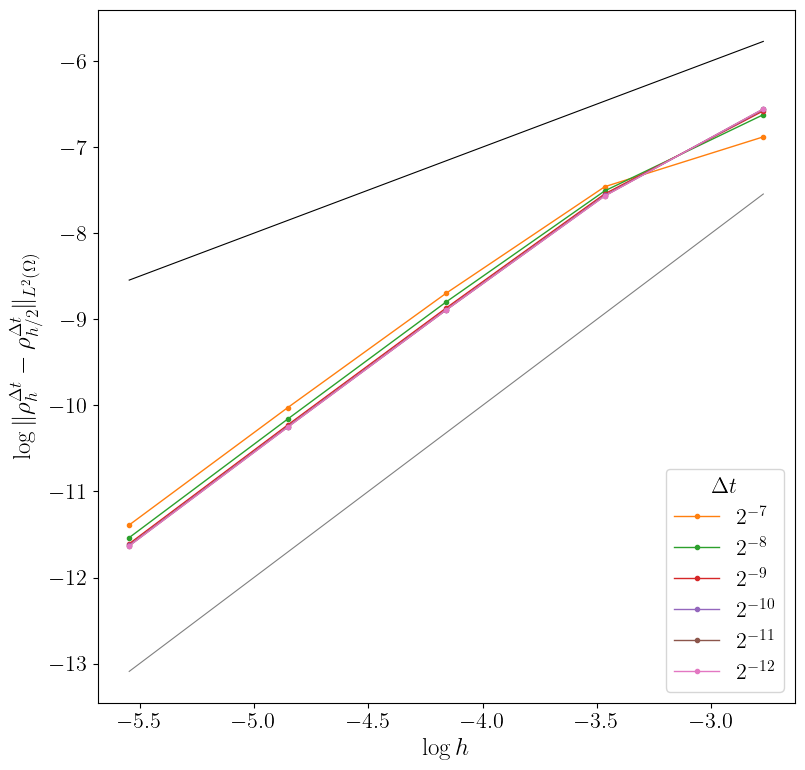}
\caption{$t=0.625$}
\end{subfigure}
\caption{Numerical solutions to the fractional porous medium equation \eqref{FullProbWSpace}, with convergence in space as $h\to 0$. The simulations were performed for the fractional order $s=0.5$, for different values of the spatial discretization parameter $h = 2^{-n}$, $n=2, \dots, 8$, and different time steps $\Delta t = 2^{-k}$, $k=6, \dots, 12$. The computational domain is $\Omega = (0,1)^2$, $\rho_0(x) \propto \exp(-|x|^2/(2 \pi \sigma))$. The parabolic regularization coefficient was $\varepsilon = 0.01$. }
\label{FigConvergence}
\end{figure}

An additional interesting question that arises from our analysis, and one that is commonly studied in the context of finite element methods, is the rate of spatial convergence as the spatial discretization parameter $h$ tends to zero. To this end,  we tested the method through a series of simulations.
Figure \ref{FigConvergence} shows that the observed rate of spatial convergence is quadratic, and that the choice of the time step does not influence the spatial speed of convergence.

\normalcolor
A further question that will require computational assessment, because of the presence of the nonlinear function $\Theta^L_\delta$ in the numerical method (cf. eq. \eqref{FullDiscWeakForm}), which is only known to be Lipschitz continuous (cf. Lemma \ref{ThetaContLemma}), is the choice of a suitable iterative method for the solution of the finite-dimensional system of nonlinear equations at a given time-level (e.g. fixed point iteration or a semismooth Newton method). 
Finally, we note that because the initial-boundary-value problem under consideration is not known to possess a \textit{unique} weak solution, the derivation of a bound on the error between a numerical solution defined by the scheme and the weak solution, to which a (sub)sequence of solutions generated by the scheme converges, is presently out of reach. This, too, is a question that will need to be explored through numerical simulations. 
\section{Conclusions}
\label{Sect6}
We have developed a fully discrete numerical method for the approximate solution of the fractional porous medium equation that involves the spectral fractional Neumann Laplacian on a bounded polytopal Lipschitz domain in two or three space dimensions, with the aim to model nonlocal effects in fluid flow through porous media, which are not captured by classical porous media flow models. The numerical method is based on a continuous piecewise linear finite element approximation on a weakly acute symplicial subdivision of the spatial domain and implicit Euler time-stepping. Our discretization of the nonlinear term appearing in the equation involves two positive cut-off parameters: $\delta \in (0,1)$, whose role is to ensure nonnegativity of the numerical solution, uniformly in the spatial and temporal discretization parameters; and $L>1$, whose role is to ensure that the sequence of numerical approximations is uniformly bounded with respect to the discretization parameters. By using a compactness argument, we have proved convergence of the scheme to a nonnegative and bounded  weak solution of the initial-boundary-value problem as the spatial discretization parameter, the lower cut-off parameter $\delta$, and the time-step all tend to zero. It transpired from the analysis that, once the upper cut-off parameter $L$ is chosen so that it exceeds the essential supremum of the (nonnegative) initial datum, its precise value is of no relevance when passing to the limit with $\delta$ and the spatial and temporal discretization parameters; conveniently, this makes passage to the limit $L \rightarrow +\infty$ redundant. We have also proved that the weak solution, to which a (sub)sequence of approximations generated by the scheme converges, satisfies an energy inequality and that the energy associated with that weak solution exhibits exponential decay in time. While the exponential decay results \eqref{eq:decay-a} and \eqref{eq:decay} did not require any additional restrictions on $\Omega$ or the range of the fractional order $s$, in order to prove the exponential decay \eqref{eq:decay-c} of the energy $E(\cdot)$ we had to assume convexity of $\Omega$ and that $s \in (1/2,1)$, so as to be able to apply Lemma \ref{W1InftyLemma}; whether or not these additional restrictions on $\Omega$ and $s$ can be relaxed or removed is an open problem. 
\section*{Acknowledgements}
This work was supported by the Advanced Grant Nonlocal-CPD (Nonlocal PDEs for Complex Particle Dynamics: Phase Transitions, Patterns and Synchronization) of the European Research Council Executive Agency (ERC) under the European Union’s Horizon 2020 research and innovation programme (grant agreement No. 883363). JAC was also partially supported by the EPSRC grant numbers EP/T022132/1 and EP/V051121/1.

\FloatBarrier

\bibliographystyle{abbrv}
\bibliography{./fractional-porous-FE.bib}
 
\renewcommand{\appendixname}{Appendix: The spectral fractional Laplacian}
\addappendix

We report here certain technical results related to the spectral fractional Neumann Laplacian (\ref{FracLapNeu}), which are of relevance for our analysis. 

For a bounded open Lipschitz domain $\Omega \subset \mathbb{R}^{d}$, $s \in (0,1)$, and $f \in \mathbb{H}_{\ast}^{-\s}(\Omega)$, consider the fractional Poisson equation involving the Neumann Laplacian:
\begin{equation} \label{FracPoiNeu}
(-\Delta_{\mathrm{N}})^{\s} u = f  \quad \textrm{in } \Omega.
\end{equation}
By a \textit{weak solution} to \eqref{FracPoiNeu} we mean a function $u \in \mathbb{H}^s_\ast(\Omega)$, $s \in (0,1)$, such that 
\[ \int_\Omega ((-\Delta_\mathrm{N} )^{s/2}u(x)) \, ((-\Delta_\mathrm{N} )^{s/2}v(x)) \dx = \langle f , v \rangle \quad \forall v \in \mathbb{H}^s_\ast(\Omega),\]
where $\langle \cdot, \cdot \rangle$ denotes the duality pairing between $\mathbb{H}^{-s}_\ast(\Omega)$ and $\mathbb{H}^s_\ast(\Omega)$, with $\mathbb{H}^{s}_\ast(\Omega) \hookrightarrow \mathbb{H}^{0}_\ast(\Omega)\equiv \mathbb{H}^{0}_\ast(\Omega)' \hookrightarrow\mathbb{H}^{-s}_\ast(\Omega)$ and $\mathbb{H}^{0}_\ast(\Omega)= L^2_\ast(\Omega)$. In particular if $f \in L^2_\ast(\Omega)$, then $\langle f, v \rangle = (f,v)$ for all $v \in \mathbb{H}^s_\ast(\Omega)$, $s \in [0,1)$.

\begin{theo}[Existence and uniqueness of a weak solution] \label{FracPoiExUniq}
Suppose that $\Omega\subset \mathbb{R}^d$ is a bounded open Lipschitz domain, $s \in (0,1)$, and $f \in \mathbb{H}^{-\s}_{\ast}(\Omega)$. Then, the boundary-value problem \eqref{FracPoiNeu} has a unique weak solution $u \in \Hsast(\Omega)$, and 
\begin{equation} \label{StabFracNeu}
\|u\|_{\Hs_{\ast}(\Omega)} \leq \|f\|_{\mathbb{H}^{-\s}_{\ast}(\Omega)}.
\end{equation}
\begin{proof}
The existence of a unique weak solution $u \in \mathbb{H}^s_\ast(\Omega)$ is a direct consequence of the Lax--Milgram lemma applied to the variational problem: find $u \in \mathbb{H}^s_\ast(\Omega)$ such that $a(u,v) = \ell(v)$ for all $v \in \mathbb{H}^s_\ast(\Omega)$, with the bilinear form 
\[ a(w,v):= \int_\Omega ((-\Delta _\mathrm{N} )^{s/2}w(x)) \, ((-\Delta_\mathrm{N} )^{s/2}v(x)) \dx,  \quad w, v \in \mathbb{H}^s_\ast(\Omega), \]
and the linear functional $\ell(v):= \langle f , v \rangle$, $v \in \mathbb{H}^s_\ast(\Omega)$. The stated stability inequality follows by noting that 
\[ \|u\|^2_{\mathbb{H}^s_\ast(\Omega)} = a(u,u) = \langle f , u \rangle \leq \|f\|_{\mathbb{H}^{-s}_\ast(\Omega)} \|u \|_{\mathbb{H}^s_\ast(\Omega)}.\]
That completes the proof of the lemma.  
\end{proof}
\end{theo}

We note in passing that the following function space interpolation inequality holds. 

\begin{lemma} \label{InterpFracSobLemma}
Let $\Omega\subset \mathbb{R}^d$ be a bounded open Lipschitz domain, $s \in (0,1)$, and $u \in \mathbb{H}^1_{\ast}(\Omega)$. Then,
\begin{equation*}
    \|u\|_{\Hs_{\ast}(\Omega)} \leq \|u\|_{L^{2}(\Omega)}^{1-\s} \; \| u \|_{\mathbb{H}^1(\Omega)}^{\s}.
\end{equation*}
\begin{proof}
The assertion of the lemma follows by applying H\"{o}lder's inequality. Indeed, we have for $\theta \in (0,1)$ that
\begin{align*}
    \|u\|_{\Hs_{\ast}(\Omega)}^{2} &= \sum_{k=1}^{\infty} \lambda_{k}^{\s} u_{k}^{2} = \sum_{k=1}^{\infty} |u_{k}|^{\theta} \lambda_{k}^{\s} |u_{k}|^{2 - \theta} \leq \Bigg( \sum_{k=1}^{\infty} |u_{k}|^{\theta p}\Bigg)^{\frac{1}{p}} \Bigg( \sum_{k=1}^{\infty} \lambda_{k}^{\s q} |u_{k}|^{(2-\theta)q} \Bigg)^{\frac{1}{q}},  
\end{align*}
for $p, q \in (1, \infty)$ such that $\frac{1}{p} + \frac{1}{q} = 1$. The stated inequality then follows by taking $p = \frac{2}{\theta}$, $q = \frac{2}{2 - \theta}$ and $\theta = 2(1 - \s)$, whereby $\theta p = 2$, $sq=1$, $(2-\theta)q=2$, $1/p=\theta/2=1-s$, and $1/q=s$.
\end{proof}
\end{lemma}

The second result that we rely on throughout our analysis is the following $H^{1}$ norm stability of the fractional Poisson equation 
involving the Neumann Laplacian on a Lipschitz domain.

\begin{lemma} \label{H1stablemma2}
    Suppose that $\Omega\subset \mathbb{R}^d$ is a bounded open Lipschitz domain, $s \in (0,1)$, and $f \in H^{1}_{\ast}(\Omega)$. Suppose further that $u \in \mathbb{H}^s_\ast(\Omega)$ is the unique weak solution of \eqref{FracPoiNeu}. Then, $u \in H^1_\ast(\Omega)$
    and there exists a positive constant $C$ such that 
    \begin{equation} \label{H1Stab}
        \| \nabla u \|_{L^{2}(\Omega)} \leq C\lambda_1^{-s} \| \nabla f \|_{L^{2}(\Omega)},
    \end{equation}
    where $\lambda_1$ is the smallest positive eigenvalue of the Neumann Laplacian on $\Omega$.  Moreover, if $s \in [1/2, 1)$, then also 
    \begin{equation} \label{H1L2Stab}
        \| \nabla u \|_{L^{2}(\Omega)} \leq C\lambda_1^{1-2s} \| f \|_{L^{2}(\Omega)},
    \end{equation}
    for a positive constant $C>0$.
\begin{proof} Let $\psi_k$, $k=1,2, \ldots$, denote the eigenfunctions of the Neumann Laplacian corresponding to the positive eigenvalues $\lambda_k$, $k=1,2, \ldots$, where $\lambda_1$ is the smallest positive eigenvalue. As $\{ \psi_k \}_{k \geq 1}$ is a complete orthonormal system in $L^2_\ast(\Omega)$ and $f \in H^1_\ast(\Omega) \subset L^2_\ast(\Omega)$, it follows that 
\[ f(x) = \sum_{k=1}^\infty f_k \psi_k(x), \quad \mbox{with } f_k := \int_\Omega f(x) \psi_k(x) \dx,\]
with the convergence of the series understood to be in the norm of $L^2(\Omega)$; i.e., the $M$-th partial sum
\[ f^M(x) := \sum_{k=1}^M f_k \psi_k(x), \quad x \in \Omega, \]
converges to $f$ strongly in the norm of $L^2(\Omega)$. Similarly, because $u \in \mathbb{H}^s_\ast(\Omega) \subset L^2_\ast(\Omega)$, we have that
\[u(x) = \sum_{k=1}^\infty u_k \psi_k(x), \quad \mbox{with } u_k := \int_\Omega u(x) \psi_k(x) \dx,\]
and the $M$-th partial sum
\[ u^M(x) := \sum_{k=1}^M u_k \psi_k(x), \quad x \in \Omega,\]
converges to $u$ strongly in the norm of $L^2(\Omega)$. Because $u$ is the unique weak solution of \eqref{FracPoiNeu}, it follows by taking $v=\psi_k$ as test function in the weak formulation that $u_k:=\lambda_k^{-s} f_k$, $k=1,2, \ldots$. 
    Furthermore,
    \begin{align}
    \|\nabla u^M\|^2_{L^2(\Omega)}= (\nabla u^M, \nabla u^M) &= \bigg(\sum_{k=1}^M u_k \nabla \psi_k , \sum_{\ell=1}^M u_\ell \nabla \psi_\ell\bigg)  = \sum_{k=1}^M \sum_{\ell=1}^M u_k u_{\ell} (\nabla \psi_k, \nabla \psi_l)\nonumber\\
    &= \sum_{k=1}^M \sum_{\ell=1}^M u_k u_{\ell} (-\Delta \psi_k, \psi_\ell) = \sum_{k=1}^M \sum_{\ell=1}^M \lambda_k u_k u_{\ell} (\psi_k, \psi_\ell) = \sum_{k=1}^M \sum_{\ell=1}^M \lambda_k u_k u_\ell \delta_{k,\ell} \nonumber \\
    &= \sum_{k=1}^M \lambda_k u_k^2 = \sum_{k=1}^M (\lambda_k)^{1-2s} f_k^2 \leq \lambda_1^{-s}  \sum_{k=1}^M (\lambda_k)^{1-s} f_k^2 \label{RefH1L2Stab1} \\
    &\leq \lambda_1^{-s}  \sum_{k=1}^\infty (\lambda_k)^{1-s} f_k^2 = \lambda_1^{-s} \|f\|^2_{\mathbb{H}^{1-s}(\Omega)}\quad \forall M \geq 1. \label{ineq-uM}
    \end{align}
As was noted in Remark \ref{RemFracSobs}, thanks to Lemma 7.1 in \cite{caffarelli2016fractional}, for $s \in (0,1)$ a function $f$ belongs to $\mathbb{H}^{1-s}_\ast(\Omega)$ if, and only if, $f \in H^{1-s}(\Omega)$ and $\int_{\Omega} f(x) \dx = 0$; in addition, the norm $\|\cdot\|_{\mathbb{H}^{1-s}_\ast(\Omega)}$ is equivalent to the seminorm $|\cdot|_{H^{1-s}(\Omega)}$, which is a norm on $\mathbb{H}^{1-s}_\ast(\Omega)$. There is therefore a positive constant $C$, independent of $f$, such that $\|f\|^2_{\mathbb{H}^{1-s}(\Omega)} \leq C |f|^2_{H^{1-s}(\Omega)}$, and because $H^1(\Omega) \hookrightarrow H^{1-s}(\Omega)$, also $\|f\|^2_{\mathbb{H}^{1-s}(\Omega)}\leq C \|f\|^2_{H^1(\Omega)}$. As $f$ has been assumed to belong to $H^1_\ast(\Omega)$, whereby $f \in H^1(\Omega)$ and $\int_\Omega f(x) \dx = 0$, Poincar\'e's inequality implies the existence of a positive constant $C$ such that 
$ \|f\|^2_{H^1(\Omega)} \leq C \|\nabla f\|^2_{L^2(\Omega)}$. Thus we have shown the existence of a positive constant $C$ (independent of $M$) such that 
\begin{equation}\label{ineq-uMa} \|\nabla u^M\|^2_{L^2(\Omega)} \leq C\lambda_1^{-s} \|\nabla f\|^2_{L^2(\Omega)}\quad \forall M \geq 1.
\end{equation}

Using an analogous calculation to the one in \eqref{ineq-uM}, for any positive integers $M$ and $N$ such that $M>N$ we have that
\begin{align*}
     \|\nabla u^M- \nabla u^N\|^2_{L^2(\Omega)}  \leq \lambda_1^{-s}
\sum_{k=N+1}^M (\lambda_k)^{1-s} f_k^2. 
\end{align*}
Because the infinite series appearing in the last line of 
\eqref{ineq-uM} is convergent, it then follows that for any $\varepsilon>0$ there exist positive integers $M$ and $N$, with $M>N$, such that 
\begin{align*}
    \|\nabla u^M- \nabla u^N\|^2_{L^2(\Omega)} \leq \lambda_1^{-s}
\sum_{k=N+1}^M (\lambda_k)^{1-s} f_k^2 < \varepsilon,
\end{align*}
i.e., $(\nabla u^M)_{M \geq 1}$ is a Cauchy sequence in the Hilbert space
$L^2(\Omega;\mathbb{R}^d)$; therefore $\nabla u^M$ converges (strongly in the norm of $L^2(\Omega;\mathbb{R}^d)$ to an element of $L^2(\Omega;\mathbb{R}^d)$. As $u^M$ converges to $u \in L^2(\Omega)$ strongly in the norm of $L^2(\Omega)$, by the uniqueness of the limit, $\nabla u^M$ must converge to $\nabla u$ in $L^2(\Omega;\mathbb{R}^d)$. Therefore,
\[ \lim_{M \rightarrow \infty} \|\nabla u^M\|^2_{L^2(\Omega)} = \|\nabla u\|^2_{L^2(\Omega)}.\]
With this information we now return to \eqref{ineq-uMa} and pass to the limit $M \rightarrow \infty$ on the left-hand side of the inequality to deduce the desired inequality \eqref{H1Stab}. 

Inequality \eqref{H1L2Stab} can be proved by following the same argument, but also using in \eqref{RefH1L2Stab1} and the assumed constraint that $s \in [1/2, 1)$, which implies that $1-2s \leq 0$; therefore,
\begin{equation*}
    \| \nabla u^{M} \|_{L^{2}(\Omega)}^{2} \leq C \lambda_{1}^{1-2s} \| f \|_{L^{2}(\Omega)}^{2} \quad \forall M \geq 1.
\end{equation*}
The desired inequality then follows by passing to the limit $M \rightarrow \infty$. This completes the proof.
    \end{proof}
\end{lemma}

In our context another useful result is a type of Sobolev--Poincar\'{e} inequality. In \cite{hurri2013fractional} it was shown that, for a bounded open  John domain $\Omega \subset \mathbb{R}^d$, the following inequality holds for $s, \tau \in (0,1)$, $p< d/s$ and $1 < p \leq q \leq \frac{dp}{d - sp}$: 
\begin{equation} \label{FracPoincare}
    \| u - \bar{u} \|_{L^{q}(\Omega)} \leq C \left( \int_{\Omega} \int_{\Omega \cap B(x, \tau \,\mathrm{dist}(x,\partial\Omega))} \frac{|u(x) - u(y)|^{p} }{|x-y|^{d+sp}} \, \dy \,\dx \right)^{\frac{1}{p}}.
\end{equation}
The seminorm appearing on the right-hand side of (\ref{FracPoincare}) is weaker than the usual Gagliardo--Slobodetski\u{\i} seminorm of the fractional-order Sobolev space $W^{s, p}(\Omega)$,
\begin{equation*}
  \left( \int_{\Omega} \int_{\Omega} \frac{|u(x) - u(y)|^{p} }{|x-y|^{d+sp}} \, \dx \,\dy \right)^{\frac{1}{p}}.  
\end{equation*}
Indeed, the two seminorms induce different fractional order Sobolev spaces, which coincide in the case of a bounded open Lipschitz domain, as was shown in \cite{dyda2006comparability}. Thus, in the case of a Lipschitz domain $\Omega$, as is the case in the setting under consideration in this paper, and for $p=q=2$, we have the following fractional Poincar\'e inequality
for $s \in (0,1)$:
\begin{subequations}
\begin{equation} \label{FracPoincareUse0}
   \| u - \bar{u} \|_{L^{2}(\Omega)} \leq C |u|_{H^{s}(\Omega)}, 
\end{equation}
where $C$ is a positive constant, independent of $u$. Therefore, in particular, 
\begin{equation} \label{FracPoincareUse1}
   \| u\|_{L^{2}_\ast(\Omega)} \leq C |u|_{H^{s}_\ast(\Omega)}.
\end{equation}
An analogous inequality holds for functions $u \in \mathbb{H}^{s}_{\ast}(\Omega)$ and $s \in (0,1)$: one has that  
\begin{equation} \label{FracPoincareUse}
    \|u\|_{L^2_\ast(\Omega)} 
    \leq \lambda_1^{-s/2} \| u \|_{\mathbb{H}^{s}_{\ast}(\Omega)}.
\end{equation}
\end{subequations}
This inequality, in conjunction with Theorem \ref{FracPoiExUniq}, ensures the $L^{2}_\ast(\Omega)$ norm stability of the weak solution to the fractional Poisson equation based on the spectral definition of the Neumann Laplacian, for $f \in L^{2}_\ast(\Omega)$.

We now prove some auxiliary results, which concern the case when the right-hand side datum of equation \eqref{FracPoiNeu} has $L^{\infty}(\Omega)$ regularity or $W^{1, \infty}(\Omega)$ regularity. Let us denote respectively by $L^{\infty}_{\ast}(\Omega)$ and $W^{1, \infty}_{\ast}(\Omega)$ the spaces of functions in $L^{\infty}(\Omega)$ and $W^{1, \infty}(\Omega)$ with zero integral average over $\Omega$

\begin{lemma} \label{EqPVSpecLemma}
Suppose that $\Omega\subset \mathbb{R}^d$ is a bounded open Lipschitz domain, $s \in (0,1)$, and $f \in L^{\infty}_{\ast}(\Omega)$. Then, we have that
\begin{equation} \label{EqPVSpec}
(-\Delta_{\mathrm{N}})^{-s} f \simeq C_{\Omega, s, d} \int_{\Omega} \frac{1}{|x-y|^{d-2s}} f(y) \, \dy,
\end{equation}
for a positive constant $C_{\Omega, s, d}$ depending only on $\Omega, s$ and $d$.
\begin{proof}
Since $L^{\infty}_{\ast}(\Omega) \subset L^2_\ast(\Omega) \subset \mathbb{H}^{-s}_{\ast}(\Omega)$, we have that
\begin{equation} \label{Eq1}
(-\Delta_{\mathrm{N}})^{-s} f = \sum_{k =1}^{\infty} \lambda_{k}^{-s} f_{k} \psi_{k}, \quad \text{where } f_{k} = \int_{\Omega} f \psi_{k} \dx.
\end{equation}
By noting that for  $\lambda > 0$ and $s \in (0, 1)$,
\begin{equation*}
\lambda^{-s} = \frac{1}{\Gamma(s)} \int_{0}^{\infty} {\mathrm e}^{-t\lambda} \frac{1}{t^{1-s}} \, \dt, 
\end{equation*}
and inserting this in (\ref{Eq1}), we get that 
\begin{align}
(-\Delta_{\mathrm{N}})^{-s} f &= \frac{1}{\Gamma(s)} \int_{0}^{\infty} \sum_{k=1}^{\infty} {\mathrm e}^{-t\lambda_{k}} f_{k} \psi_{k} \frac{1}{t^{1-s}} \, \dt = \frac{1}{\Gamma(s)} \int_{0}^{\infty} \frac{{\mathrm e}^{t\Delta_{\mathcal{N}} f}}{t^{1-s}} \, \dt \nonumber \\
&= \frac{1}{\Gamma(s)} \int_{0}^{\infty} \int_{\Omega} W_{t}(x, y) \frac{1}{t^{1-s}} f(y) \, \dy \, \dt,\label{FracPoiHeatNeuKer}
\end{align}
where $\{ {\mathrm e}^{t\Delta_{\mathcal{N}}} \}_{t >0}$ stands for the heat semigroup associated to $\Delta_{\mathrm{N}}$, i.e. $e^{t\Delta_{\mathcal{N}}} u$ solves problem \eqref{HeatNeuSemiGroup} and $W_{t}(x, y)$  is the distributional heat kernel for $-\Delta$ subject to a homogeneous Neumann boundary condition on $\partial\Omega$.

The following bounds are known to hold for the distributional heat kernel (cf. Theorem 3.1 in \cite{SaloffCoste} or Section 7 in \cite{caffarelli2016fractional}):
\begin{equation*}
{c}_{1} \frac{e^{-|x-y|^{2}/(c_{3}t)}}{t^{\frac{d}{2}}} \leq W_{t}(x, y) \leq c_{2} \frac{{\mathrm e}^{-|x-y|^{2}/(c_{4}t)}}{t^{\frac{d}{2}}}, \quad x, y \in \Omega, \quad t > 0, 
\end{equation*}
for positive constants $c_{1}, c_{2}, c_{3}$ and $c_{4}$. Therefore, we have that 
\begin{align*}
\int_{0}^{\infty} \frac{1}{t^{1-s}} W_{t}(x, y) \, \dt &\simeq \int_{0}^{\infty} \frac{1}{t^{\frac{d}{2} + 1 - s}} {\mathrm e}^{-|x-y|^{2}/t} \, \dt \\
&=\frac{1}{|x-y|^{d-2s}}\int_{0}^{\infty} {\mathrm e}^{-r} r^{-1-s+\frac{d}{2}} \, \textrm{d}r = \Gamma\Big(\frac{d}{2} - s\Big) \frac{1}{|x-y|^{d-2s}}, 
\end{align*}
where in the passage from the first line to the second line we used the change of variable $r = |x-y|^{2}/t$.
This, together with expression \eqref{FracPoiHeatNeuKer}, gives us the desired result.
\end{proof}
\end{lemma}

\begin{cor} \label{LInftyPoiStabCor}
Suppose that $\Omega\subset \mathbb{R}^d$ is a bounded open Lipschitz domain, $s \in (0,1)$, and $f \in L^{\infty}_{\ast}(\Omega)$. Let $u \in \mathbb{H}^{s}_{\ast}(\Omega) $ be the unique solution of \eqref{FracPoiNeu}. Then $u \in L^{\infty}(\Omega)$ and
\begin{equation} \label{LInftyPoiStab}
    \| u \|_{L^{\infty}(\Omega)} \leq K_{\Omega, s, d} \| f \|_{L^{\infty}(\Omega)},
\end{equation}
for a positive constant $K_{\Omega, s, d}$, depending only on $\Omega, s$ and $d$.
\begin{proof}
    The proof is a direct consequence of \eqref{EqPVSpec} and the fact that the function $1/|z|^{d - 2s}$ is integrable for $s \in (0,1)$ over any bounded ball in $\mathbb{R}^d$. 
\end{proof}
\end{cor}

The following result is an extension to the fractional order setting of a result by Maz'ya \cite{mazya2009}, concerning the boundedness of the gradient of the solution to the Neumann Laplace problem in a convex domain.

\begin{lemma} \label{W1InftyLemma}
    Suppose that $\Omega\subset \mathbb{R}^d$ is a bounded open convex domain, $s \in (1/2,1 )$, and $f \in L^{ \infty}_{\ast}(\Omega)$.
     Let $u \in \mathbb{H}^{s}_{\ast}(\Omega) $ be the unique solution of \eqref{FracPoiNeu}. Then,  $u \in W^{1, \infty}_{\ast}(\Omega)$ and 
     \begin{equation*}
        \| u \|_{W^{1, \infty}(\Omega)} \leq C \| f \|_{L^{\infty}(\Omega)},
     \end{equation*}
     for a positive constant $C>0$.
     \begin{proof}
    By Corollary \ref{LInftyPoiStabCor} we have in first place that $u \in L^{\infty}(\Omega)$, with the $L^{\infty}(\Omega)$ stability bound \eqref{LInftyPoiStab}. 

    Let us now prove that $u \in W^{1, \infty}(\Omega)$. We start from expression \eqref{FracPoiHeatNeuKer} and we use the following inequality for the gradient of the heat kernel:
    \begin{equation*}
        | \nabla_{x} W_{t}(x, y) | \leq C_{1} \frac{1}{t^{\frac{d+1}{2}}} {\mathrm e}^{-|x-y|^{2}/(C_2 t)},\quad x, y \in \Omega, \quad t > 0;
    \end{equation*}
    this inequality is known to hold in bounded open convex domains (cf.  Lemma 3.1 in \cite{Wang2013} or Theorem 1.2 in \cite{yan2010gradient}) for positive constants $C_{1}$ and $C_{2}$.
    We then have the following inequality for $x \in \Omega$ and for some positive constant $C>0$:
    \begin{equation*}
        |\nabla u (x)| \leq C \frac{1}{\Gamma(s)} \int_{0}^{\infty} \int_{\Omega} \frac{1}{t^{\frac{d}{2} + \frac{3}{2} - s}} {\mathrm e} ^{-|x-y|^{2}/t} |f(y)| \dy \dt. 
    \end{equation*}
By using the change of variable $r=|x-y|^2/t$ in the above equality together with the fact that $f \in L^{\infty}(\Omega)$, we then have that
    \begin{align*}
    |\nabla u(x) | &\leq C \frac{\| f \|_{L^{\infty}(\Omega)}}{\Gamma(s)} \int_{\Omega} \frac{1}{|x-y|^{d-2s+1}} \int_{0}^{\infty} {\mathrm e}
^{-r} r^{-\frac{1}{2} - s + \frac{d}{2}} \text{d}r \dy \\
&= C \frac{\| f \|_{L^\infty(\Omega)}}{\Gamma(s)}\, \Gamma \bigg( \frac{d}{2}+ \frac{1}{2} - s\bigg) \int_\Omega \frac{1}{|x-y|^{d-2s+1}}  \dy.
    \end{align*}
The final assertion of the lemma then follows by noting that the function $1/|z|^{d-2s+1} $ is integrable for $s\in (1/2, 1)$ over any bounded ball in $\mathbb{R}^d$. 
     \end{proof}
\end{lemma}
\color{black}

\label{sec:appendix}

\end{document}